\documentclass[leqno,a4paper]{amsart}
\usepackage{egastyle}
\usepackage{mymacros}
\usepackage[all,tips]{xy}
\usepackage{url}
\usepackage{bm}
\usepackage{verbatim}
\usepackage{psfrag}
\usepackage{graphicx}
\usepackage{caption}

\CompileMatrices

\def\B{\bB}

\def\G{\bG}

\DeclareMathOperator\trop{trop}

\DeclareMathOperator\slope{slope}
\renewcommand\red{\operatorname{red}}

\newcommand\an{{\operatorname{an}}}
\renewcommand\div{{\operatorname{div}}}

\newcommand\res{{\operatorname{res}}}

\newcommand\can{{\operatorname{can}}}
\renewcommand\cc{{\circ\circ}}

\renewcommand{\Div}{{\rm div}}
\newcommand{\cyc}{{\rm cyc}}
\newcommand{\ve}{{\varepsilon}}

\newcommand{\Acal}{{\mathscr A}}

\newcommand{\Ccal}{{\mathscr C}}

\newcommand{\Fcal}{{\mathscr F}}

\newcommand{\Mcal}{{\mathscr M}}

\newcommand{\Ocal}{{\mathscr O}}

\newcommand{\Scal}{{\mathscr S}}

\newcommand{\Ucal}{{\mathscr U}}

\newcommand{\Xcal}{{\mathscr X}}

\newcommand{\Ycal}{{\mathscr Y}}

\newcommand{\Zcal}{{\mathscr Z}}

\newcommand{\rdop}{{\mathbf R}}

\newcommand{\bdop}{{\mathbf B}}

\newcommand{\zdop}{{\mathbf Z}}

\renewcommand\proof\pf

\newcommand{\Val}{{\rm Val}}
\newcommand{\Trop}{{\rm Trop}}

\newcommand{\mb}{{\mathbf m}}

\newcommand{\Deltabar}{{\overline{\Delta}}}
\newcommand{\Ccalbar}{{\overline{\Ccal}}}

\newcommand{\Xan}{{X^{\rm an}}}

\newcommand{\Tan}{{T^{\rm an}}}
\newcommand{\Uan}{{U^{\rm an}}}

\newcommand{\relint}{{\rm relint}}

\newcommand{\str}{{\rm str}}
\newcommand{\kcirc}{K^{\circ}}
\newcommand{\ktilde}{\widetilde{K}}

\newcommand{\In}{{\rm in}}

\newcommand{\tropbar}{{\overline{\trop}}}
\newcommand{\inn}{\operatorname{in}}
\newcommand{\aff}{{\operatorname{aff}}}

\newcommand{\Gr}{\mathrm{Gr}}

\makeatletter
\@namedef{subjclassname@2010}{%
 \textup{2010} Mathematics Subject Classification}
\makeatother

\title{Skeletons and tropicalizations}

\author{Walter Gubler}
\address{Walter Gubler, Fakult\"at f\"ur Mathematik,  Universit\"at Regensburg,
Universit\"atsstra{\ss}e 31, D-93040 Regensburg}
\email{walter.gubler@mathematik.uni-regensburg.de}

\author{Joseph Rabinoff}
\address{Joseph Rabinoff, School of Mathematics, Georgia Institute of Technology, Atlanta GA 30332-0160, USA}
\email{jrabinoff@math.gatech.edu}

\author{Annette Werner}
\address{Annette Werner, Institut f\"ur Mathematik, Goethe-Universit\"at
  Frankfurt, Robert-Mayer-Stra{\ss}e 8, D-60325 Frankfurt a.M.}
\email{werner@math.uni-frankfurt.de}

\begin{document}

\begin{abstract}
  Let $K$ be a complete, algebraically closed non-archimedean field with ring of
  integers $K^\circ$ and let $X$ be a $K$-variety. 
  We associate to the data of a strictly semistable $K^\circ$-model $\sX$ of $X$ plus a suitable
  horizontal divisor $H$ a skeleton $S(\sX,H)$ in the analytification of 
  $X$. This generalizes Berkovich's original construction by
  admitting unbounded faces in the directions of the components of $H$.
  It also generalizes constructions by Tyomkin and Baker--Payne--Rabinoff from
  curves to higher dimensions. Every such skeleton has an integral 
  polyhedral structure. We show that the valuation of a non-zero rational
  function is piecewise linear on $S(\sX, H)$. For such functions we define
  slopes along codimension one faces and prove a slope formula expressing a
  balancing condition on the skeleton. Moreover, we obtain a multiplicity
  formula for skeletons and tropicalizations in the spirit of a well-known result by
  Sturmfels--Tevelev. 
  We show a faithful tropicalization result saying roughly
  that every skeleton can be seen in a suitable tropicalization. We also prove a
  general result about existence and uniqueness of a continuous section to the
  tropicalization map on the locus of tropical multiplicity one. 
\end{abstract}

\maketitle

\section{Introduction}

\paragraph
Throughout this paper, $K$ denotes an algebraically closed 
field which is complete with respect
to a non-trivial, non-archimedean valuation $v: K\to\R\cup\{\infty\}$. The
corresponding valuation ring is denoted by $\kcirc$ and the value group by 
$\Gamma \coloneq v(K^\times)\subset\R$.

\paragraph[Tropicalizations]
Let $X$ be a  $K$-variety, i.e.\ an integral, separated $K$-scheme of finite
type.  Suppose that $\phi: X\to T = \Spec(K[x_1^{\pm 1},\ldots,x_n^{\pm 1}])$ 
is a closed immersion of  $X$ into a
multiplicative torus.  To 
$\phi$ we associate the \emph{tropicalization} $\Trop(X)$ of $X$.  As a
set, $\Trop(X) = \trop\circ\phi(X^\an)$, where
$\trop:T^\an\to\R^n$ is the valuation map given by
\[ \trop(p) = \big(-\log|x_1(p)|,\ldots,-\log|x_n(p)|\big), \]
and $(\scdot)^\an$ denotes analytification in the sense of Berkovich.  
By the Bieri--Groves theorem and work of Speyer--Sturmfels, the tropicalization $\Trop(X)$ can be enriched with the structure of a
balanced, weighted, integral polyhedral complex of pure dimension 
$d = \dim(X)$ (see ~\parref{par:tropicalization} for details).
Tropicalizations have proven to be interesting objects to study: on the
one hand they are combinatorial in nature, and as such are amenable to explicit
calculations; on the other hand, they are rich enough objects to be used
as a tool to study the original variety $X$.  
An excellent introduction to the subject can be found in
the book of Maclagan and Sturmfels~\cite{maclagan_sturmfels:book}.

As $\Trop(X)$ depends on the
embedding $\phi$, for our purposes we will sometimes call $\Trop(X)$ an
\emph{embedded} or \emph{parameterized tropicalization} of $X$.

\paragraph[Skeletons] 
Now suppose that $X$ is a proper, smooth $K$-variety with a strictly semistable
$\kcirc$-model $\sX$.  This is a proper, flat scheme over $\kcirc$ with generic
fibre $X$ such that the special fibre is a simple normal crossing divisor (see
Definition \ref{strictly semistable pairs}). Berkovich introduces the skeleton
$S(\Xcal)$ of $\Xcal$ as a closed subset of $\Xan$ in
\cite{berkovich:locallycontractible1}. He shows that $S(\Xcal)$ is a piecewise
linear space of dimension bounded by $\dim(X)$ which is covered by canonical
simplices reflecting the stratification of the special fibre $\sX_s$. In
particular, the vertices are in bijective correspondence with the irreducible
components of $\sX_s$.  The skeleton is in a canonical way a proper strong
deformation retraction of $\Xan$. For details and generalizations to the
analytic setting and to pluristable models, we refer
to~\cite{berkovich:locallycontractible1,berkovich:locallycontractible2}. 

The piecewise linear structure carried by $S(\sX)$ is strongly analogous to that of a
tropicalization $\Trop(X)$.  For the purposes of the introduction, we will
regard $S(\sX)$ as an \emph{intrinsic tropicalization} of the variety $X$.
 
\paragraph[The case of curves]
In the special case of a smooth projective curve $X$
the skeleton $S(\sX)$ is a metric graph whose underlying graph is the incidence
graph of the special fibre $\sX_s$.  
In this case, both the skeleton and tropicalizations are metrized graphs. 
However, the skeleton is bounded while $\Trop(X)$ is unbounded, which makes
direct comparisons between the two awkward.
To remedy this, Tyomkin \cite{tyomkin} introduces the skeleton $S(\sX,H)$ of a
marked curve $(X,H)$, by adding a ray in the direction of every marked point to
$S(\sX)$. 
With this tool, Tyomkin obtains 
an algebraic proof and a generalization of Mikhalkin's correspondence theorem. 
The latter is the key in Mikhalkin's pioneering work on Gromov--Witten invariants on 
the plane. Mikhalkin's  original proof in~\cite{Mikhalkin:enumerative}
is based on complex analytic and symplectic techniques which are fundamentally different 
from the non-archimedean techniques used in our paper. 

The program of comparing skeletons and tropicalizations
was launched in a systematic way by Baker, Payne and Rabinoff 
in~\cite{bpr:trop_curves,bpr:analytic_curves}
who work more generally over an algebraically closed field $K$ with a
non-trivial, non-archimedean, complete, real valuation.  (Tyomkin considers a
complete discretely valued ground field; one recovers his construction from the one 
in~\cite{bpr:trop_curves,bpr:analytic_curves} after base extension.)  For a
projective smooth curve $X$ with a set $H$ of marked points reducing to
distinct smooth points in the special fibre of a given semistable model, the
skeleton $S(\sX,H)$ is realized in~\cite{bpr:analytic_curves} as
a metrized unbounded graph which is a subset of $\Xan$ and is the target of a canonical
retraction map $\tau:\Xan \setminus H \rightarrow S(\sX,H)$. If $X$ is embedded
in a toric variety with dense torus $T$ such that $X\setminus H\subset X\cap T$,
then comparison theorems relating the skeleton $S(\sX,H)$ and the
tropicalization $\Trop(X \cap T)$ are proved in~\cite{bpr:trop_curves}.

\paragraph[Goals]
The overall goal of this paper is a careful study of the relationship
between intrinsic and  parameterized tropicalizations of a variety
$X$.  We generalize a substantial part of the results in~\cite{bpr:trop_curves,bpr:analytic_curves} to higher dimensions. 
We hope that they will be useful for correspondence theorems in higher 
dimensions and for applications to arithmetic geometry as for example for the development of a non-archimedean Arakelov theory.

Let us now describe our main results. 

\paragraph[Skeletons for strictly semistable pairs]
Suppose that $X$ is a proper, smooth $K$-variety.  A 
\emph{strictly semistable pair} roughly consists of a strictly semistable
$\kcirc$-model $\sX$ of $X$ along with a Cartier divisor $H$ on $X$ such that
$H$ plus the special fibre $\sX_s$ of $\sX$ is a simple normal crossings divisor.
See~\S\ref{Strictly semistable pairs} for a precise definition.  To such a
pair we associate a \emph{skeleton} $S(\sX,H)$ of $X$.  The skeleton is a
closed subset of the analytification of $U \coloneq X\setminus\Supp(H)$ and its dimension
is bounded above by $d = \dim(X)$.  It turns out that $S(\sX,H)$ is a piecewise linear space whose combinatorics reflect the
stratification of $\sX_s$ associated to $D = H + \sX_s$.  This means more precisely that for every stratum $S \subset \sX_s$ arising 
from a finite intersection of components of $H$ and $\sX_s$, there is associated
a canonical integral $\Gamma$-affine polyhedron $\Delta_S$ in $X^\an$ and 
such polyhedra form an atlas for the skeleton $S(\sX,H)$. Note that $S$ is obtained from a stratum $T$ of $\sX_s$ by intersecting 
with  horizontal components $H_{i_1}, \dots ,H_{i_p}$ of $H$. We get $\Delta_S$ by expanding $\Delta_T$ in linearly 
independent directions corresponding to $H_{i_1}, \dots ,H_{i_p}$ and hence we have $\Delta_S \cong \Delta_T \times \R_+^p$. 
In other words, the skeleton of a strictly semistable pair generalizes the
skeleton of a strictly semistable model in Berkovich's sense by allowing
unbounded faces.

We refer to \S \ref{Strictly semistable pairs}--\S \ref{section:functoriality} for a detailed study of these skeletons. In particular, 
we describe the closure $\hat{S}(\sX,H)$ of $S(\sX,H)$ in $\Xan$ which we call the {\it compactified skeleton}. The main 
result from these sections is the following.

\newtheorem*{deformation retraction}{Theorem \ref{deformation retraction}}
\begin{deformation retraction}
Let $(\sX,H)$ be a strictly semistable pair and let $X$ be the generic fibre of $\sX$. Then there is a canonical retraction map $\tau$ from 
$\Xan \setminus H^{\an}$ onto the skeleton $S(\sX,H)$ which extends to a 
 proper strong deformation retraction 
$\hat\tau$ from $\Xan$ onto the compactified skeleton $\hat{S}(\sX,H)$.  
\end{deformation retraction}

We formulate and prove this theorem in the setting of Raynaud's admissible formal schemes over $\kcirc$. The proof follows closely Berkovich's 
proof of the corresponding fact for $S(\sX)$ in \cite{berkovich:locallycontractible1}  taking  the unbounded part of our building blocks for the skeleton into account.  

If $X$ is a curve, the integral $\Gamma$-affine structure on the canonical
polyhedra in $S(\sX,H)$  amounts
to a metric structure on the edges and rays of $S(\sX,H)$.
In this case, the edge lengths are contained in $\Gamma$ and are induced by the
(logarithmic) modulus of various associated generalized open annuli in $X^\an$.

\paragraph[Slope formula] In his thesis \cite{thuillier:thesis}, Thuillier develops a non-archimedean potential theory on curves and proves an analogue of the Poincar\'e--Lelong equation. 
In ~\cite[Theorem~5.15]{bpr:analytic_curves}, an interpretation of the Poincar\'e--Lelong equation in terms of slopes on the skeleton is given. We generalize this slope formula to higher dimensions.

For a non-zero rational function $f$ on $X$, we show in Proposition
\ref{lem:piecewise.affine} that the skeleton $S(\sX,H)$ can be covered by
finitely many integral $\Gamma$-affine polyhedra $\Delta$ such that the
restriction $F$ of $-\log|f|$ to each $\Delta$ is integral $\Gamma$-affine.  The
latter means that $F|_{\Delta}$ is an affine function whose linear part is
given by a row vector with $\Z$ coefficients and that the constant term is in
the value group $\Gamma$.  Suppose now that $\div(f)$ is supported on the
boundary divisor $H$.  In this case $F$ is integral $\Gamma$-affine on each
canonical polyhedron $\Delta_S$ by Proposition~\ref{piecewise linear}.
For a canonical polyhedron $\Delta_S$ of dimension
$d\coloneq \dim(X)$, we define the slope $\slope(F;\Delta_T, \Delta_S)$ of $F$
at $\Delta_S$ along a codimension $1$ face $\Delta_T$. If $d=1$, this amounts to
the na\"ive outgoing slope along the edge or ray $\Delta_S$ emanating from the
point $\Delta_T$, relative to its metric. In higher dimensions however, it is
not clear in which direction in $\Delta_S$ one should measure the slope of $F$.
We define a canonical direction using some intersection numbers on the special
fibre $\sX_s$ (see Definition~\ref{slope-definition}).  With this in hand, we
define the divisor of $F$ as the formal sum
$$\hat\div(F) \coloneq \sum_{\Delta_T} \sum_{\Delta_S \succ \Delta_T } \slope(F;\Delta_T, \Delta_S) \, \Delta_T,$$
where $\Delta_T$ ranges over all  $d-1$-dimensional canonical polyhedra of $S(\sX,H)$ and where $\Delta_S$ ranges over all  $d$-dimensional canonical polyhedra containing $\Delta_T$. 
Then we  show the following slope formula for $S(\sX,H)$:

\newtheorem*{PL for pair}{Theorem \ref{PL for pair}}
\begin{PL for pair}
  Let $f\in K(X)^\times$ be a rational function such that
  $\supp(\div(f))\subset\supp(H)_\eta$ and let
  $F = -\log|f|\big|_{S(\sX,H)}$.  Then $F$ is continuous and 
  integral $\Gamma$-affine on each canonical polyhedron of $S(\sX,H)$, and we have 
  \[ \hat\div(F) = 0. \]
\end{PL for pair}

This is a kind of balancing condition
on $F$ which is a direct analogue of the balancing condition for tropical
varieties. The proof is based on the refined intersection theory of cycles with
Cartier divisors on admissible formal schemes over $\kcirc$ given
in~\cite{gubler:local_heights,gubler:lchs}. As the reader might be unfamiliar
with this intersection theory in non-noetherian situations, which we use at
several places in our paper, we  recall it in 
Appendix~\ref{Refined intersection theory with Cartier divisors}. In the end,
the slope formula  follows from the 
basic fact that the degree of a principal divisor intersected with the curve
given by the stratum closure of $T$ has degree $0$. 

From Theorem \ref{PL for pair} we  deduce a slope formula for the bounded skeleton $S(\sX)$ (see Theorem \ref{PL for bounded}).  
This formula is inspired by and
generalizes work of Cartwright~\cite{cartwright:tropical_complexes} on
tropical complexes, as well as the slope formula for curves
as formulated in~\cite[Theorem~5.15]{bpr:analytic_curves}.

A different higher-dimensional generalization of Thuillier's Poincar\'e--Lelong formula is given by Chambert--Loir and Ducros in~\cite[Theorem~4.6.5]{chambert_loir_ducros:forms_currents}.
It is formulated in terms of differential forms and currents on Berkovich spaces
using tropical charts and hence it is not directly related to our skeletal
approach.  The work of Chambert--Loir and
Ducros does not rely on a skeletal theory, and therefore applies to
essentially arbitrary analytic spaces; in contrast, our skeletal
version is quite explicit and is amenable to calculations (see below).

\subparagraph[A two-dimensional example] 
In \S \ref{A two-dimensional example}, we  illustrate skeletons and the slope formula in a non-trivial two-dimensional example 
which is obtained from the abelian variety $A=E^2$ for a Tate elliptic curve $E$. We choose a  regular triangulation of the canonical skeleton of $A$ leading 
to a $\kcirc$-model $\Acal$ of $A$ by Mumford's construction. Blowing up the closure of the origin $0$ of $A$ in $\Acal$, we obtain a strictly semistable pair $(\sX,H)$, where 
$H$ has five components given by the exceptional divisor and the strict transforms of the diagonal, the anti-diagonal, $E \times \{0\}$ and $\{0\} \times E$. Then we illustrate the 
slope formula for a certain rational function on the generic fibre $X$ with
support in the boundary divisor $H$. It is interesting to compare intersection
numbers on $\sX_s$ with the combinatorics of the skeleton of $(\sX,H)$.

\paragraph[Sturmfels--Tevelev multiplicity formula] 
Let $\varphi:U \rightarrow U'$ be a dominant generically finite morphism of varieties over $K$. We suppose that 
$U'$ (resp. $U$) is a closed subvariety of a multiplicative torus $T'$ (resp. $T$) and that 
$\varphi$ is the restriction of a homomorphism $T \rightarrow T'$. The  original Sturmfels--Tevelev multiplicity formula 
relates the tropical multiplicities of $\Trop(U)$ and $\Trop(U')$. It is proved in~\cite[Theorem~1.1]{sturmfels_tevelev:elimination} for
fields with a trivial valuation and
in~\cite[Corollary~8.4]{bpr:trop_curves} in general. The Sturmfels--Tevelev multiplicity formula is widely used in tropical geometry. For example, it is the basis for integration of differential forms on Berkovich spaces in~\cite{chambert_loir_ducros:forms_currents} and it is important for implicitization results (see \cite[\S 5]{sturmfels_tevelev:elimination}).

In Section \ref{Alterations},  we develop a similar formula relating the skeleton $S(\sX,H)$ of a strictly semistable pair $(\sX,H)$ as above to the tropical variety $\Trop(U')$ in the situation when 
$\varphi:U \coloneq X \setminus H \rightarrow U'$ is a dominant generically finite morphism to a closed subvariety $U'$ of a multiplicative torus $T$ with cocharacter group $N$. We prove in Proposition  \ref{cor:map.to.torus} that 
the map  $\trop\circ\phi:U^\an\to N_\R$ factors through the retraction 
$\tau:U^\an\to S(\sX,H)$ and that its restriction to $S(\sX,H)$ induces a
piecewise linear map $\varphi_{\rm aff}:S(\sX,H) \rightarrow N_\R$ with image
$\Trop(U')$. Moreover, the restriction of $\varphi_{\rm aff}$ to any canonical
polyhedron $\Delta_S$ of $S(\sX,H)$ is an integral $\Gamma$-affine map, i.e.\ it
is obtained from a linear map defined over $\Z$ and a translation by a
$\Gamma$-rational vector in $N_\R$.  

We consider a \emph{regular} point $\omega$ of $\Trop(U')$, which means that
$\omega$ has an integral $\Gamma$-affine polyhedron $\Delta$ as a neighbourhood
in $\Trop(U')$ such that the tropical multiplicity $m_\Trop(\Delta)$ of $\Delta$
in $\Trop(U')$ is well-defined.  See~\parref{par:tropicalization} for the
definition.  We assume that $\omega$ is not contained in a polyhedron
$\varphi_{\rm aff}(\Delta_S)$ of dimension $<d$ for any canonical polyhedron
$\Delta_S$ of $S(\sX,H)$.  Note that such points are dense in $\Trop(U')$. If
$\Delta_S$ is any canonical polyhedron of $S(\sX,H)$ with
$\omega \in \varphi_{\rm aff}(\Delta_S)$, then our assumptions imply that the
linear part of $\varphi_{\rm aff}$ induces an injective map
$N_{\Delta_S} \rightarrow N_\Delta$ between the underlying lattices of the
corresponding polyhedra.  Since $\dim(\Delta_S) = d$, the cokernel is finite and
hence we get a lattice index which we denote by
$[N_\Delta:N_{\Delta_S}] \coloneq \#\coker(N_{\Delta_S} \rightarrow N_\Delta)$.
(See also~\S\ref{par:polyhedra}).  We
prove the following variant of the Sturmfels--Tevelev multiplicity formula.

\newtheorem*{Sturmfels--Tevelev multiplicity formula for alterations}{Theorem \ref{Sturmfels--Tevelev multiplicity formula for alterations}}
\begin{Sturmfels--Tevelev multiplicity formula for alterations}
Under the hypotheses above, we have 
  \[ [U:U']\, m_\Trop(\Delta) = \sum_{\Delta_S}  [N_\Delta:N_{\Delta_S}], \]
  where the sum ranges over all canonical polyhedra $\Delta_S$ of the
  skeleton $S(\Xcal,H)$ with 
  $\relint(\Delta_S) \cap \phi_\aff^{-1}(\omega) \neq \emptyset$.
\end{Sturmfels--Tevelev multiplicity formula for alterations}

It follows that $\Trop(U')$ as a weighted polyhedral complex is essentially
determined by $S(\sX,H)$ and $\phi_\aff$.  The proof relies on similar
techniques from non-archimedean analytic geometry as the proof of the torus case
given in ~\cite[Corollary~8.4]{bpr:trop_curves}.

A special case of 
Theorem~\ref{Sturmfels--Tevelev multiplicity formula for alterations}
is proved for a smooth curve embedded as a
closed subscheme of a torus in~\cite[Corollary~6.9]{bpr:trop_curves}.
This formula relates the tropical multiplicity of an edge $e$ in the
tropicalization of the curve to the amount that the tropicalization map
``stretches'' the edges of the skeleton mapping to
$e$. Cueto~\cite[Theorem 2.5]{cueto:implicitization} also proves a version of Theorem 
\ref{Sturmfels--Tevelev multiplicity formula for alterations} for a closed
subvariety of a torus over a trivally valued field in characteristic $0$ with a
compactification whose boundary has  simple normal crossings.

\paragraph[Faithful tropicalization] The results outlined above allow one to compute tropicalizations in terms
of skeletons; those outlined below show that in certain situations, one
can do the reverse.  
The following {faithful tropicalization} result 
roughly says that a given skeleton can be ``seen'' in a suitable
tropicalization. 

\newtheorem*{thm:faithful.ss}{Theorem \ref{thm:faithful.ss}}
\begin{thm:faithful.ss}
 Let $(\sX,H)$ be a strictly semistable pair with generic fibre $X$.
 Then there exists a dense open subset $U$ of
$X$ and a morphism $\phi:U\to T = \bG_{m,K}^n$
such that the restriction $\phi_\aff$
of $\trop\circ\phi$ to $S(\sX,H)$ is a homeomorphism onto its image in $\R^n$ 
and is unimodular on every polyhedron of $S(\sX,H)$.
\end{thm:faithful.ss}

Note that $U$ may be a proper subset of $X \setminus H$ and hence $\varphi_{\rm aff}$
will not necessarily be affine on 
 canonical polyhedra.  The unimodularity condition roughly means that
 $\phi_\aff$ preserves the piecewise integral $\Gamma$-affine
 structure of the skeleton.  More formally, $\phi_\aff$ is 
{\it unimodular} provided that $S(\sX,H)$ has a finite covering by
integral $\Gamma$-affine polyhedra $\Delta$ which are  
contained in canonical polyhedra such that $\varphi_{\rm aff}$
restricts to an integral $\Gamma$-affine bijective map
$\Delta\rightarrow \Delta'$ 
onto an integral $\Gamma$-affine polyhedron $\Delta'$ of $N_\R = \R^n$ with
$[N_{\Delta'}:N_{\Delta}]=1$.  This is a local condition.
In the proof, one first uses local equations for strata on $\sX_s$ to produce a
$\phi$ such that $\phi_\aff$ is unimodular but not necessarily
globally injective.  When $\sX$ is quasiprojective, one 
can separate generic points of strata on $\sX_s$ using finitely many rational
functions; these functions in addition to $\phi$ give an injective unimodular
map.  In general one reduces to the quasiprojective case using Chow's lemma.

If $X$ is a curve, Theorem~\ref{thm:faithful.ss} says that there exists a
morphism $\phi:X\setminus H'\to\bG_{m,K}^n$ for a finite set of closed points
$H'\supset H$ such that $\phi_\aff$ is a homeomorphism and a local isometry from
$S(\sX,H)$ onto its image in $\R^n$, with respect to the lattice length on the target.
The unimodularity condition in this case translates into the local
isometry condition.
Such a result is proven in \cite[Theorem 6.22]{bpr:trop_curves}.

\paragraph[Section of Tropicalization]
One consequence of the Sturmfels--Tevelev multiplicity
formula is that if $\varphi:U \rightarrow U'$ is a birational morphism and if a polyhedron $\Delta$ of dimension $d = \dim(X)$ in
$S(\sX,H)$ maps to a $d$-dimensional polyhedron 
$\phi_\aff(\Delta)$ with tropical multiplicity one, then $\Delta$ is the
only maximal polyhedron mapping to $\phi_\aff(\Delta)$, and the
restriction of $\phi_\aff$ to $\Delta$ is unimodular.  From this it follows that
$\phi_\aff$ has a continuous partial section defined on $\phi_\aff(\Delta)$
which is also an integral $\Gamma$-affine map.   

Motivated by this observation, we prove the following general result on sections of tropicalization maps,
which makes no reference to semistable models or to skeletons.
Let $U$ be an (irreducible) very affine variety together with a closed immersion
$\varphi: U \hookrightarrow T \cong\bG_{m,K}^n$, and let $Z\subset\Trop(U)$ be a subset such that
every point of $Z$ has tropical multiplicity one.  Set 
$\trop_\phi = \trop\circ\phi^\an:U^\an\to N_\R$.

\newtheorem*{thm:section}{Theorem \ref{thm:section}}
\begin{thm:section}
  For every $\omega \in Z$, the affinoid space 
  $\trop_\varphi^{-1}(\omega)$ has a unique Shilov boundary
  point $s(\omega)$, and $\omega\mapsto s(\omega)$ defines a continuous partial
  section $s:Z\to U^\an$ of the tropicalization map 
  $\trop_\varphi : \Uan \rightarrow \Trop(U)$ on the subset $Z$.  
  Moreover, if $Z$ is contained in the closure of its interior in $\Trop(U)$, then $s$ is the unique continous section of $\trop_\phi$ defined on $Z$.
\end{thm:section}

An affinoid space has a unique Shilov boundary point if and only if the supremum
seminorm on its affinoid algebra is multiplicative, in which case the
Shilov boundary point is equal to the supremum seminorm.
In order to show that the resulting
section $s$ is continuous, we reduce to the case of a torus using a toric
Noether normalization argument, i.e.\ by choosing a
homomorphism $\alpha: T\to\G_m^d$ such that $\alpha\circ\phi$ is finite.

In the case of curves, such a result is proven in 
\cite[Theorem~6.24]{bpr:trop_curves}.
As a higher-dimensional example, the case of the Grassmannian $\Gr(2,n)$ of planes in $n$-space is studied in \cite{CHW}. The Pl\"ucker embedding of $\Gr(2,n)$ into projective space gives rise to a tropical Grassmannian $\mathcal{T} \Gr(2,n)$ in tropical projective space, which is an example of an extended tropicalization in the sense of \cite{payne:analytification}. Then \cite[Theorem 1.1]{CHW} states that the tropicalization map $\Gr(2,n)^{\an} \rightarrow \mathcal{T} \Gr(2,n)$ has a continuous section. 
Incidentally, the  construction of the section implies an algebraic result on the structure of the boundary components of the Grassmannian  \cite[Lemma 5.3]{CHW}. Note that Theorem \ref{thm:section} does not imply the continuity of the section on the whole tropical Grassmannian $\mathcal{T} \Gr(2,n)$.
Draisma and Postinghel \cite{dp} reprove this result with different techniques and use torus actions to obtain sections of the tropicalization map in other explicit situations. 

When there is a strictly semistable pair
$(\sX,H)$ such that $U = X\setminus\supp(H)$, where $X$ is the generic
fibre of $\sX$, we show that the image of the section $s$ is contained in
$S(\sX,H)$.   It follows that $s(Z)$ maps homeomorphically onto $Z$ under
$\phi_\aff$, and that $\phi_\aff$ is unimodular on $s(Z)$ in a suitable
sense: see Proposition~\ref{prop:section-skeleton}.
In other words, in this case one ``sees'' the skeleton in the
tropicalization.

\paragraph[Acknowledgements] Work  on this project started during an open
problem session of the Simons symposium on {\it Non-archimedean and Tropical
  Geometry} in 2013. The authors are very grateful to the Simons foundation for
financial support of this meeting and to the organizers, Matt Baker and Sam
Payne for the inspiring program. The authors also thank Dustin Cartwright,
Filippo Viviani and the referee for helpful comments. The first author was
partly supported by the collaborative research center SFB 1085 funded by the
Deutsche Forschungsgemeinschaft.  The second author was sponsored by the
National Security Agency Young Investigator Grant number H98230-15-1-0105.

\section{Preliminaries}

\paragraph[Notation and conventions] \label{par:notation}
An inclusion $A \subset B$ of sets allows the case $A= B$. The complement
of $A$ in $B$ is denoted by $B \setminus A$.  The sets $\N$ and $\R_+$ include
$0$. 

If $R$ is a ring with $1$, then the group of multiplicative units is denoted by $R^\times$. 

Throughout the paper, $K$ denotes an algebraically closed field endowed  with a non-trivial, non-archimedean, complete absolute value $|\scdot|$. Then $v\coloneq -\log|\scdot|$ is the corresponding valuation 
on $K$ with valuation ring $\kcirc \coloneq  \{ \alpha \in K \mid  |\alpha| \leq 1\}$, residue field $\ktilde$ and value group $\Gamma \coloneq  v(K^\times)$.  The maximal 
ideal $\{\alpha \in K \mid |\alpha| < 1 \}$ is denoted by $K^{\circ\circ}$. The corresponding point in $\Spec(\kcirc)$ is called the special point  $s$. We have $\Spec(\kcirc)=\{\eta,s\}$, where the generic point $\eta$ corresponds to the trivial ideal $\{0\}$. 

By an \emph{analytic space} we mean a $K$-analytic
space in the sense of Berkovich~\cite[\S 1.2]{berkovich:etalecohomology}. All analytic spaces which occur in this paper are good 
and hence we may also use the more restricted definition in~\cite{berkovich:analytic_geometry}. The
analytification functor from finite-type $K$-schemes to analytic spaces is
denoted $(\scdot)^\an$.  We distinguish between \emph{affinoid algebras} and \emph{strictly affinoid algebras} as in~\cite{berkovich:analytic_geometry} where they are called 
$K$-affinoid algebras and strictly $K$-affinoid algebras. Note that this is in contrast to \cite{bpr:trop_curves} and to the literature in rigid geometry as in~\cite{bgr:nonarch}, where affinoid means strictly affinoid. 
The Berkovich spectrum of an affinoid algebra $\sA$ is
denoted $\sM(\sA)$. 
Let $Y = \sM(\sA)$ be an affinoid space and let    
$\sA^\circ\subset \sA$ be the subring of power-bounded elements.
If $\sA$ is strictly affinoid, then the \emph{canonical model} of $Y$ is the $K^\circ$-formal scheme
$\Spf(\sA^\circ)$; this is an affine admissible formal scheme when $\sA$
is reduced by~\cite[Theorem~3.17]{bpr:trop_curves}.

A {\it variety} is an irreducible, reduced, and 
separated scheme of finite type over the base. A {\it very affine variety}
over a field is a variety which is isomorphic to a closed subvariety of a
multiplicative torus.  

If $X$ is a scheme over a ring $R$ and $R'$ is an $R$-algebra, the
extension of scalars is denoted $X_{R'} = X\tensor_R R'$.  Similarly, if
$X$ is a scheme over a base scheme $S$ and $S'\to S$ is a morphism, the
base change is denoted $X_{S'} = X\times_S S'$.

For a scheme $\sX$ over $\kcirc$, the fibre over $\eta$ is called 
{\it the generic fibre} and is denoted by $\sX_\eta$, and the fibre over
$s$ is called the {\it special fibre} and is denoted by $\sX_s$. Usually
we assume that 
$\sX$ is flat. Note that flatness for a $K^\circ$-variety $\sX$ is equivalent to $\sX_\eta
\neq \emptyset$.  
In this situation 
we call $\Xcal$ an {\it algebraic $\kcirc$-model} of the generic fibre $\Xcal_\eta$. 
For a Cartier divisor $D$ on a variety $\sX$, there is an
associated Weil divisor $\cyc(D)$ and an intersection theory with cycles
on $\sX$. As the variety $\sX$ need not be noetherian, this intersection
theory is not standard and we recall it in Appendix \ref{Refined
  intersection theory with Cartier divisors}. Here, we want to emphasize
that $\div(f)$ denotes the Cartier divisor associated to a non-zero
rational function $f$ on $\sX$ and $\cyc(f)$ is the associated Weil
divisor. If $\sX$ is proper over $\kcirc$, then we have a reduction map
$\red:\sX_\eta \rightarrow \sX_s$. 

Similarly, for an admissible formal scheme $\fX$ over $K^\circ$, we let
$\fX_s$ denote its special fibre and $\fX_\eta$ its generic fibre. We refer to 
\cite[Section 3.5]{bpr:trop_curves} for an expository treatment of admissible formal schemes.  The
generic fibre is an analytic space and the special fibre is a 
$\td K$-scheme. We call $\fX$ a {\it formal $\kcirc$-model of $\fX_\eta$}.  There is a canonical reduction map 
$\red:\fX_\eta\to\fX_s$.  If $\sX$ is a flat $K^\circ$-scheme of finite
type then its completion $\fX = \hat\sX$ with respect to any nonzero
element of $K^\cc$ is an admissible formal scheme, and $\sX_s = \fX_s$.
If $\sX$ is proper then $\fX_\eta = \sX_\eta^\an$.
We also use the notation $\cyc(\scdot),\div(\scdot)$ for Cartier and Weil
divisors on admissible formal schemes. 

We will generally denote schemes over $K^\circ$ using calligraphic letters 
$\sX,\sY,\ldots$ and admissible formal schemes over $K^\circ$ using German
letters $\fX,\fY,\ldots$  We will use Roman letters $X,Y,\ldots$ for both
schemes and analytic spaces over $K$.

The Tate algebra of restricted (i.e.\ convergent) power series in
indeterminates $x_1,\ldots,x_n$ 
with coefficients in $K$ (resp.\ $K^\circ$) is denoted
$K\angles{x_1,\ldots,x_n}$ (resp.\ $K^\circ\angles{x_1,\ldots,x_n}$).
The closed unit ball, viewed as an analytic space, is denoted by
$\bB$. It is the Berkovich spectrum $\sM(K\angles x)$. 
The formal ball $\Spf(\kcirc \langle x \rangle)$ will be denoted by $\fB$.  

If $M\cong\Z^n$ is a finitely generated free abelian group and $G$ is a
non-trivial additive subgroup of $\R$ then we set 
$M_G \coloneq M\tensor_\Z G$, regarded as a subgroup of the vector space
$M_\R$.

\paragraph[Integral $\Gamma$-affine polyhedra] \label{par:polyhedra}
Let $M\cong\Z^n$ be a finitely generated free abelian group and let 
$N = \Hom(M,\Z)$.  Let $\angles{\scdot,\scdot}:M\times N_\R\to\R$ denote
the canonical pairing.  An \emph{integral $\Gamma$-affine polyhedron} in 
$N_\R$ is a subset of $N_\R$ of the form
\[ \Delta = \big\{ v\in N_\R \mid \angles{u_i, v} + \gamma_i \geq 0 
\text{ for all } i = 1,\ldots,r \big\} \]
for some $u_1,\ldots,u_r\in M$ and $\gamma_1,\ldots,\gamma_r\in\Gamma$.
Any face of an integral $\Gamma$-affine polyhedron $\Delta$ is again integral
$\Gamma$-affine.  The relative interior of $\Delta$ is denoted $\relint(\Delta)$.  
A bounded polyhedron is called 
a {\it polytope}. An \emph{integral $\Gamma$-affine polyhedral complex} in $N_\R$ is a polyhedral
complex whose faces are integral $\Gamma$-affine. 
The \emph{dimension} of an  integral
$\Gamma$-affine polyhedral complex $\Sigma$ is 
$\dim(\Sigma)\coloneq \max\{\dim(\Delta)\mid\Delta\in\Sigma\}$.  We say
that $\Sigma$ has \emph{pure dimension $d$} provided that every maximal
polyhedron of $\Sigma$ (with respect to inclusion) has dimension $d$.

An \emph{integral $\Gamma$-affine function} on $N_\R$ is a function of the
form
\[ v \mapsto \angles{u,v} + \gamma ~:~ N_\R \to \R \]
for some $u\in M$ and $\gamma\in\Gamma$.  
More
generally, let $M'$ be a second finitely generated free abelian group and
let $N' = \Hom(M',\Z)$.  An \emph{integral $\Gamma$-affine map} from 
$N_\R$ to $N'_\R$ is a function of the form $F = \phi^* + v$, where
$\phi:M'\to M$ is a homomorphism, $\phi^*:N_\R\to N'_\R$ is the dual
homomorphism extended to $N_\R$, and $v\in N'_\Gamma$.  If $N' = M' = \Z^m$
and $F = (F_1,\ldots,F_m) : N_\R\to\R^m$ is a function, then $F$ is
integral $\Gamma$-affine if and only if each coordinate $F_i:N_\R\to\R$ is
integral $\Gamma$-affine.  

An \emph{integral $\Gamma$-affine map} from an integral $\Gamma$-affine
polyhedron $\Delta\subset N_\R$ to $N'_\R$ is by definition the restriction to
$\Delta$ of an integral $\Gamma$-affine map $N_\R\to N'_\R$.  If
$\Delta'\subset N'_\R$ is an integral $\Gamma$-affine polyhedron then a
function $F:\Delta\to\Delta'$ is \emph{integral $\Gamma$-affine} if the
composition $\Delta\to\Delta'\inject N'_\R$ is integral $\Gamma$-affine. 

Let $\Delta\subset N_\R$ be an integral $\Gamma$-affine polyhedron.
Let $(N_\R)_\Delta$ be the linear span of $\Delta-v$ for any $v\in\Delta$
and let $N_\Delta = N\cap (N_\R)_\Delta$.  This is a saturated subgroup of
$N$.  If $F:\Delta\to N'_\R$ is an integral $\Gamma$-affine map as
above then the image of $F$ is an integral $\Gamma$-affine polyhedron
$\Delta'$ in $N'_\R$.  Extending $F$ to $N_\R$, by definition the linear
part of $F$ takes $N$ into $N'$, hence induces a homomorphism
$F_\Delta: N_\Delta \to N'_{\Delta'}$ which is independent of the
extension of $F$ to $N_\R$.  We define the \emph{lattice index}
of $F$ to be 
\[ [N'_{\Delta'} ~:~ N_{\Delta}] \coloneq 
[N'_{\Delta'} ~:~ F_\Delta(N_\Delta)]
= \#\coker(F_\Delta). \]
We say that $F$ is \emph{unimodular} provided that it satisfies any of
the following equivalent conditions:
\begin{enumerate}
\item $F$ is injective and its lattice index is $1$.
\item Every integral $\Gamma$-affine function $f:\Delta\to\R$ is of the
  form $f'\circ F$ for an integral $\Gamma$-affine function 
  $f':\Delta'\to\R$.
\item $F$ is injective and the image of $N_\Delta$ under the linear part of $F$ is saturated in
  $N'$. 
\item $F$ is injective and the inverse function $\Delta'\to\Delta$ is
  integral $\Gamma$-affine.  (Recall that $\Delta'\coloneq F(\Delta)$.)
\end{enumerate}

\paragraph[Tropicalization] \label{par:tropicalization}
Let $T\cong\bG_{m,K}^n$ be an algebraic $K$-torus, let $M\cong\Z^n$ be the
character lattice of $T$, and let $N = \Hom(M,\Z)$ be its cocharacter
lattice.  For $u\in M$ the corresponding character of $T$ is denoted
$\chi^u$.  The \emph{tropicalization map} is the continuous, proper
surjection $\trop:T^\an\to N_\R$ defined by
\[ \angles{\trop(p), u} = -\log|\chi^u(p)|, \]
where $\angles{\scdot,\scdot}$ is the pairing between $M$ and
$N_\R$.  Choose a basis $u_1,\ldots,u_n$ for $M$, 
let $x_i = \chi^{u_i}$, and identify $N_\R$ with
$\R^n$ using the dual basis.  Then 
$K[M] = K[x_1^{\pm 1},\ldots,x_n^{\pm 1}]$, and we have
\[ \trop(p) = \big(-\log|x_1(p)|, \ldots, -\log|x_n(p)|\big). \]

If $\varphi: U \hookrightarrow T$ is a closed subscheme, 
we put $\trop_\varphi = \trop \circ \varphi^{\an}$, and we define the 
\emph{tropicalization} of $U$ to be the subset
\[ \Trop(U) \coloneq \trop_\varphi (U^\an) \subset N_\R. \]
The tropicalization map restricts to a continuous, proper surjection
$\trop_\varphi: U^\an\to\Trop(U)$.
By the Bieri--Groves theorem~\cite[Theorem~3.3]{gubler:guide},
$\Trop(U)$ is the support of an integral $\Gamma$-affine polyhedral complex
$\Sigma$ in $N_\R$.  If $U$ is a variety of
dimension $d$ then $\Sigma$ has pure dimension $d$.

Via the closed embedding $\varphi$ we have 
$ U \cong K[x_1^{\pm1}, \ldots, x_n^{\pm1}] / \mathfrak{a}$ for some ideal
$\mathfrak{a}$ in the Laurent polynomial ring. Fix 
$\omega = (\omega_1, \ldots, \omega_n) \in \rdop^n$,
and put $r_i = \exp(- \omega_i) \in \rdop$.  We write 
$r = (r_1, \ldots, r_n)$, $x = (x_1, \ldots, x_n)$ and use multi-index notation where
convenient. In particular, we put 
$K\angles{r^{-1} x, r x^{-1}} 
= K \angles{ r_1^{-1} x_1, \ldots, r_n^{-1} x_n, r_1 x_1^{-1}, \ldots, r_n x_n^{-1}}$.
The poly-annulus 
$\Mcal(K\angles{r^{-1} x, r x^{-1}}) = \trop\inv(\omega)$ is an affinoid subdomain of 
$T^{\an}$, which is strict if and only if
$\omega_1,\ldots,\omega_n\in\Gamma = v(K^\times)$.  The Banach norm on $K\angles{r^{-1} x, r x^{-1}}$ is denoted by $\| \sum_I a_I x^I\|_r = \max_I |a_I| r^I$. Assume that $\omega \in \Trop(U)$, and put 
\[ A_\omega = K\angles{r^{-1} x, r x^{-1}} / \mathfrak{a} K\angles{ r^{-1} x, r x^{-1} }.\]
  Then $U_\omega\coloneq \trop_\varphi^{-1}(\omega) $
can be identified with the affinoid subdomain $\Mcal(A_\omega)$ of
$U^{\an}$.

Choose an algebraically closed, complete valued
extension field $L/K$ whose value group $\Gamma_L$ is large enough so that
$\omega\in N_{\Gamma_L}$.  Choose $t\in T(L)$ such that 
$\trop(t) = \omega$.  The \emph{initial degeneration} 
$\In_\omega(U)$ of $U$ at $\omega$ is the special fibre of the schematic
closure of  
$t\inv U_{L}$ in the $L^\circ$-torus $T_{L^\circ} = \Spec(L^\circ[M])$.
This means that $\In_\omega(U)$ is a closed subscheme of 
$T_{\td L} = \Spec(\td L[M])$.  See~\cite[\S5]{gubler:guide} for a discussion of
initial degenerations and the dependence on $L$ and $t$.

 Let $g$ be a non-zero element of 
$L\angles{r^{-1}x, r x^{-1}}$.  Since each $r_i = \exp(- \omega_i)$ is contained in the value       
group of $L$, we have $\|g\|_r = |c|$ for $c \in L$. As  
$|t_i| = r_i$ for
all $i$, the Laurent series $c^{-1} g(t x)$ is an element of
$L\angles{x, x^{-1}}^\circ$. Its image under the reduction map
\[ L\angles{x, x^{-1}}^\circ \longrightarrow \widetilde{L}[x, x^{-1}] \] 
is the \emph{initial form} of $g$, which we denote by $\In_\omega (g)$.  

The \emph{initial ideal} of $\mathfrak{a}\subset K[x_1^{\pm 1},\ldots,x_n^{\pm1}]$ is defined as the ideal generated by 
all initial forms of polynomials in $\mathfrak{a}$:
\[ \In_\omega(\mathfrak{a}) 
= \big(\In_\omega(g)\mid g \in \mathfrak{a} \big) \subset \td L[x,x\inv]. \] 
The closed subscheme of 
$T_{\td L} = \Spec(\td L[M])$ given by the initial ideal $  \In_\omega(\mathfrak{a})$ is the initial degeneration $\In_\omega(U)$.
The initial degeneration is well-defined up to translation by elements
of $T_{\td L}(\td L)$, and up to the choice of the field
$L$.  See~\cite[Definition~10.6]{gubler:guide}.  If $\omega\in N_\Gamma$
we will always take $L = K$.

Let $\omega\in\Trop(U)$.  The \emph{tropical multiplicity} 
$m_\Trop(\omega)$ of the point $\omega$ is the number of irreducible
components of $\In_\omega(U)$, counted with multiplicity.  This quantity
is independent of all choices involved in the definition of
$\In_\omega(U)$.  
Suppose now that $U$ is a variety of dimension $d$.  
Let $\Sigma$ be an integral $\Gamma$-affine polyhedral complex with support
$\Trop(U)$ and let $\Delta\in\Sigma$ be a maximal (i.e.\ $d$-dimensional) polyhedron.
If $\Sigma$ is sufficiently fine then for $\omega\in\relint(\Delta)\cap N_\Gamma$
the initial degeneration $\inn_\omega(U)$ is isomorphic to 
$Y\times\bG_{m,\td K}^d$ for a dimension-zero $\td K$-scheme $Y$.  
Different choices of $\omega\in\relint(\Delta)\cap N_\Gamma$ give rise to
isomorphic initial degenerations. By \cite[Section 3.3]{maclagan_sturmfels:book} the multiplicity $m_\Trop(\omega)$ is constant on the relative interior of $\Delta$, hence we call it 
the \emph{tropical multiplicity $m_\Trop(\Delta)$} of $\Delta$.
Equivalently $m_\Trop(\Delta)$  is the length of the dimension-zero scheme $Y$:
this is explained in~\cite[Theorem~4.29]{bpr:trop_curves}.

More generally, if $\phi:U\to T$ is any morphism then we define $\Trop(U)$
to be the tropicalization of the schematic image $U'$ of $U$.  Initial
degenerations and tropical multiplicities are all defined with respect to
$U'$.

\section{Strictly semistable pairs}  \label{Strictly semistable pairs}

In this section, we give a variant of de Jong's notion 
of a strictly semistable pair $(\Xcal,D)$ in the case of a base $\Spec(K^\circ)$ for the
valuation ring $K^\circ$ of our valued field $K$. 
Roughly speaking, a strictly semistable pair $(\Xcal,D)$ over $K^\circ$
consists of a strictly semistable proper scheme $\Xcal$ over $K^\circ$ and
a divisor $D$ on $\Xcal$ with strictly normal crossings which includes the
special fibre.  It is convenient to only include the horizontal part of
the divisor as part of the data, as it is Cartier whereas the vertical
part may not be.

\begin{defn} \label{strictly semistable pairs}
  A {\it strictly semistable pair} $(\Xcal, H)$ consists of an irreducible proper flat
  scheme $\Xcal$ over $K^\circ$ and a sum
  $H = H_1 + \cdots + H_S$ of distinguished effective Cartier divisors
  $H_i$ on $\sX$ such that $\Xcal$ is covered by open  
  subsets $\Ucal$ which admit an \'etale morphism 
  \begin{equation}
    \label{eq:local.sss.pair}
    \psi ~:~ {\Ucal} \longrightarrow {\Scal} \coloneq
    \Spec\big( K^\circ[ x_0, \dots, x_d ] / \angles{x_0 \cdots x_r - \pi} \big) 
  \end{equation}
  for $r \leq d$ and $\pi \in K^\times$ with $|\pi|<1$.  We assume that
  each  $H_i$ has irreducible support and that the restriction of $H_i$ to $\sU$ is either trivial or defined by 
  $\psi^*(x_{j})$ for some $j \in \{r+1, \dots, d\}$.
\end{defn}

\paragraph \label{comments for strictly semistable} 
The underlying scheme
$\Xcal$ is called a {\it strictly semistable scheme} over $\kcirc$. Note
that such a scheme is a variety over $\kcirc$ in the sense of~\parref{par:notation}.
Indeed, reducedness follows from \cite[Proposition 17.5.7]{egaIV_4}.
Note that  $\Xcal$ is not noetherian, but the underlying topological space is noetherian and of topological 
dimension $d+1$. This follows from the fact that the generic fibre and the special fibre are both noetherian. 

Similarly, de Jong \cite{dejong:alterations} defined strictly semistable schemes and strictly semistable pairs over a complete discrete valuation ring $R$,
where $\pi$ is a uniformizer. Suppose that $R$ is a subring of our algebraically closed field $K$ and that the discrete valuation of $R$ 
extends to our given complete valuation $v$. 
By  \cite[2.16]{dejong:alterations}, the base change of a strictly semistable scheme over $R$ (in the sense of de Jong) to the valuation ring $\kcirc$ is a strictly semistable scheme as above. By  \cite[6.4]{dejong:alterations}, the same applies for strictly semistable pairs if we neglect the vertical components of the divisor.

\begin{notn} \label{notn:ssp.notation}
  We fix the following notation for a strictly semistable pair
  $(\sX, H)$.  The generic fibre of $\sX$ is denoted by $X$.
  Let $\sH_i$ be  the closed subscheme of $\sX$ locally cut out
  by a defining equation for $H_i$, and let $\sH = \bigcup_{i=1}^S\sH_i$.  
   Let $V_1,\ldots,V_R$ be the irreducible components of the special fibre $\sX_s$.  Define
  $D_i = V_i$ for $i\leq R$ and $D_i = \sH_{i-R}$ for 
  $R < i \leq N\coloneq R+S$, and set $D = \sum_{i=1}^N D_i$.  This is a
  Weil divisor on $\sX$.  We call the $\sH_i$ (resp.\ $V_j$) the 
  \emph{horizontal components} (resp.\ \emph{vertical components}) of $D$.
\end{notn}

\paragraph \label{components of D} 
Let $(\sX, H)$ be a strictly semistable pair.
It is clear that the generic fibre of each $\sH_i$ is smooth.
Since $H_i$ is Cartier on $\sX$, the special fibre $(\sH_i)_s$ is the
support of a Cartier divisor on $\sX_s$, so it has pure codimension one in
$\sX_s$. 
Each $V_j$ is a smooth $\td K$-scheme because, locally, we have
$\psi^*(x_k) = 0$ on $V_j\cap\sU_s$ for some $k$, so $V_j\cap\sU_s$ is
\'etale over an affine space. We may regard $V_j$ as a Weil divisor on
$\sX_s$, but it is not necessarily the support of a Cartier divisor on
$\sX$: see Proposition~\ref{prop:comps.are.cartier}.

The generic fibre $X$ of $\sX$ is smooth as the generic fibre of the
scheme $\sS$ of~\eqref{eq:local.sss.pair} is smooth.
It is clear that $d=\dim(X)$, but $r, s$ may depend on the choice of $\Ucal$.
In this generality, $v(\pi)$ may also depend on the choice of $\sU$; when
$X$ is a curve, this reflects the fact that the edges of the skeleton
of $\sX$ may have different lengths.  This is related to the fact that the
$V_j$ may not be Cartier; again see Proposition~\ref{prop:comps.are.cartier}.

\begin{eg}
  Let $(\sX,H)$ be a strictly semistable pair of relative dimension one.
  Then its generic fibre $X = \sX_\eta$ is a smooth, proper, connected $K$-curve, 
  and its special fibre $\sX_s$ has smooth irreducible components and at
  worst nodal singularities.  The Cartier divisor $H$ amounts to a finite
  collection of points in $X(K)$ which reduce to distinct smooth points of
  $\sX_s(\td K)$.
\end{eg}

\begin{rem} \label{inductive pairs} 
  Let $(\sX,H)$ be a strictly semistable pair and let $\sH_k$ be a
  horizontal component of $D$.
  Then $(\sH_k,H|_{\sH_k})$ is again a strictly semistable pair, where
  $H|_{\sH_k} = \sum_{j\neq k} H_j|_{\sH_k}$.
  This is immediate from Definition~\ref{strictly semistable pairs}.
  Note however that $H_j|_{\sH_k}$ does not necessarily have irreducible
  support; it must be broken up into a sum of irreducible Cartier divisors.
\end{rem}

\medskip
It is also useful to have a notion of a strictly semistable pair in the
category of admissible formal schemes.  
Definition~\ref{strictly semistable pairs} carries over verbatim.
The definition is constructed so as not to allow the horizontal components to
intersect themselves in the generic fibre --- this condition is
not local in the analytic topology.

\begin{defn} \label{defn:formal.ss.pair}
  A {\it formal strictly semistable pair} $(\fX, H)$ consists of a
  connected quasi-compact 
  admissible formal $K^\circ$-scheme $\fX$ and a sum
  $H = H_1 + \cdots + H_S$ of distinguished effective Cartier divisors
  on $\fX$ such that $\fX$ is
  covered by formal open subsets $\fU$ which admit an \'etale morphism 
  \begin{equation}
    \label{eq:formal.sss.pair}
    \psi ~:~ \fU \To 
    \Spf\big( K^\circ\angles{x_0, \dots, x_d} / \angles{x_0 \dots x_r - \pi} \big) 
  \end{equation}
  for $r \leq d$ and $\pi \in K^\times$ with $|\pi|<1$.  
  We assume that each $H_i|_{\fX_\eta}$ has irreducible support and that 
  $H_i|_\fU$ is defined by $\psi^*(x_j)$ for some $j > r$ unless it is
  trivial. 
\end{defn}

\paragraph \label{par:formal.ssp.notn}
We use notation analogous to~\parref{notn:ssp.notation} for formal
strictly semistable pairs.  That is, we let $X=\fX_\eta$ be the generic fibre of
$\fX$ and $\fX_s$ its special fibre.  We define $\fH_i$ as the admissible
formal closed subscheme of 
$\fX$ locally cut out by a defining equation for $H_i$.  
Its generic fibre
$(\fH_i)_\eta$ is an irreducible Weil divisor on $X$ and its special fibre
is a Weil divisor on $\fX_s$.   We also set $\fH = \bigcup_{i=1}^s\fH_i$.
Let $V_1,\ldots,V_R$ be the irreducible
components of $\fX_s$, and define $D_i = V_i$ for $i\leq R$ and
$D_i = \fH_{i-R}$ for $R<i\leq N\coloneq R+S$.  Set
$D = \sum_{i=1}^N D_i$.  We may regard $D$ as a Weil divisor on $\fX$ in
the sense of~\cite[\S3]{gubler:local_heights}, with horizontal components
$\fH_i$ and vertical components $V_j$.  (In \textit{loc.\ cit}.\ the
horizontal divisors live on the generic fibre $X$, but for our purposes it
is convenient to remember the special fibre of the $\fH_i$.) 

The remarks made in~\parref{components of D} apply to formal strictly
semistable pairs as well.

\paragraph \label{formal strictly semistable scheme}
Let $\fX$ be a {\it  strictly semistable formal scheme} over $\kcirc$ which 
means that $(\fX,0)$ is a formal strictly semistable pair. We consider a 
formal open subset $\fU$ of $\fX$. Then $\fU$ is formal affine if and 
only if $\fU_s$ is an affine open subscheme of $\fX_s$. Indeed, the special 
fibre of $\fX_s$ is reduced and hence we may use \cite[Proposition 1.11]{gubler:local_heights}
to deduce the claim from a theorem of Bosch 
\cite[Theorem 3.1]{bosch:rigid_raume} about formal analytic varieties.
Note that in this case, $\red\inv(\fU_s)$ is an affinoid domain in $\fX_\eta$.

\begin{eg} \label{standard pair} Let $r,s,d \in \N$ with $r+s \leq d$. In the
  case $r>0$, we choose $\pi \in K^{\times}$ with $|\pi|<1$. If $r=0$, then we
  always take $\pi \coloneq 1$ to avoid ambiguities.  Let
  $\fB = \Spf(K^\circ\angles x)$ be the formal ball of radius $1$ as
  in~\S\ref{par:notation} and let
  $\fU_{\Delta(r,\pi)} = \Spf(K^\circ\angles{x_0,\ldots,x_r}/(x_0\cdots
  x_r-\pi))$
  be the canonical model of the polytopal domain~\cite[\S6]{gubler:guide} in the
  hyperplane $x_0 \cdots x_r = \pi$ of $\bdop^{r+1} \coloneq \fB_\eta^{r+1}$
  associated to the simplex
\begin{equation} \label{standard simplex}
\Delta(r,\pi)\coloneq\{\bv \in \rdop_+^{r+1} \mid v_0 + \dots + v_{r} = v(\pi) \}.
\end{equation}
We consider the strictly semistable formal scheme $\fS\coloneq \fU_{\Delta(r,\pi)} \times \fB^{d-r}$.
Note that 
\[ \fS \cong \Spf\big(K^\circ\angles{x_0,\ldots,x_d}/\angles{x_0\cdots x_r-\pi}\big). \]
Then $(\fS, H(s))$ is a formal
strictly semistable pair, where $H(s)$ is the principal Cartier divisor
defined by $x_{r+1}\cdots x_{r+s}$. We call  $(\fS,H(s))$ a {\it standard pair}. 
The isomorphism 
class of the formal scheme $\fS$ of a standard pair is uniquely determined by $(r,d,v(\pi))$.
\end{eg}

\begin{rem} \label{pull-back of standard pairs}
By the definitions, any formal strictly semistable pair $(\fX,H)$ is covered by formal open subsets $\fU$ with an \'etale morphism 
\begin{equation} \label{etale morphism and standard}
\psi ~:~ \fU \To 
    \fS =\Spf\big( K^\circ\angles{x_0, \dots, x_d} / \angles{x_0 \dots x_r - \pi} \big) 
\end{equation}
to the formal scheme $\fS$ of a standard pair $(\fS,H(s))$ such that $H|_\fU = \psi^*(H(s))$. We have only to note that  the morphism in \eqref{eq:formal.sss.pair} 
is not changed by our 
choice $\pi=1$ in case of $r=0$.
\end{rem}

\paragraph \label{completion}
For a proper, flat $K^\circ$-scheme we let $\hat\sX$ be its completion
with respect to a non-zero element of $K^\cc$.  
Let $\sU\subset\sX$ be an open subset which admits an \'etale morphism
$\psi:\sU\to\sS$ as in~\eqref{eq:local.sss.pair}.  Taking completions, we
have an \'etale morphism in the category of admissible formal
$K^\circ$-schemes 
\[ \psi ~:~ {\hat\sU} \To {\fS} \coloneq
\Spf\big( K^\circ\angles{x_0, \dots, x_d} / \angles{x_0 \dots x_r - \pi} \big). \]
Since $\sX$ is proper, the analytification of $\sX_\eta$ is naturally
identified with the analytic generic fibre of $\hat\sX$.
A Cartier divisor $H$ on $\sX$ naturally induces a Cartier divisor 
$\hat H$ on $\hat\sX$, and the analytification of an irreducible closed
subscheme of $\sX_\eta$ is an irreducible Zariski-closed subspace of
$\sX_\eta^\an$ \cite[Theorem 2.3.1]{conrad:irredcomps}.  Hence we have shown:

\begin{prop} \label{prop:formal.alg.sss.pair}
  If $(\sX,H)$ is a strictly semistable pair then $(\hat\sX,\hat H)$ is a
  formal strictly semistable pair. 
\end{prop}

\begin{rem} \label{formal model}
  The converse to Proposition~\ref{prop:formal.alg.sss.pair} is also true:
  if $\sX$ is an irreducible proper flat $K^\circ$-scheme and $H$ is an effective
  Cartier divisor on $\sX$ such that $(\hat\sX,\hat H)$ is a formal
  strictly semistable pair then $(\sX,H)$ 
  is a strictly semistable pair.  

We sketch the argument. Let $x \in \sX_s$. There is a formal open subset  $\fU=\Spf(A)$ of $\hat\sX$ containing $x$ with a 
morphism $\psi$ with the same properties as in \ref{eq:formal.sss.pair}. Since $\sX$ is algebraic, we may assume that $\fU$ is 
the formal completion of an open subscheme $\sU=\Spec(B)$ of $\sX$. For $i=0,\dots,d$, we approximate $\gamma_i:=\psi^*(x_i) \in A$ 
by sufficiently close $b_i \in B$. Using $\gamma_0 \cdots \gamma_r = \pi$, we   get $b_0 \cdots b_r = \pi +r$ with 
$r \in B \cap \rho A$ for some $\rho \in K$ with $|\rho|<|\pi|$. Replacing $b_0$ by $b_0(1+r/\pi)^{-1}$ and shrinking the open neighbourhood $\Ucal$ of $x$ to make this a regular function,
we may assume that $b_0 \cdots b_r = \pi$. Then the map $y \mapsto (b_0(y), \dots , b_d(y))$ induces a morphism $\psi':\Ucal \to \Scal$ as in \eqref{eq:local.sss.pair}. 
Since $\psi_s'=\psi_s$ is \'etale, the morphism $\psi'$ is \'etale on an open neigbhourhood of $\Ucal_s=\fU_s$ in $\Ucal$. This yields easily the claim.
\end{rem}

\paragraph \label{stratification} 
The Weil divisor $D = \sum_{i=1}^N D_i = \sum_{i=1}^R V_i + \sum_{j=1}^S \sH_j$
of a strictly semistable pair $(\Xcal,H)$ has a stratification, where a {\it stratum} $S$ is given
as an irreducible component of 
$\bigcap_{i \in I} D_i \setminus \bigcup_{i \not \in I} D_i$ for any 
$I \subset \{1, \dots, N\}$. This is a special
case of a more general definition of strata given in
\cite{berkovich:locallycontractible1}, \S 2.  It is easy to see that the
closure $\overline{S}$ of a stratum $S$ is a strata subset, i.e.\ a
disjoint union of strata. Note 
also that two strata are either disjoint or one of them is contained in
the closure of the other. 
The set of strata of $D$ is partially ordered by $S \leq T$ if and only if $\overline{S} \subseteq \overline{T}$.

If $S$ is contained in the generic fibre, then
it follows from Remark~\ref{inductive pairs} that 
$(\overline{S},H|_{\bar S})$ is a strictly semistable pair.
We let $\str(\sX_s,H)$ denote the set of vertical strata, i.e.\ the strata
contained in the special fibre $\sX_s$. Note that the closure of a vertical stratum is smooth over $\ktilde$. 

\subparagraph 
The Weil divisor $D = \sum_{i=1}^R V_i + \sum_{j=1}^S \fH_j$
of a formal strictly semistable pair $(\fX,H)$ has a stratification
defined in the same way, where we consider $\fH_j$ as a disjoint union of
its generic and special fibres.  If we handle the horizontal strata with some care, then
everything from above applies. We denote the set of vertical strata by $\str(\fX_s,H)$.
Let us explain how we treat a horizontal stratum $S$ of $D$. The closure of $S$ in  the generic fibre $\fX_\eta$
is a closed analytic subvariety $Y$ of $\fX_\eta$. The closure of $Y$ in $\fX$,
defined similarly to the scheme-theoretic closure  
in algebraic geometry, is an admissible formal scheme $\fY$ over $\kcirc$ which is a closed formal subscheme of $\fX$ and which has 
generic fibre $Y$ (see \cite[Proposition 3.3]{gubler:local_heights}).  Setting $\overline{S}\coloneq \fY$,  we can show as above that 
$(\overline{S},H|_{\bar S})$ is a formal strictly semistable pair. For the definition of the partial ordering $\leq$, we always  view $\overline{S}$ as the strata subset given by the disjoint union of $\fY_\eta$ and $\fY_s$.

\section{The skeleton of a strictly semistable pair} 
\label{The skeleton of a strictly semistable pair} 

In this section, $(\fX,H)$ denotes a formal strictly semistable pair.  As always
we use the associated notation introduced in~\parref{par:formal.ssp.notn}. 
We will define the
skeleton $S(\fX,H)$
which generalizes the skeleton for strictly
semistable formal schemes over $K^\circ$ introduced by Berkovich in
\cite{berkovich:locallycontractible1} and
\cite{berkovich:locallycontractible2}. The horizontal divisor $H$ is the
new ingredient here.   The skeleton $S(\fX,H)$ is well-known in the
case of curves (see e.g.~\cite{tyomkin} and~\cite{bpr:analytic_curves}). 
In particular, we  will obtain the skeleton of an algebraic strictly semistable pair 
$(\Xcal,H)$ by using the formal completion $(\hat{\Xcal},\hat{H})$ as
in~\ref{completion} and setting $S(\Xcal,H)\coloneq S(\hat{\Xcal},\hat{H})$.  

We will define $S(\fX,H)$ in such a way that the vertical strata
$S\in\str(\fX_s,H)$ are in one-to-one correspondence with the open faces of the
skeleton.  (By ``open face'' we mean the relative interior of a face --- see
Remark~\ref{relation between faces}.)  We will do so first on the standard pairs
$(\fS,H(s))$ from~\ref{standard pair} and then on $\fX$ using the \'etale
morphisms $\psi:\fU\to\fS$ of~\eqref{etale morphism and standard}.  For this, it
is easier if there is a unique stratum $S$ on $\fU_s$ which maps to the minimal
vertical stratum of $(\fS,H(s))$; in this case, the part of the skeleton of
$(\fX,H)$ contained in $\fU_\eta$ will map homeomorphically onto the skeleton of
$(\fS,H(s))$.  The following proposition states that such $\fU$ exist for each
stratum $S\in\str(\fX_s,H)$. Recall from \ref{par:formal.ssp.notn} that $D$ is
the Weil divisor on $\fX$ given as the sum of the Weil divisor induced by $H$
and the Weil divisor induced by the special fibre of $\fX$.

\begin{prop} \label{normalization of covering} 
  Let $(\fX,H)$ be a formal strictly semistable pair.  Any open covering of
  $\fX_s$ admits a refinement $\{\fU_s\}$ by affine open subsets
  $\fU_s$ satisfying the following properties:
  \begin{itemize}
  \item[(a)] The formal open subscheme $\fU$ of $\fX$ with underlying set $\fU_s$
    admits an \'etale morphism $\psi:\fU\to\fS = \fU_{\Delta(r,\pi)}\times\fB^{d-r}$ as
    in~\eqref{etale morphism and standard} such that $(\fU,H|_{\fU})$ is the pull-back of the standard pair $(\fS,H(s))$ for some $s \in \{0,\dots, d-r\}$. 
  \item[(b)] There is a distinguished vertical stratum $S$ of $D$ associated to
    $\fU_s$ such that for any stratum $T$ of $D$, we have 
    $S \subset \overline{T}$ if and only if $\fU_s \cap \overline{T} \neq \emptyset$.
  \item[(c)]  The distinguished stratum $S$ from (b) is given on $\fU_s$ 
    by $\psi^{-1}(x_0=\dots=x_{r+s}=0)$ in terms of the \'etale morphism
    $\psi$ in~(a).
  \item[(d)] Every vertical stratum of $D$ is the distinguished stratum of
    a suitable $\fU_s$. 
  \end{itemize} 
\end{prop}

\smallskip
The arguments are similar as in the proof of Proposition 5.2 in
\cite{gubler:canonical_measures}. We leave the details to the reader. 

When choosing a covering of $\fX$ as in 
Definition~\ref{defn:formal.ss.pair}, we will always pass to a refinement
whose special fibre $\{\fU_s\}$ satisfies 
Proposition~\ref{normalization of covering}.  
The formal  open
subschemes $\fU$ induced on the open subsets $\fU_s$ are called the
\emph{building blocks} of the formal strictly semistable pair $(\fX,H)$. By \ref{formal strictly semistable scheme}, building blocks 
are formal affine open subschemes of $\fX$. For any building block $\fU$, the \'etale morphism $\psi$ from (a) induces a bijective correspondence between $\str(\fU_s,H|_{\fU})$ and the vertical strata of the standard pair $(\fS,H(s))$ of~\parref{standard pair}.

\paragraph[The skeleton of a standard pair] \label{par:skeleton.standard} 
We start by defining the skeleton of the standard 
pair $(\fS,H(s))$ of~\parref{standard pair}, where
$\fS = \fU_{\Delta(r,\pi)} \times \fB^{d-r} 
= \Spf(K^\circ\angles{x_0,\ldots,x_d}/\angles{x_0\cdots x_r-\pi})$. 
For $\ve\in |K^\times|$ with $0<\ve <1$,
consider the affinoid annulus $U_\ve \coloneq \{ x \in\bB^1 \mid |x| \geq \ve \}$.
The canonical model $\fU_\ve$ of $U_\ve$ is an admissible formal affine scheme over
$K^\circ$ which is strictly semistable; 
indeed, we have 
$\fU_\ve\cong\Spf(K^\circ\angles{y_0,y_1}/\angles{y_0y_1-a_\epsilon})$, 
where $a_\epsilon\in K^\times$ is any element with $|a_\epsilon| = \epsilon$.
We conclude that 
$\fU_{\Delta(r,\pi)} \times \fU_\ve^{s} \times\fB^{d-r-s}$ is a strictly
polystable formal scheme and we define its 
skeleton $S(\fU_{\Delta(r,\pi)} \times \fU_\ve^{s} \times\fB^{d-r-s})$ as
in \cite[\S5]{berkovich:locallycontractible1}, or in 
\cite[\S4]{berkovich:locallycontractible2}.  
This is a closed subset of the generic fibre
$U_{\Delta(r,\pi)}\times U_\ve^s\times\bB^{d-r-s}$ which 
is homeomorphic to $\Delta(r,\pi) \times [0, -\log \ve]^s$. 
It has the following explicit description.
First note that projection onto $x_1,\ldots,x_d$ induces an isomorphism of 
$U_{\Delta(r,\pi)} \times U_\ve^{s} \times\bB^{d-r-s}$ onto
the affinoid domain
\[\begin{split}
  U = \big\{ p \in\bB^d \mid & -\log|x_1(p)| - \cdots - \log|x_r(p)| \leq v(\pi), \\
  &\,\text{$ -\log|x_j(p)| \leq -\log\epsilon$ for $j=r+1, \dots,r+s$} \big\} 
\end{split}\]
in $\bB^d$, and similarly that projection onto $v_1, \dots, v_{r+s}$ 
maps $\Delta(r,\pi)\times[0,-\log\epsilon]^s$ homeomorphically onto
\[ S = \big\{ (v_1,\ldots,v_{r+s}) \in \R_+^{r+s} \mid 
v_1 + \cdots + v_r \leq v(\pi), ~ 
 \text{$v_j \leq -\log\epsilon$ for $j=r+1, \dots,r+s$}  \big\}. \]
For $\bv = (v_1,\ldots,v_{r+s})\in S$ we define a bounded multiplicative norm $\|\cdot\|_\bv\in U$ by the formula
\[ \bigg\|\sum_{\mb} a_\mb x^\mb\bigg\|_\bv
= \max_\mb |a_\mb|\, \exp(-m_1 \, v_1 - \ldots - {m_{r+s}} \, v_{r+s}), \]
where $\mb = (m_1,\ldots,m_d)$ ranges over $\Z^{r+s} \times \N^{d-r-s}$. 
This gives a continuous inclusion $S\inject U$, hence an inclusion
$\Delta(r,\pi)\times[0,-\log\epsilon]^s\inject 
U_{\Delta(r,\pi)} \times U_\ve^{s} \times\bB^{d-r-s}$.  The image is the
skeleton $S(\fU_{\Delta(r,\pi)} \times \fU_\ve^{s} \times\fB^{d-r-s})$.
Note that 
\[ p \mapsto (-\log|x_0(p)|,\ldots,-\log|x_{r+s}(p)|) ~:~
S(\fU_{\Delta(r,\pi)} \times \fU_\ve^{s} \times\fB^{d-r-s}) \To
\Delta(r,\pi)\times[0,-\log\epsilon]^s \]
is a continuous left inverse to the above inclusion.

We define 
$S(\fS,H(s))\subset \fS_\eta = U_{\Delta(r,\pi)}\times\bB^{d-r}$ 
as the union of the skeletons 
$S(\fU_{\Delta(r,\pi)} \times \fU_\ve^{s} \times\fB^{d-r-s})$ as 
$\ve \to 0$.  We see by the above description that $S(\fS,H(s))$ is a closed
subset of 
$$U_{\Delta(r,\pi)} \times (\bB \setminus{\{0\}})^{s} \times\bB^{d-r-s}
= \bigcup_{\epsilon\to 0} U_{\Delta(r,\pi)}\times U_\ve^s\times\bB^{d-r-s}$$ which is
homeomorphic to $\Delta(r,\pi) \times \rdop_+^s$. The homeomorphism is given by
the restriction of the map 
\[ \Val ~:~ U_{\Delta(r,\pi)} \times (\bB \setminus{\{0\}})^{s} \times \bB^{d-r-s}
\To \Delta(r,\pi) \times \rdop_+^{s}, \quad p \mapsto    
(-\log |x_0(p)|, \dots , -\log |x_{r+s}(p)|) \]
to $S(\fS,H(s))$.  Composing $\Val$ with the inverse homeomorphism
$\Delta(r,\pi)\times\R_+^s\isom S(\fS,H(s))$ yields a map
\[ \tau ~:~ U_{\Delta(r,\pi)} \times (\bB \setminus{\{0\}})^{s} \times \bB^{d-r-s} 
\To S(\fS,H(s)) \]
which is a proper strong deformation retraction. The latter follows from \cite[Theorem 5.2]{berkovich:locallycontractible1}. By construction, we have 
that $\Val=\Val\circ\tau$, i.e.\ that $\Val$ factors through the
retraction to the skeleton.

\paragraph[The skeleton of a building block]
\label{par:skel.building.block}
Next we consider a formal affine open building block $\fU=\Spf(A)$
as in Proposition~\ref{normalization of covering} and let  
\[ \psi ~:~ \fU\To\fS = \fU_{\Delta(r,\pi)} \times \fB^{d-r} 
= \Spf(K^\circ\angles{x_0,\ldots,x_d}/\angles{x_0\cdots x_r-\pi}) \]
be the \'etale map from (a).
Recall that $\fH_\eta$ is the support of $H|_{\fX_\eta}$. Define 
\[ \Val \coloneq \Val\circ\psi ~:~ U\setminus \fH_\eta \To 
\Delta(r,\pi)\times \R_+^s, \]
where $U = \fU_\eta$.  Then we have
\[ \Val(p) =
\big(-\log |\psi^*(x_0)(p)|,\ldots,-\log|\psi^*(x_{r+s})(p)|\big). \]
We define 
the skeleton $S(\fU,H|_\fU)$ as the preimage of $S(\fS,H(s))$ under $\psi$.
This is a closed subset of $U\setminus\fH_\eta$. Using Berkovich's results in \cite[\S5]{berkovich:locallycontractible1} and the
$\ve$-approximation argument from \ref{par:skeleton.standard}, one can show that:
\begin{enumerate}
\item $\psi$ induces a homeomorphism of $S(\fU,H|_\fU)$ onto 
  $S(\fS,H(s))$, so $\Val$ restricts to a homeomorphism
  $S(\fU,H|_\fU)\isom \Delta(r,\pi) \times \rdop_+^s$;
\item composing $\Val$ with the inverse homeomorphism 
  $\Delta(r,\pi)\times\R_+^s\isom S(\fU,H|_\fU)$ yields a proper strong
  deformation retraction
  $ \tau :~ U\setminus\fH_\eta\To S(\fU,H|_\fU)$; 
\item $S(\fU,H|_\fU)$ is intrinsic to $(\fU,H)$ and does not depend on
  the choice of $\psi$.
\end{enumerate}
We wish to emphasize that the above map $\Val$ factors through the
retraction to the skeleton.  Moreover, $\Val$ is essentially
intrinsic to $(\fU,H|_\fU)$ and is independent of
$\psi$ up to reordering the coordinates, as we now prove.  

\begin{lem} \label{lem:comps.are.cartier}
  In the notation above, suppose that $\fU_s$ is not irreducible.
  Let $V$ be an irreducible component of $\fU_s$.
  There exists $i\leq r$ such that the zero set of 
  the restriction of $f \coloneq  \psi^*(x_i) \in A$ to $\fU_s$ is equal to $V$.  Moreover,
  if $f'\in A$ is any other function whose restriction to $\fU_s$ has zero set contained in $V$ and whose restriction to the generic fibre is a unit, then
  $f' = uf^n$ for some $u\in A^\times$ and $n \in \N$.
  We have $\cyc(f)= v(\pi)\, V$.
\end{lem}

\pf
This follows from \cite[Proposition~2.11(c)]{gubler:tropical}.
See Appendix~\ref{Refined intersection theory with Cartier divisors} for
the definition of the Weil divisor $\cyc(f)$ associated to the Cartier divisor $\div(f)$.  \qed 

\smallskip 

Note that  $\fU_s$ is not irreducible if and only if  $r>0$. 
We claim that $\Val: U\setminus\fH_\eta\to \Delta(r,\pi)\times\R_+^s$ is
intrinsic to $(\fU,H)$ up to reordering the coordinates. 
If $r>0$, then Lemma~\ref{lem:comps.are.cartier} shows that 
 the functions $\psi^*(x_i)$ for $i=0,\ldots,r$ are intrinsic to $(\fU,H)$
up to units on $\fU$. Since $H = \psi^*H(s)$ is a Cartier divisior on $\fU$ with distinguished components $H_i$, the functions 
$\psi^*x_i$ for $i=r+1,\ldots,r+s$ are also intrinsic up to units on $\fU$.  This proves that
$\Val: U\setminus\fH_\eta\to \Delta(r,\pi)\times\R_+^s$ is
well-defined up to reordering the coordinates.  
If $r=0$ then $\Delta(0,\pi) = \{0 \}\subset\R$ by our choice $\pi\coloneq 1$ in the definition of a standard pair (see Example \ref{standard pair}), hence there is no additional ambiguity in the definition
of $\Val$. 

\paragraph \label{par:skel.block.2}
We make one final remark about the skeleton $S(\fU,H|_\fU)$.
Using~\cite[Theorem~5.2(iv)]{berkovich:locallycontractible1}, one sees
that for every $x\in U\setminus\fH_\eta$, $\red(\tau(x))$ is equal to the generic point of
the stratum of $D$ containing $\red(x)$.  This implies that if
$\fU'\subset\fU$ is a formal open subset intersecting the minimal stratum
and if $U'$ is the generic fibre of $\fU'$, 
then $S(\fU',H|_{\fU'})$ is equal to $S(\fU,H|_\fU)$.  In particular,
$S(\fU,H|_\fU)\subset U'$.  When $\fU$ is a building block of a formal strictly
semistable pair $(\fX,H)$, this proves that the skeleton
$S(\fU,H|_\fU)\subset X$ depends only on the minimal
stratum and not the choice of building block.

\paragraph[The skeleton of a formal strictly semistable pair] \label{canonical polyhedra}
Let $(\fX,H)$ be a formal strictly semistable pair.  Recall that $X = \fX_\eta$
in our standard notation~\parref{par:formal.ssp.notn}.
In~\parref{par:skel.block.2} we associated to every building block 
$\fU$ of $(\fX,H)$ the skeleton
$S(\fU,H|_\fU)\subset X$, and we showed that this skeleton
depends only on the distinguished stratum $S$ of $\fU$ and not on the
choice of $\fU$.  For clarity we decorate the map $\Val$ associated to
this stratum with the subscript $S$, i.e.\ we have
$\Val_S:\fU_\eta\setminus\fH_\eta\to\Delta(r,\pi)\times\R_+^s$.
We call $\Delta_S\coloneq S(\fU,H|_\fU)\subset X$
the \emph{canonical polyhedron} of $S$.  
The name is justified by the fact
that $\Val_S$ restricts to an identification 
$\Delta_S\isom\Delta(r,\pi) \times \rdop_+^s\subset\rdop_+^{r+s+1}$ which is canonical
up to reordering the coordinates; the range is a polyhedron with a single
maximal bounded face $\Delta(r,\pi)$.
We call $\Val_S^{-1}(\Delta(r,\pi)) \cap \Delta_S$ the \emph{finite part} of $\Delta_S$  and we call $r$ the 
\emph{dimension of the finite part} of $\Delta_S$.  It is equal to the
number of irreducible components of $\fX_s$ containing $S$ minus $1$.
We call $v(\pi)$ the \emph{length} of $\Delta_S$.
Note that $\dim(\Delta_S)$ is equal to
the codimension of $S$ in $\fX_s$.  The canonical polyhedra satisfy the
following properties:
\begin{itemize}
\item[(a)] For a vertical stratum $S$ of $D$, the map $T \mapsto \Delta_T$
  gives a bijective order reversing correspondence between vertical strata
  $T$ of $D$ with $S \subset \overline{T}$ and closed faces of $\Delta_S$.
\item[(b)] For vertical strata $R, S $ of $D$, $\Delta_{R} \cap \Delta_S$
  is the union of all $\Delta_T$ with $T$ ranging over all vertical strata
  of $D$ such that $\overline{T} \supset R \cup S$.
\end{itemize}
We define the \emph{skeleton} of $(\fX,H)$ to be 
\[ S(\fX,H) \coloneq \bigcup_{S\in\str(\fX_s,H)} \Delta_S. \] 
This is a
closed subset of $X \setminus\fH_\eta$ which depends only on the formal strictly
semistable pair $(\fX,H)$.  The above incidence relations endow
$S(\fX,H)$ with the canonical structure of 
a piecewise linear space whose charts are integral $\Gamma$-affine polyhedra and whose transition functions are integral $\Gamma$-affine maps.
More precisely, 
if $\Delta_S\cong\Delta(r,\pi)\times\R_+^s\subset\R_+^{r+s+1}$ and
$\Delta_{S'}\cong\Delta(r',\pi')\times\R_+^{s'}\subset\R_+^{r'+s'+1}$ are
canonical polyhedra associated to strata $S,S'$ with $S'\subset\bar S$, we have 
the following description: If $r>0$, 
then after potentially reordering the coordinates $v_0, \dots, v_{r'+s'}$, the polyhedron $\Delta_S$ is the intersection of
$\Delta_{S'}$ with the linear subspace
$\{v_{r+1}=\cdots=v_{r'}=0, v_{r'+s+1}=\cdots=v_{r'+s'}=0\}$.
In particular, 
we have $v(\pi) = v(\pi')$; as we 
noted after Lemma~\ref{lem:comps.are.cartier}, $v(\pi)$ is intrinsic
to $S$ and $S'$ in this case. The case $r=0$ is trivial as $\Delta_S$ is a vertex of $\Delta_{S'}$.

As a consequence, we see that if two canonical polyhedra
$\Delta_S\cong\Delta(r,\pi)\times\R_+^s$ and
$\Delta_{S'}\cong\Delta(r',\pi')\times\R_+^{s'}$
with $r,r' \geq 1$ are connected by a chain of canonical polyhedra 
$\Delta_T$ whose finite part has positive dimension, then 
$v(\pi) = v(\pi')$.  In this case we say that $\Delta_S,\Delta_{S'}$ are
\emph{connected by finite faces of positive dimension}.  The skeleton
$S(\fX,H)$ canonically decomposes into components which are connected by
finite faces of positive dimension. 

\subparagraph
When $H = 0$ the skeleton $S(\fX,0)$ coincides with the skeleton $S(\fX)$
of the strictly semistable formal scheme $\fX$ in the sense of
Berkovich~\cite{berkovich:locallycontractible1}.  We will use the
notations $S(\fX)$ and $S(\fX,0)$ interchangeably.

\begin{rem} \label{rem:abhyankar}
  Every point $x$ of the skeleton is Abhyankar: that is, 
  \[ \rank_\Z(|\sH(x)^\times|/|K^\times|) + \trdeg(\td\sH(x)/\td K) =
  \dim(X). \]
  Indeed, any point of the skeleton can be interpreted as a monomial
  valuation with respect to some system of local coordinates, and it is
  easy to see that any monomial valuation is Abhyankar.
  Moreover, $x$ induces a valuation on the function field of $X$. As a consequence, 
  the skeleton $S(\fX,H)$ is contained in the complement of every  closed analytic 
  subvariety $Y \neq X$
  (see~\cite[Proposition~9.1.3]{berkovich:analytic_geometry}), hence is a
  ``birational'' feature of $X$.
\end{rem}

\begin{rem} \label{relation between faces}
  A {\it face} of the skeleton $S(\fX,H)$ is the same as a canonical polyhedron. The relative interior 
of a canonical polyhedron is called an {\it open face}. We use the partial ordering $\Delta_T \preccurlyeq \Delta_S$
for canonical polyhedra meaning  $\Delta_T$ is a face of $\Delta_S$. If additionally $\Delta_T \neq \Delta_S$,
 then we write $\Delta_T \prec \Delta_S$. We use this partial ordering also for open faces by applying it to the closures. 
\end{rem}

\paragraph[Retraction to the skeleton] \label{Retraction to the skeleton}
The retractions onto the skeletons of the building blocks
$\tau:\fU_\eta\setminus\fH_\eta \to S(\fU,H|_\fU)$ glue to give a
proper strong deformation retraction 
\begin{equation} \label{eq:retract.to.skel}
  \tau ~:~ X \setminus \fH_\eta \To S(\fX,H). 
\end{equation}
Since $\tau$ is continuous and surjective, $S(\fX,H)$ is connected.
We can use the retraction map to give a different description of the
correspondence between $\str(\fX_s,H)$ and the faces of $S(\fX,H)$. 
In the following {\it orbit-face correspondence}, we use the restriction 
$\red_{H^c}: X \setminus \fH_\eta \rightarrow \fX_s$ of the  reduction map $\red:X \rightarrow \fX_s$.  
Recall also the partial orders $\preccurlyeq$ defined on the open faces of $S(\fX,H)$  and $\leq$ defined on the strata of $D$ (see \ref{stratification}, Remark \ref{relation between faces}).

\begin{prop} \label{stratum-face correspondence} 
  There is a bijective order-reversing 
  correspondence between open faces $\sigma$ of the skeleton $S(\fX,H)$
  and vertical strata $S$ of $D$, given by
  \begin{equation} \label{orbit correspondence}
    S = \red_{H^c}\left(\tau^{-1}(\sigma)\right), \quad \sigma = \tau\left(\red_{H^c}^{-1}(Y)\right),
  \end{equation}
  where  $Y$ is any non-empty subset of $S$.
  We have $\dim(S) + \dim(\sigma) = \dim(X)$.
\end{prop}

\proof This is completely analogous to Proposition 5.7 in
\cite{gubler:canonical_measures}. \qed 

\begin{cor} \label{cor:stratum-genericpt}
  There is a bijective correspondence between open faces $\sigma$ of 
  $S(\fX,H)$ and the generic points $\zeta_S$ of vertical strata $S$ of
  $D$, given by
  \begin{equation} \label{eq:genericpt.correspondence}
    \{\zeta_S\} = \red(\Omega), \quad \sigma = \red\inv(\zeta_S)\cap S(\fX,H),
  \end{equation}
  where $\Omega$ is any non-empty subset of $\sigma$.
\end{cor}

\pf As remarked in~\parref{par:skel.block.2}, for any 
$p\in X \setminus \fH_\eta $, $\red(\tau(p))$ is equal to the generic
point of the stratum containing $\red(p)$.  For
$p\in\relint(\Delta_S)$ we then have $\red(p) = \zeta_S$ by
Proposition~\ref{stratum-face correspondence}.\qed

\begin{rem} \label{Compactification of the skeleton} 
  We can compactify the
  skeleton $S(\fX,H)$ of a formal strictly semistable pair $(\fX,H)$ by taking
  its closure in $X$. We get a compact subset $\hat{S}(\fX,H)$ of
  $X$ whose boundary has the following interpretation: We note that for
  $k=1, \dots , S$, we have canonical formal strictly semistable pairs 
  $(\fH_k, H|_{\fH_k})$ (see Remark~\ref{inductive pairs}) and hence we get
  corresponding skeletons $S(\fH_k, H|_{\fH_k})$. Proceeding inductively, we
  get formal strictly semistable pairs $(\overline{T},H|_{\bar T})$ for every horizontal
  stratum $T$ of $D$. Then it follows from the construction that
  $\hat{S}(\fX,H)$ contains the  disjoint union of $S(\fX,H)$ and of all these skeletons
  $S(\overline{T},H|_{\bar T})$. The retraction map $\tau$ extends to a  map
  $\hat{\tau}:X \rightarrow \hat{S}(\fX,H)$ such that the restriction of $\hat{\tau}$ to any horizontal stratum $T$ is the canonical retraction map to $S(\overline{T},H|_{\bar T})$. 
  It will follow from Theorem \ref{deformation retraction} that $\hat{\tau}$ is a continuous map. We conclude that $\hat{\tau}(X)$ is compact which means that $\hat{S}(\fX,H)$ is indeed the closure of $S$ and that $\hat{\tau}$ is surjective.  
\end{rem}

\begin{thm} \label{deformation retraction}
The map $\hat{\tau}:X \rightarrow \hat{S}(\fX,H)$ is a proper strong deformation retraction.
\end{thm}

\proof By 
\cite[Theorem 5.2]{berkovich:locallycontractible1} and the $\ve$-approximation  from \ref{par:skeleton.standard}, 
we have a map $\Phi: X \times [0,1] \rightarrow  X$. The restriction of  $\Phi$ to $X\setminus \fH_\eta \times [0,1)$ 
gives the homotopy leading to the proper strong deformation retraction $\tau:X \setminus \fH_\eta \rightarrow S(\fX,H)$
from \ref{Retraction to the skeleton}. Moreover, for every horizontal stratum $T$ of $D$ as in Remark \ref{Compactification of the skeleton}, 
the restriction of $\Phi$ to $T$ is the homotopy giving the proper strong deformation retraction of $T$ to $S(\overline{T},H|_{\bar T})$. 
It remains to show $\Phi$ is continuous (properness is then obvious from compactness of $X$). This can be done completely similar as in Berkovich's proof 
of \cite[Theorem 5.2]{berkovich:locallycontractible1} along the following lines:

Let $\fG_m = \Spf\kcirc\langle \Z\rangle$ be the formal affine torus of rank $1$ over  $\kcirc$. 
Recall that we have considered  standard pairs $(\fS\coloneq \fU_{\Delta(r,\pi)} \times \fB^{d-r},  H(s))$ in \ref{standard pair}. 
We first 
check the claim for the slightly restricted standard pair  $(\fS'\coloneq \fU_{\Delta(r,\pi)} \times \fG_m^s \times \fB^{d-r-s},  H(s)|_{\fS'})$. 
Note that any $\fS$ is a finite open union of such $\fS'$ and we may use them as well for building blocks in  the construction of the skeleton which fits better to Berkovich's setting.  

Let $\fG^{(r)}$ be the kernel of the multiplication map $\fG_m^{r} \rightarrow \fG_m$. Then $\fG^{(r)}$ acts canonically on $ \fU_{\Delta(r,\pi)}$ 
and hence $\fG\coloneq \fG^{(r)} \times \fG_m^{d-r}$ acts canonically on $\fS'$. It is clear that $\fG$ is a formal affine torus of rank $d$ and 
the  generic fibre $G$ of $\fG$ is a formal affinoid torus of rank $d$. In Step 2, Berkovich gives a canonical continuous 
map $[0,1] \rightarrow G$, mapping $t \in [0,1)$ to the Shilov boundary point $g_t$ of the closed disk in $T$ with center $1$ and radius $t$.  
Morover, $t=1$ is mapped to the Shilov boundary point of $G$. By~\cite[\S 5.2]{berkovich:analytic_geometry}, the points $g_t$ are peaked 
and the group action induces well-defined points $g_t * x$ on $\fS_\eta'$ for every $x \in \fS_\eta'$. In this way we get a continuous 
homotopy $\fS_\eta' \times [0,1] \rightarrow \fS_\eta'$, given by $(x,t) \mapsto g_t * x$.
 Note that the action by $g_t$ leaves any horizontal stratum $T$ invariant and acts there in the same way as the corresponding peaked point $g_t'$ for 
the formal affinoid torus $G'$ of rank $\dim(T)$ (apply ~\cite[Proposition 5.2.8(ii)]{berkovich:analytic_geometry} with $X=X'=T$ and $\phi$ the projection from $G$ onto $G'$).
By construction, the homotopy agrees with $\Phi$ 
and leads to the strong deformation retraction $\tau:\fS_\eta' \rightarrow S(\fS',H(s)) \cong \Delta(r,\pi)\times\R_+^s$. 

The general case is deduced from the above case with the same arguments as in the proof of \cite[Theorem 5.2]{berkovich:locallycontractible1}. 
The special shape of the building blocks is not used there. \qed
 
\paragraph[Fibres of the retraction] \label{Polytopal domains from the skeleton} 
Let $(\fX,H)$ be a formal strictly semistable pair
and let $S\in\str(\fX_s,H)$ be a zero-dimensional stratum, so 
$S = \{x\}$ for $x\in\fX_s(\td K)$.  Let $\Delta_S\subset S(\fX,H)$ be the
corresponding canonical polyhedron and let $\omega$ be a $\Gamma$-rational
point in the relative interior of $\Delta_S$.  In the proof of Theorem~\ref{Sturmfels--Tevelev multiplicity formula for alterations} it will be
important to understand the analytic subdomain 
$X_\omega\coloneq\tau\inv(\omega)\subset X$.  We will prove
that $X_\omega$ is isomorphic to the affinoid torus
$T = \sM(K\angles{x_1^{\pm 1},\ldots,x_d^{\pm 1}})$, where
$d = \dim(X) = \dim(\Delta_S)$.

Let $\fU\subset\fX$ be a building block with
distinguished stratum $S$.  Let
$\psi:\fU\to\fS=\fU_{\Delta(r,\pi)}\times\fB^{d-r}$ be an \'etale morphism
as in Proposition~\ref{normalization of covering} and let 
$y = \psi(x) = (0,0,\ldots,0)$, so
$\{y\}$ is the minimal stratum of $(\fS,H(d-r))$.
By Proposition~\ref{stratum-face correspondence}
we have $\red(X_\omega) = \red(\tau\inv(\omega)) = \{x\}$, so 
$X_\omega$ is contained in the formal fibre
$\red\inv(x)$.  In particular, $X_\omega\subset\fU_\eta$.  
Since $\psi$ is \'etale at $x\in\fU_s(\td K)$, the induced map on formal
fibres $\red\inv(x)\to\red\inv(y)$ is an isomorphism
by~\cite[Proposition~2.9]{gubler:tropical}.  

Let $v = \Val_S(\omega)\in \Delta(r,\pi)\times\R^{d-r}_+\subset\R^{d+1}$.
The coordinates of $v$ are contained in $\Gamma$, so 
$\Val_S\inv(v)\subset U_{\Delta(r,\pi)}\times(\B\setminus\{0\})^{d-r}
\subset\G_m^{d+1,\an}$ is non-canonically isomorphic to the affinoid
torus $T$ defined above.
Since $\psi$ commutes with the retraction maps 
$\tau:\fU_\eta\setminus\fH_\eta\to\Delta_S$ and 
$\tau:U_{\Delta(r,\pi)}\times(\B\setminus\{0\})^{d-r}\to S(\fS,H(d-r))$,
we have 
\[ X_\omega = \tau\inv(\omega) = \psi\inv(\tau\inv(\psi(\omega))) 
= \psi\inv(\Val_S\inv(v)). \]
On the other hand, as above $\tau\inv(\psi(\omega))$ is contained in the
formal fibre $\red\inv(y)$, so as $\red\inv(x)\to\red\inv(y)$ is an
isomorphism, the same is true for the map
$\psi: X_\omega\isom\Val_S\inv(v)\cong T$.  Hence we have proved:

\begin{prop} \label{prop:Xomega.is.torus}
  Let $(\fX,H)$ be a formal strictly semistable pair, let $S\in\str(\fX_s,H)$ be a
  zero-dimensional stratum, and let $\omega\in\Delta_S$ be a
  $\Gamma$-rational point contained in the relative interior of
  $\Delta_S$.  Then $X_\omega\coloneq\tau\inv(\omega)$ is isomorphic to
  the affinoid torus $T = \sM(K\angles{x_1^{\pm 1},\ldots,x_d^{\pm 1}})$,
  where $d = \dim(X) = \dim(\Delta_S)$.
\end{prop}

\smallskip
Let $Y = \sM(A)$ be an affinoid space and let 
$A^\circ\subset A$ be the subring of power-bounded elements.
Recall from~\parref{par:notation} that
the \emph{canonical model} of $Y$ is the $K^\circ$-formal scheme
$\Spf(A^\circ)$; this is an affine admissible formal scheme when $A$
is reduced by~\cite[Theorem~3.17]{bpr:trop_curves}.

\begin{cor} \label{torus corollary}
  With the notation in Proposition~\ref{prop:Xomega.is.torus}, the
  analytic subdomain $X_\omega$ is strictly affinoid, and its
  canonical model $\fX_\omega$ is isomorphic to a
  formal torus over $\kcirc$ of rank $d$.
\end{cor}

\smallskip

At this point we are able to formulate a precise statement about which
vertical components of $D$ are Cartier, and where.  
See
\ref{associated Weil divisor}   for
the definition of  the Weil divisor $\cyc(C)$ associated to the Cartier divisor $C$.

\begin{prop} \label{prop:comps.are.cartier}
  Let $V$ be an irreducible component of $\fX_s$ and let
  $v$ be the
  vertex of $S(\fX,H)$ associated to the
  open stratum in $V$.
  \begin{enumerate}
  
    \item If $\fX_s=V$, then   any  vertical effective Cartier divisor on $\fX$ is equal to $\div(\lambda)$ for some non-zero $ \lambda \in \kcirc$. 

    \item If $\fX_s \neq V$ and  if all  canonical polyhedra containing $v$ with positive-dimensional finite part have the same length $v(\pi)$,   then there is a unique effective Cartier divisor $C$ on $\fX$ with $\cyc(C) = v(\pi)\, V$, and  any  effective Cartier divisor  with support contained in $V$ is equal to $n \,C$ for a unique $n \in \N$.

    \item If $v$ is adjacent to canonical polyhedra $\Delta_{S_1},\Delta_{S_2}$
    with positive-dimensional finite part and
    whose lengths $v(\pi_1),v(\pi_2)$ are not
    commensurable, then $V$ is not the support of 
    a Cartier divisor on $\fX$.
  \end{enumerate}
\end{prop}

\pf Note that~(1) follows from the fact that  a vertical Cartier divisor on an admissible formal scheme over $\kcirc$ with reduced special fibre is 
uniquely determined by the associated Weil divisor  (see Proposition \ref{vertical Cartier divisors}).

To prove~(2), we assume that $\fX_s \neq V$. Since $X$ is connected, continuity and surjectivity of the retraction map $\tau$ yield that $S(\fX,H)$ is connected. 
We conclude from the stratum-face correspondence in Proposition \ref{stratum-face correspondence} 
that there is at least one canonical polyhedron  of $S(\fX,H)$
containing $v$
with positive-dimensional finite part. By assumption, all such canonical polyhedra have the same length $v(\pi)$ for some
non-zero $\pi \in K^{\circ \circ}$. 
Choose a cover of $\fX$ by building blocks $\fU_i$.  If
$(\fU_i)_s\cap V = \emptyset$ then a local equation for $C$ on $\fU_i$ is
$f_i = 1$.  If $V$ is the distinguished stratum of $\fU_i$ then a local equation
for $C$ on $\fU_i$ is $f_i = \pi$.  Otherwise the special fibre of $\fU_i$ is
not irreducible; we choose the function $f_i = f$ of 
Lemma~\ref{lem:comps.are.cartier} as the local equation for $C$ on
$\fU_i$.  As the Weil divisor associated to $\div(f_i)$ is $v(\pi)\,(V\cap(\fU_i)_s)$ for all $i$, these
indeed define an effective Cartier divisor $C$. 
This and uniqueness follows again from Proposition \ref{vertical Cartier divisors}.

To prove the last statement in~(2), let $C'$ be an  effective Cartier divisor  with support contained in $V$. Then the associated Weil divisor is equal to $v(\lambda) \, V$ for  $\lambda \in K^\circ$. We may apply Lemma \ref{lem:comps.are.cartier} to a buidling block $\fU_i$ whose special fibre is not irreducible. The latter is equivalent to the property that the finite part of the canonical polyhedron of the distinguished stratum is positive-dimensional. This proves that $v(\lambda)= n \, v(\pi)$ for some non-zero $n \in \N$. Then  Lemma \ref{lem:comps.are.cartier} again and part~(1) applied to the buidling blocks with irreducible special fibre prove that $C'= n \, C$. Since $n$ is the multiplicity of $C'$ in $V$, it is unique.

In the situation of~(3), we may shrink $\fX$ to assume that it is covered
by two building blocks $\fU_1,\fU_2$ with distinguished strata $S_1$ and
$S_2$.  It follows from Lemma~\ref{lem:comps.are.cartier} that any Cartier
divisor $C_i$ on $\fU_i$ with support equal to $(\fU_i)_s\cap V$  has
$\cyc(C_i) = n_iv(\pi_i)\,(V\cap(\fU_i)_s)$ for $n_i\in\Z\setminus\{0\}$.
Since $n_1v(\pi_1)\neq n_2v(\pi_2)$ for $n_1,n_2\neq 0$ there does not exist
a Cartier divisor $C$ on $\fX$ such that $C|_{\fU_i} = C_i$ for
$i=1,2$.\qed 

\section{Functoriality}  \label{section:functoriality}

In this section, $(\sX,H)$ is a strictly semistable pair.
As always we use the associated notation introduced in~\parref{par:formal.ssp.notn};
in particular, $X\coloneq \sX_\eta$ is the generic fibre.  Let $d\coloneq
\dim(X)$.  Consider a non-zero rational function $f$ on $X$.  This induces a meromorphic function on 
$\Xan$. The goal of this section is to show that the restriction of $-\log|f|$ to the skeleton $S(\sX,H)$
is an everywhere-defined piecewise linear function. If the  support of $\div(f)$ is contained in the boundary $\supp(H)_\eta$, then  the restriction to any canonical polyhedron of $S(\sX,H)$ is integral $\Gamma$-affine for the value group $\Gamma=v(K^\times)$. This basic result will be used in later sections. As a further consequence, we will show functoriality of the retraction to the skeleton.
We will deduce piecewise linearity of the restriction of $-\log|f|$ to the skeleton $S(\sX,H)$ without boundary assumptions. Note 
that in this generality, the restriction is not necessarily integral $\Gamma$-affine on canonical polyhedra as above.

\paragraph \label{par:deco1}
Let $f$ be a non-zero rational function on $X$ such that the support of
$\div(f)$ is contained in the generic fibre of $\supp(H)$. 
As in ~\parref{notn:ssp.notation}, let $D_1=V_1, \dots , D_R=V_R$ be the irreducible components of the special fibre $\sX_s$,
let $D_{R+1}=\sH_1, \dots , D_{R+S}=\sH_S$ be the prime components of the 
(horizontal) Weil divisor $\sH$ on $\sX$ associated to $H$ and let $D=D_1 + \dots +D_{R+S}$. 
We have
\begin{equation}
  \label{deco1}
  \cyc(f) = \sum_{i=1}^{R+S} \ord(f,D_i)\,D_i
  = \sum_{i=1}^R \ord(f,V_i)\,V_i + \sum_{j=1}^S \ord(f,\sH_j)\,\sH_j. 
\end{equation}
for the associated Weil divisor $\cyc(f)$ on $\sX$. Here $\ord(f,\sH_j)\in\Z$ is the usual order of vanishing of $f$ along 
$(\sH_j)_\eta$ in $X$ and $\ord(f,V_i) = -\log|f(\xi_i)|$, where $\xi_i\in X^\an$ is
the unique point reducing to the generic point of $V_i$; see Appendix~\ref{Refined intersection theory with Cartier divisors}.

The next result generalizes parts~(1) and~(2) of~\cite[Theorem~5.15]{bpr:analytic_curves}.

\begin{prop} \label{piecewise linear} Let $f$ be a non-zero rational function on $X$ such that the support of
$\div(f)$ is contained in the generic fibre of $\sH=\supp(H)$. 
  We consider $F = -\log|f|$ as a function 
  $F:U^\an \coloneq X^\an\setminus\sH_\eta^\an\to\R$. 
  Then 
$F$ factors through the
  retraction map $\tau:U^\an\to S(\sX,H)$. 
  Moreover, the restriction of $F$ to $S(\sX,H)$ is an integral
  $\Gamma$-affine
  function on each canonical polyhedron.
\end{prop}

\proof
Let $\fX=\hat\sX$ and let $\fU\subset\fX$ be a building block with
distinguished stratum $S$.  Choose an \'etale morphism
$\psi:\fU\to\fS=\fU_{\Delta(r,\pi)}\times\fB^{d-r} = 
\Spf(K^\circ\angles{x_0, \dots, x_d} / \angles{x_0 \dots x_r - \pi})$ 
as in Proposition~\ref{normalization of covering}.  This means that 
the restriction of the formal completion of $(\sX,H)$ to $\fU$ is equal
to the pull-back of the standard pair $(\fS,H(s))$ with respect to $\psi$ 
for some $s \in \{0, \dots ,d-r\}$. 
Since the support of $\div(f)$ on the generic fibre is contained in the generic fibre of the horizonal divisor, there exist integers $n_{r+1}, \ldots, n_{r+s}$ such that $f|_\fU / \prod_{i = r+1}^{r+s} \psi^{\ast}(x_i)^{n_i}$ is a unit on $\fU_\eta$. 
By~\cite[Proposition~2.11]{gubler:tropical}, we find
\begin{equation} \label{deco2}
f|_\fU=\lambda u\, \psi^*(x_0)^{n_0} \cdots \psi^*(x_{r+s})^{n_{r+s}},
\end{equation}
with $\lambda \in K^\times$, $u \in \Ocal(\fU)^\times$ and 
$n_0, \ldots , n_{r} \in \zdop$.
It follows that for $p\in\fU_\eta$ we have
\[ F(p) = v(\lambda) - \sum_{j=0}^{r+s} n_j\log|\psi^*x_j(p)|.\]
By definition $\Val_S:\fU_\eta\setminus\sH_\eta^\an\to\R^{r+s+1}$ 
is given by 
\[ \Val_S(p) = \big(-\log|\psi^*x_0(p)|,\ldots,-\log|\psi^*x_{r+s}(p)|\big), \]
so the restriction of $F$ to $\fU_\eta \setminus \sH_\eta^\an$ is the composition of $\Val_S$ followed by
$$(v_0,\ldots,v_{r+s})\mapsto v(\lambda) + \sum_{j=0}^{r+s} n_j v_j.$$
We conclude that $F$ factors through the retraction $\tau$. 
We also see that the restriction of $F$ to $\Delta_S$ is integral $\Gamma$-affine. \qed

\begin{rem}   \label{order and representation}
The function $F$ is completely determined by the Weil divisor $\cyc(f)$ in \eqref{deco1}. 
From the proof of Proposition \ref{piecewise linear}, we deduce the following explicit formula 
for $F$ in terms of the multiplicities $(\ord(f,D_i))_{i=1,\dots,R+S}$: 
For $j=0, \dots, s$, the zero set of $\psi^*(x_j)$ on the building block $\fU$ is equal to $V_{i_j}$ for a unique 
$i_j \in \{1, \dots , R\}$. 
For $j=r+1, \dots ,r+s$, the coordinate $\psi^*(x_j)$ on  $\fU$ is a local equation for a unique $\hat{H}_{i_j}$ with $i_j \in \{1, \dots, S\}$. 
By Lemma~\ref{lem:comps.are.cartier} and \eqref{deco1}, we have the relations  $\ord(f,V_{i_j}) = v(\lambda) + n_jv(\pi)$ for
$j\in \{0, \dots,  r\}$ and $\ord(f,\sH_{i_j}) = n_j$ for $j \in \{r+1, \dots, r+s\}$.
For $p \in \fU_\eta$, we get
\[\begin{split}
  F(p) &= v(\lambda) - \sum_{j=0}^{r+s} n_j\log|\psi^*x_j(p)| \\
  &= v(\lambda) - \sum_{j=0}^{r} \frac{\ord(f,V_{i_j})-v(\lambda)}{v(\pi)}\,\log|\psi^*x_j(p)|
  - \sum_{j=r+1}^s \ord(f,\sH_{i_j})\,\log|\psi^*x_j(p)| \\
  &= - \frac 1{v(\pi)}\sum_{j=0}^{r} \ord(f,V_{i_j})\,\log|\psi^*x_j(p)|
  - \sum_{j=r+1}^s \ord(f,\sH_{i_j})\,\log|\psi^*x_j(p)|,
\end{split}\]
where the last equality holds because 
$-\sum_{j=0}^r\log|\psi^*x_j(p)| = v(\pi)$.  
We conclude that the restriction of $F$ to $\fU_\eta \setminus \sH_\eta^\an$ is the composition of $\Val_S:\fU_\eta \setminus \sH_\eta^\an \to \R^{r+s+1}$ followed by the linear form
\begin{equation} \label{ord representation}
 (v_0,\ldots,v_{r+s})\mapsto\sum_{j=0}^r \frac{\ord(f,V_{i_j})}{v(\pi)}\,v_j 
+ \sum_{j=r+1}^s \ord(f,\sH_{i_j})\,v_j. 
\end{equation}
This is the desired formula in terms of the multiplicities. 

By the stratum--face correspondence in Proposition 
\ref{stratum-face correspondence}, a vertex $u$ of the canonical polyhedron $\Delta_S$ corresponds to a vertical 
component $V_i$ of $D$. Then we deduce from \eqref{ord representation} that 
$$F(u)=\ord(f,V_i).$$
Alternatively, by Corollary~\ref{cor:stratum-genericpt}, $u$ is the
unique point of $S(\sX)$ reducing to the generic point of $V_i$, so
$\ord(f,V_i) = -\log|f(u)| = F(u)$ by definition.
\end{rem}

\begin{rem} \label{generalization to formal pairs}
The above results hold also for a non-zero meromorphic 
function $f$ on a formal strictly semistable pair $(\fX,H)$ assuming 
that the generic fibre of the support of $\div(f)$ is contained in the 
generic fibre of $\supp(H)$. The same proof applies without 
any change. 
\end{rem}

\smallskip
Next, we will show {\bf functoriality} of the retraction to the skeleton.
This holds for formal strictly semistable pairs as well, but we will restrict 
our attention to the algebraic case.

\begin{prop} \label{Functoriality}
Let $(\sX,H), (\sX',H')$ be  strictly semistable pairs and let $\varphi: \Xcal' \rightarrow \Xcal$ be a morphism with generic 
fibre $\varphi_\eta: X' \rightarrow X$. We assume that $\varphi_\eta^{-1}(\supp(H)_\eta) \subset \supp(H')_\eta$. Then there is a unique map $\varphi_{\rm aff}:S(\sX',H') \rightarrow S(\sX,H)$ with 
$$\varphi_{\rm aff} \circ  \tau' = \tau \circ \varphi_\eta$$ on $X' = \sX_\eta'$, where $\tau$ (resp. $\tau'$) is the retraction 
to the skeleton $S(\sX,H)$ (resp. $S(\sX',H')$) from \ref{Retraction to the skeleton}. The map $\varphi_{\rm aff}$ is continuous. For any  canonical polyhedron $\Delta_{S'}$ of $S(\sX',H')$, there is a canonical polyhedron $\Delta_S$ of $S(\sX,H)$ with $\varphi_{\rm aff}(\Delta_{S'}) \subset \Delta_S$. Moreover, the induced map $\Delta_{S'} \rightarrow \Delta_S$ is 
integral $\Gamma$-affine.
\end{prop}

\proof Put $\fX=\hat\sX$ and consider a covering by building blocks $\fU_i$
as in Proposition~\ref{normalization of covering}.  
Now refine the covering $\varphi_s^{-1}(\fU_{i,s})$ of $\fX'_s$ according to 
Proposition~\ref{normalization of covering}.
In this way we get a covering of $\fX'$ by building blocks $\fU'_{i,k}$ such that each $\fU'_{i,k}$ maps to a building block $\fU_i$ of $\fX$ via the formal completion $\widehat{\varphi}$ of $\varphi$. 
Let $\fU'$ be one of those building blocks of $\fX'$ with associated stratum $S'$ such that $\widehat{\varphi}$ maps $\fU'$  to a building block $\fU$ of $\fX$ with associated stratum $S$. 
Choose  an \'etale morphism
$\psi:\fU\to\fS=\fU_{\Delta(r,\pi)}\times\fB^{d-r} = 
\Spf(K^\circ\angles{x_0, \dots, x_d} / \angles{x_0 \dots x_r - \pi})$ 
as in Proposition~\ref{normalization of covering}~(a).

Consider the analytic functions $f_i = \varphi_\eta^{\ast} \psi_\eta^\ast x_i $ on $\fU'_\eta$ for $i = 0, \ldots, r+s$.  Since $\varphi_\eta^{-1}(\supp(H_\eta)) \subset \supp(H')_\eta$,
all $f_i$ are units on 
$\fU'_\eta \backslash (\sH'_\eta)^{\rm an}$. 
Choose an \'etale morphism
$\psi':\fU'\to\fS'=\fU_{\Delta(r',\pi')}\times\fB^{d'-r'} = 
\Spf(K^\circ\angles{x'_0, \dots, x'_{d'}} / \angles{x'_0 \dots x'_{r'} - \pi'})$ 
as in Proposition~\ref{normalization of covering} (a).
Now we argue simililarly as in the proof of Proposition ~\ref{piecewise linear}. It follows from~\cite[Proposition~2.11]{gubler:tropical} that 
\begin{equation}
f_i =\lambda u\, \psi'^{*}(x'_0)^{n_0} \cdots \psi'^*(x'_{r'+s'})^{n_{r'+s'}},
\end{equation}
with $\lambda \in K^\times$, $u \in \Ocal(\fU')^\times$ and 
$n_0, \ldots , n_{r'+s'} \in \zdop$. 
Hence the map $-\log |f_i|$ factors through an
integral $\Gamma$-affine map on $\Delta_{S'} = S(\fU',H'|_{\fU'})$.

Note that the map $(\tau \circ \varphi_\eta)|_{\fU'_\eta \backslash (\sH'_\eta)^{\rm an}}$ is given by $(-\log|f_0|,\ldots,  -\log| f_{r+s}|)$.
Hence there exists a unique map $\varphi_{\aff}: \Delta_{S'} \to \Delta_S $ satisfying 
$\varphi_{\aff} \circ \tau' = \tau \circ \varphi_\eta$. By construction, it is integral $\Gamma$-affine. As these maps fit together on the intersection 
of building blocks, we get a well-defined map  $\varphi_{\rm aff}:S(\sX',H') \rightarrow S(\sX,H)$ with the required properties. 
Uniqueness is obvious from surjectivity of $\tau'$. \qed

\smallskip

\begin{prop} \label{lem:piecewise.affine}
  Let $(\sX,H)$ be a strictly semistable pair and 
  let $f$ be a non-zero rational function on $X = \sX_\eta$.  Then $S(\sX,H)$ can be
  covered by finitely many integral $\Gamma$-affine polyhedra $\Delta$ such that 
  $-\log|f|$ restricts to an integral $\Gamma$-affine function on  $\Delta$. 
\end{prop}

\proof Locally, $\sX$ is given by an open subset $\sU$ with an \'etale morphism 
    \begin{equation} 
    \psi ~:~ {\Ucal} \longrightarrow {\Scal} \coloneq
    \Spec\big( K^\circ[ x_0, \dots, x_d ] / \angles{x_0 \cdots x_r - \pi} \big) 
  \end{equation}
   and $H |_\sU =\div(\psi^*(x_{r+1} \dots x_{r+s}))$. 
Then the formal completion $\fU : = \hat{\sU}$ with respect to a non-zero $\pi \in K$ with $|\pi| <1$ together with the Cartier divisor $H_\fU$ induced by $\psi^*(x_{r+1} \dots x_{r+s})$ is the pull-back of the standard pair $(\fS,H(s))$ as in Example  \ref{standard pair}. We may assume that $\fU$ is a building block with distinguished stratum $S$. Then we get an \'etale morphism
 \begin{equation} 
    \hat\psi ~:~ {\fU} \longrightarrow {\fS} \coloneq
    \Spf\big( K^\circ \langle x_0, \dots, x_d \rangle / \angles{x_0 \cdots x_r - \pi} \big) 
  \end{equation}
of affine formal schemes. The corresponding face $\Delta_S$ of $S(\sX,H)$ is mapped homeomorphically onto $S(\fS,H(s)) = \Delta(r,\pi) \times \R_+^s$ by $\psi^{\rm an}$. In the following, we will use this to identify $\Delta_S$ with $\Delta(r,\pi) \times \R_+^s$.

We have to show that $\Delta_S$ can be covered by finitely many integral $\Gamma$-affine polyhedra $\Delta$ such that $(-\log|f|)|_\Delta$ is integral $\Gamma$-affine. Let $H'$ be the Cartier divisor of $\fU$ given by $\hat\psi^*(x_{r+1} \dots x_{d})$ and let $\fH'$ be the associated Weil divisor on $\fU$. Then $(\fU,H')$ is a formal affine strictly semistable pair with skeleton $\Delta(r,\pi) \times \R_+^{d-r}$ containing $\Delta_S = \Delta(r,\pi) \times \R_+^s \times \{0\}$ as a closed face. Moreover, the retraction
\[ \tau_d ~:~ \fU_\eta \setminus \fH_\eta' \To \Delta(r,\pi) \times
\R_+^{d-r} \subset \R_+^{r+1} \times \R_+^{d-r}=\R_+^{d+1} \]
composed with the projection $\R_+^{d+1}=\R_+^{r+s+1} \times \R_+^{d-r-s} \rightarrow \R_+^{r+s+1}$ is the restriction of the retraction $\tau: \Xan \setminus \sH_\eta \rightarrow S(\sX,H)$ to $\fU_\eta \setminus \fH_\eta'$. 

Let $Y$ be the support of the Cartier divisor $\div(f)$ on $X$. We claim that $\tau_d(Y \cap \fU_\eta \setminus \fH')$ is contained in the support of an integral $\Gamma$-affine polyhedral complex in $\Delta(r,\pi) \times \R_+^{d-r}$ of dimension at most $d-1$. This follows from the fact that $\tau_d(Y^{\rm an} \cap \fU_\eta \setminus \fH_\eta')$ is contained in the tropicalization of the closure of $\psi_\eta(Y \cap \sU_\eta)$ in the torus $\Spec(K[x_0^{\pm 1}, \dots, x_d^{\pm 1}])$. We conclude that there is an integral $\Gamma$-affine polyhedral subdivision $\Sigma$ of $\Delta(r,\pi) \times \R_+^{d-r}$ such that $\tau_d(Y^{\rm an} \cap \fU_\eta \setminus \fH_\eta')$ is contained in lower dimensional faces of $\Sigma$. Since $\Delta_S$ is a face of $S(\fU,H') = \Delta(r,\pi) \times \R_+^{d-r}$, it is enough to show that the restriction of $-\log|f|$ to any maximal face $\sigma$ of $\Sigma$ is integral $\Gamma$-affine. By  density of the value group $\Gamma$ and by continuity, it is enough to prove that $(-\log|f|)|_\Delta$ is integral $\Gamma$-affine for any simplex $\Delta$ contained in the interior of $\sigma$. Here, $\Delta$ is assumed to be unimodular to $\Delta(d,\pi')$ for some non-zero $\pi' \in K$ with $|\pi'| \leq |\pi|$. It follows that $\tau_d^{-1}(\Delta)$ is the generic fibre of a strictly semistable formal open subscheme $\fU_\Delta$ of $\fU$. By construction, we have $S(\fU_\Delta,H'|_{\fU_\Delta}) = \Delta$. Since $\tau_d(Y)$ does not meet $\Delta$, we conclude that $f$ is invertible on $(\fU_\Delta)_\eta$. By Proposition \ref{piecewise linear} and Remark \ref{generalization to formal pairs}, the function $(-\log|f|)|_\Delta$ factorizes through the retraction to $S(\fU_\Delta, H'|_{\fU_\Delta})=\Delta$ and the restriction of $-\log|f|$ to $\Delta$ is an integral $\Gamma$-affine function. \qed

\section{The slope formula for skeletons} \label{The slope formula for skeletons}

In this section we assume that $(\Xcal,H)$ is a strictly semistable pair
with generic fibre $X$ of dimension $d\coloneq\dim(X)$.
We will give a generalization of the slope
formula in~\cite[Theorem~5.15]{bpr:analytic_curves} to this higher dimensional
situation. In this section, we will use the divisoral intersection theory on $\sX$ 
recalled in  Appendix~\ref{Refined intersection theory with Cartier divisors}.

Let $f$ be a non-zero rational function on $X$ such that the support of
$\div(f)$ is contained in the generic fibre of $\supp(H)$. 
As in ~\parref{par:deco1}, let $D_1=V_1, \dots , D_R=V_R$ be the irreducible components of the special fibre $\sX_s$,
let $D_{R+1}=\sH_1, \dots , D_{R+S}=\sH_S$ be the prime components of the 
(horizontal) Weil divisor $\sH$ on $\sX$ associated to $H$ and let $D=D_1 + \dots +D_{R+S}$. 
We have 
\begin{equation} \label{decomposition of Weil divisor to f}
  \cyc(f) = \sum_{i=1}^{R+S} \ord(f,D_i)\,D_i
  = \sum_{i=1}^R \ord(f,V_i)\,V_i + \sum_{j=1}^S \ord(f,\sH_j)\,\sH_j. 
\end{equation}
for the associated Weil divisor $\cyc(f)$ on $\sX$. 
We recall from Proposition \ref{piecewise linear}  that the restriction of
$-\log|f|$ to the skeleton $S(\Xcal,H)$ is a continuous function which is
integral $\Gamma$-affine on any canonical polyhedron. The goal of this section is to study the slopes of this piecewise linear function.

\paragraph[The retraction of the divisor to the compactified skeleton] 
\label{tropical divisor at infinity} 
We recall from
Remark~\ref{Compactification of the skeleton} that we can compactify the
skeleton $S(\Xcal,H)$ to a subset $\hat{S}(\Xcal,H)$ of $\Xan$ whose
boundary is the disjoint union of skeletons $S(\overline{T},H|_{\bar T})$ for
the strictly semistable pairs $(\overline{T},H|_{\bar T})$ associated to the
horizontal strata $T$ of $D$. 
Any $(d-1)$-dimensional canonical polyhedron $\Delta_S'$ of the boundary
is contained in a unique skeleton $S(\sH_k, H|_{\sH_k})$ for 
some $k=1, \dots , S$; for such $\Delta_S'$ we set
$m(\hat\tau( f),\Delta_S')\coloneq\ord(f,\sH_k)$.
We define the {\it retraction $\hat{\tau}(f)$ of $\cyc(f)$ to $S(\sX,H)$}  as the
formal sum 
\[ \hat{\tau}(f)\coloneq \sum_{\Delta_S'} m(\hat\tau( f),\Delta_S')\, \Delta_S', \]
where the sum ranges over all $(d-1)$-dimensional canonical polyhedra
$\Delta_S'$ in the boundary of $\hat S(\sX,H)$. 
The support of $\hat{\tau}(f)$ is either of pure dimension $d-1$ or empty.

Let $\Delta_S'$ be a $(d-1)$-dimensional canonical polyhedron in the
boundary of $\hat S(\sX,H)$.  Then $\Delta_S'$ corresponds to a
zero-dimensional vertical stratum $S$ of $D|_{\sH_k}$ for some $k$.
Hence $S$ is a vertical stratum of $D$ contained in
$\sH_k$, so it is a component of the intersection of
$(\sH_k)_s$ with the closure of a unique one-dimensional vertical stratum
$T$ of $D$.  We have $\Delta_S = \Delta_T\times\R_+$, with the direction
$(0,1)$ corresponding to the divisor $H_k$, and the canonical polyhedron
$\Delta'_S$ is naturally identified with $\Delta_{T}\times\{\infty\}$. 
By Remark \ref{order and representation} we can recover the multiplicity
$m(\hat\tau( f),\Delta_S')$ of $\Delta_S'$ in $\hat\tau( f)$ as the slope
of the restriction of $F = -\log|f|$ to $\Delta_S = \Delta_T\times\R_+$ in
the direction $(0,1)$.

\paragraph \label{tropical divisor} 
We define the {\it retraction $\tau(f)$ of $\cyc(f)$ to $S(\sX)$}
 as the push-forward of
$\hat{\tau}( f)$ with respect to the canonical retraction 
$ S(\Xcal, H)\to S(\Xcal)$.  The push-forward is defined as
follows.  If $\Delta_S' = \Delta_{T}\times\{\infty\}$ is a
dimension-$(d-1)$ canonical polyhedron of the boundary of 
$\hat S(\sX,H)$ as in~\parref{tropical divisor at infinity}, we define its
push-forward to be $\Delta_{T}$ if $\Delta_{T}\subset S(\sX)$, i.e.\ if
$\Delta_{T}$ is bounded, and we define it to be zero otherwise.
The corner locus $\tau(f)$ is then the formal sum of the push-forwards of
the $(d-1)$-dimensional canonical polyhedra of $\hat\tau(f)$.  We write
\[ \tau( f) = \sum_{\Delta_T} m(\tau( f),\Delta_T)\,\Delta_T, \]
where the sum ranges over all $(d-1)$-dimensional canonical polyhedra of
$S(\sX)$.
Note that, assuming $\Delta_{T}$ is bounded of dimension $d-1$, the coefficient 
$m(\tau(f),\Delta_T)$ of
$\Delta_{T}$ in $\tau(f)$ is given by
$$m(\tau(f),\Delta_T)= \sum_k \ord(f, \sH_k) \cdot \#(\bar T\cap(\sH_k)_s),$$
where $k$ ranges over all numbers in $\{1, \dots, S\}$ such that $\sH_k$ does not contain $T$.

\paragraph[Cartwright's $\alpha$-numbers] \label{cartwright}
For the moment, we assume that we are in the following situation often occurring in number theory: 
Let $R$ be a complete discrete valuation ring 
with uniformizer $\pi$. We assume that $R$ is a subring of $K^\circ$ and that the discrete valuation of $R$ extends to our given valuation $v$ on the 
algebraically closed field $K$. 
Suppose that $\Xcal$ is the base change of a strictly semistable scheme over $R$ in the sense of de Jong (see \ref{comments for strictly semistable}).  
In such a situation, the strictly semistable scheme $\Xcal$ can be  covered by open subsets $\sU$ which admit an \'etale morphism 
 \begin{equation*}
    \psi ~:~ \fU \To 
    \Spf\big( K^\circ\angles{x_0, \dots, x_d} / \angles{x_0 \dots x_r - \pi} \big) .
  \end{equation*}
Let $u$ be a vertex of $S(\sX,H)$ corresponding to the irreducible component $V_u$ of $\Xcal_s$ by the stratum--face correspondence in Proposition \ref{stratum-face correspondence}. 
Since we have the same $\pi$ for every chart $\fU$ and since we assume $v(\pi)=1$, there is a unique Cartier divisor $C_u$ on $\sX$ with $\cyc(C_u)= V_u$ (see Proposition \ref{prop:comps.are.cartier}). 

Let $T$ be a one-dimensional vertical stratum of $D$, so that the corresponding canonical polyhedron $\Delta_T$ has dimension $d-1$. 
Then we have the intersection number $C_u\cdot\bar T$ and we set 
$$\alpha(u,\Delta_T)\coloneq -C_u\cdot\bar T .$$
 Since everything 
is defined over $R$, this is a usual intersection number on a proper regular 
noetherian integral scheme over $R$ and hence $\alpha(u,\Delta_T) \in \Z$. 

Dustin Cartwright uses the numbers $\alpha(u,\Delta_T)$ to endow the compact skeleton
$S(\sX)$ with the structure of a tropical complex in the sense
of~\cite[Definition~1.1]{cartwright:tropical_complexes}.  Note that he imposes an 
additional local Hodge condition which plays no role here.  

\paragraph \label{par:slopes2}
Our next goal is to generalize Cartwright's $\alpha$-numbers to any strictly
semistable pair $(\Xcal,H)$ without additional assumptions.  The resulting
objects might be called weak tropical complexes. We have to deal with the
problem that the  irreducible component $V_u$ of $\Xcal_s$ corresponding to the
vertex $u$ of $S(\Xcal,H)$ is not necessarily the  
support of a Cartier divisor (see Proposition \ref{prop:comps.are.cartier}). 
As above,  a one-dimensional vertical stratum $T$ of $D$ corresponds to  a $(d-1)$-dimensional canonical polyhedron $\Delta_T$ of
Since $\sX$ is a proper
scheme, the curve $\bar T$ is projective.  For every point
$x\in\bar T(\td K)$ there is a neighbourhood $\fU_x$ of $x$
in $\fX=\hat\sX$ which is a building block
with distinguished stratum $T$ or $\{x\}$.  We let 
$\fU = \bigcup_{x\in\bar T(\td K)} \fU_x$.  This is a formal open subset of 
$\fX$ containing $\bar T$. 
Note that $\fU$ is a strictly semistable formal scheme and the Cartier divisor $H$ on $\sX$ 
induces a formal strictly semistable pair $(\fU,\hat{H}|_\fU)$ (use Proposition \ref{prop:formal.alg.sss.pair}).
By construction, we have 
$$\fU_\eta\cap S(\sX,H) = S(\fU,\hat{H}|_\fU)=\bigcup_{\Delta_S \succcurlyeq \Delta_T} \Delta_S,$$
where $\Delta_S$ ranges over all canonical polyhedra of $S(\sX,H)$ which contain $\Delta_T$. 
In particular, the vertices $u$ of $S(\sX,H)$ contained in $\fU_\eta$
correspond to the irreducible components $V_u$ of $\sX_s$ intersecting $\bar T$.

Let $\Delta_S$ be a $d$-dimensional canonical polyhedron of $S(\sX,H)$ containing
$\Delta_T$.  Then $S$ is a component of the intersection of $\bar T$ with
a vertical component $V_u$ if and only if the finite part of
$\Delta_S$ is strictly larger than the finite part of $\Delta_T$, in which
case we say that $\Delta_S$ 
\emph{extends $\Delta_T$ in a bounded direction}.
Otherwise $S$ is a component of the intersection of $\bar T$ with a
horizontal component of $D$; in this case we say that $\Delta_S$
\emph{extends $\Delta_T$ in an unbounded direction}.
We define $\deg_b(\Delta_T)$ (resp.\ $\deg_u(\Delta_T)$) to be the number
of canonical polyhedra $\Delta_S$ extending $\Delta_T$ in a bounded
(resp.\ unbounded) direction. 

For each vertex $u\in\Delta_T$ we define an integer $\alpha(u,\Delta_T)$ as
follows. 

\begin{enumerate}
\item Suppose that $\Delta_T$ has positive-dimensional finite part, i.e.\
  that $T$ lies on at least two irreducible components of $\sX_s$, and let
  $v(\pi)$ be the length of $\Delta_T$.
  Then the same is true of any $\Delta_S \succcurlyeq \Delta_T$, so by
  Proposition~\ref{prop:comps.are.cartier}, for every vertex $u$ of
  $\Delta_T$, there is a unique effective Cartier
  divisor $C_u$ on $\fU$ with $\cyc(C_u) = v(\pi)\,(V_u\cap\fU_s)$.
  For every such vertex $u$ we let
  $-\alpha(u,\Delta_T)\in\Z$ denote the 
  intersection number $C_u\cdot\bar T$.  This is by definition the degree of
  the pull-back to $\bar T$ of the line bundle on $\fU$ associated to $C_u$.

\item If $\Delta_T$ has zero-dimensional finite part $\{u\}$, then we set 
   $\alpha(u,\Delta_T)\coloneq\deg_b(\Delta_T)$.
\end{enumerate}

Note that $\alpha(u,\Delta_T)$ can be calculated in any neighbourhood $\fU$ of
$\bar T$ as above, hence is intrinsic to $u$ and $T$.  We will also need to
remember intersection numbers with horizontal divisors, which we think of as
``data at infinity''.  Recall that the \emph{recession cone} of a polyhedron
$\Delta\subset\R^n$ is the polyhedral cone $\rho(\Delta)\subset\R^n$ consisting
of vectors $v\in\R^n$ such that $\Delta+rv\subset\Delta$ for all $r \geq 0$.
The rays (unbounded one-dimensional faces) of the recession cone
$\rho(\Delta_T)$ are in bijective correspondence with the
horizontal components of $D$ containing $T$ in their special fibre.  For such a
ray $r$ we let $H_r$ be the corresponding horizontal component, and we define
\[ \alpha(r,\Delta_T) \coloneq -H_r\cdot\bar T. \]
There is no problem computing this intersection product, as $H_r$ is a Cartier
divisor on all of $\sX$.

\begin{lem} 
  Let $\Delta_T$ be a $(d-1)$-dimensional canonical polyhedron of
  $S(\sX,H)$.  Then
  \begin{equation} \label{eq:deg.alphas}
    \deg_b(\Delta_T) = \sum_{u\in\Delta_T} \alpha(u,\Delta_T).
  \end{equation}
\end{lem}

\pf This is definitional if $\Delta_T$ has zero-dimensional finite part.
Otherwise, since $\bar T$ is a complete curve we have
\[ 0 = \div(\pi)\cdot\bar T = \sum_u C_u\cdot\bar T 
= \sum_{u\notin\Delta_T} C_u\cdot\bar T - \sum_{u\in\Delta_T} \alpha(u,\Delta_T) \]
where the sums are taken over vertices of $S(\sX,H)$ contained in
$\fU_\eta$, and all intersection products are taken in $\fU$.
Since all intersections are transverse, if $u\notin\Delta_T$ then
$C_u\cdot\bar T$ is the number of points in $V_u\cap\bar T$, i.e.\
the number of $d$-dimensional canonical polyhedra
containing both $u$ and $\Delta_T$.\qed

\paragraph[The divisor of a piecewise-affine function] 
We will define the divisor of a piecewise integral $\Gamma$-affine function $F$ on
$S(\sX,H)$ as a formal sum of $(d-1)$-dimensional canonical polyhedra.  In
analogy with the slope formula for
curves~\cite[Theorem~5.15]{bpr:analytic_curves}, we first define outgoing
slopes along a $d$-dimensional canonical polyhedron.

\begin{defn} \label{slope-definition}
  Let $\Delta_T$ be a $(d-1)$-dimensional
  canonical polyhedron of $S(\sX)$ and let $\Delta_S$ be a $d$-dimensional
  canonical polyhedron of $S(\sX,H)$ containing $\Delta_T$.  
  Let $F:\Delta_S\to\R$ be an integral $\Gamma$-affine function.
  \begin{enumerate}
  \item If $\Delta_S$ extends $\Delta_T$ in a bounded direction then there is a
    unique vertex $w$ of $\Delta_S$ not contained in $\Delta_T$, and we
    define the \emph{slope of $F$ at $\Delta_T$ along $\Delta_S$} to be the
    quantity 
    \begin{equation} \label{eq:slope.bounded}
      \slope(F;\Delta_T,\Delta_S) 
      \coloneq \frac 1{v(\pi)} \bigg(F(w) - 
      \frac 1{\deg_b(\Delta_T)}\sum_{u\in\Delta_T} \alpha(u,\Delta_T)\,F(u)\bigg), 
    \end{equation}
    where the sum is over all vertices $u$ of $\Delta_T$ and $v(\pi)$ is the
    length of $\Delta_S$.  

  \item If $\Delta_S$ extends $\Delta_T$ in an unbounded direction then
    there is a unique ray $s$ of $\rho(\Delta_S)$ not contained in
    $\rho(\Delta_T)$, and we define the
    \emph{slope of $F$ at $\Delta_T$ along $\Delta_S$} to be 
    \begin{equation} \label{eq:slope.unbounded}
      \slope(F;\Delta_T,\Delta_S) \coloneq d_{s} F
      - \frac 1{\deg_u(\Delta_T)}\sum_{r\subset\rho(\Delta_T)}
      \alpha(r,\Delta_T)\, d_r F.
    \end{equation}
    Here the second sum is over all rays of the recession cone $\rho(\Delta_T)$, and 
    for a ray $r$ we denote by $d_rF$ the derivative of $F$ along
    the primitive vector in the direction of $r$.
  \end{enumerate}

\end{defn}

\subparagraph \label{par:explain.slopes}
Definitions~\eqref{eq:slope.bounded}
and~\eqref{eq:slope.unbounded} require some 
explanation.  First we treat~\eqref{eq:slope.bounded}.
If $X$ is a curve, or more generally if $\Delta_T$ has
zero-dimensional finite part $\{u\}$, then 
$\alpha(u,\Delta_T) = \deg_b(\Delta_T)$ by definition, so
in this case $\slope(F;\Delta_T,\Delta_S) = (F(w)-F(u))/v(\pi)$.
This is the difference between the values of $F$ at the
endpoints of the edge in $\Delta_S$, divided by the length of $\Delta_S$.
If $\Delta_T$ has positive-dimensional finite part, suppose for simplicity
that $\Delta_T$ is bounded.
The problem with defining the slope in this situation 
is that the na\"ive slope $(F(w)-F(u))/v(\pi)$ may depend on the vertex
$u\in\Delta_T$. If all of the quantities
$\alpha(u,\Delta_T)$ were nonnegative then 
\begin{equation} \label{eq:m.T}
  m_T\coloneq\frac 1{\deg_b(\Delta_T)}
  \sum_{u\in\Delta_T} \alpha(u,\Delta_T)\,u 
\end{equation}
would be a point of $\Delta_T$ by~\eqref{eq:deg.alphas}, so since $F$ is
affine-linear on $\Delta_T$, we have
$\slope(F;\Delta_T,\Delta_S) = (F(v)-F(m_T))/v(\pi)$.  Interpreting $m_T$
as a weighted midpoint of $\Delta_T$, and declaring that all of $\Delta_T$
has distance $v(\pi)$ from $v$, we are again able to interpret
$\slope(F;\Delta_T,\Delta_S)$ as a slope of $F$ along a line segment.  In
general $\alpha(u,\Delta_T)$ need not be nonnegative, so $m_T$ cannot be
interpreted as a point of $\Delta_T$, and hence this explanation is more
of a heuristic. 

Now consider~\eqref{eq:slope.unbounded}.  Again if all of the 
$\alpha(r,\Delta_T)$ were nonnegative then the vector
\[ u_T \coloneq \frac 1{\deg_u(\Delta_T)} \sum_{r\subset \rho(\Delta_T)} \alpha(r,\Delta_T)\, r_0 \]
would be contained in  $\rho(\Delta_T)$, where $r_0$ denotes the primitive
vector on the ray $r$.  In this case $\slope(F;\Delta_T,\Delta_S)$ would be the
derivative of $F$ in the direction $s_0 - u_T$.  The primary argument for
the reasonableness of this definition is that it is the obvious
``linearized'' analogue of~\eqref{eq:slope.bounded}.

Note that in either case, $\slope(F;\Delta_T,\Delta_S)$ does not change if
we replace $F$ by $F+c$ for $c\in\R$.

\begin{remsub}
  When $X$ is a curve, $\slope(F;\Delta_T,\Delta_S)$ is always an integer.
  In the bounded case this follows from the fact that $F$ is integral
  $\Gamma$-affine on each canonical polyhedron, and in the unbounded case
  it follows from the fact that $d_s F\in\Z$.
  In higher dimensions the slopes need not be integers. See
  Remark~\ref{rem:in.terms.of.slopes}.  
\end{remsub}

\begin{defn} \label{defn:div.PL.func}
  Let $F:S(\sX,H)\to\R$ be a continuous function which is integral
  $\Gamma$-affine on   each canonical polyhedron.  For every
  $(d-1)$-dimensional polyhedron $\Delta_T$ of $S(\sX,H)$ we define
  \begin{align} 
    m(\div(F),\Delta_T) &\coloneq 
    \sum_{\substack{\Delta_S \succ \Delta_T\\\text{bounded}}} \slope(F;\Delta_T,\Delta_S) 
    \label{eq:ord.divF.DeltaT} \\
    m(\hat\div(F),\Delta_T) &\coloneq 
    \sum_{\Delta_S \succ \Delta_T} \slope(F;\Delta_T,\Delta_S), 
    \label{eq:ord.hatdivF.DeltaT}
  \end{align}
  where the first (resp.\ second) sum ranges over all bounded (resp.\
  bounded and unbounded) $d$-dimensional canonical polyhedra $\Delta_S$ of
  $S(\sX,H)$ containing $\Delta_T$.  Note that $m(\div(F),\Delta_T) = 0$ if
  $\Delta_T$ is unbounded.
  We define
  \begin{align*}
    \div(F) &\coloneq \sum_{\Delta_T} m(\div(F),\Delta_T)\,\Delta_T \\
    \hat\div(F) &\coloneq \sum_{\Delta_T} m(\hat\div(F),\Delta_T)\,\Delta_T,
  \end{align*}
  where the first (resp.\ second) sum ranges over the bounded (resp.\
  bounded and unbounded) $(d-1)$-dimensional canonical
  polyhedra $\Delta_T$ of $S(\sX,H)$.
\end{defn}

Let $\Delta_T$ be a  $(d-1)$-dimensional canonical polyhedron of
$S(\sX,H)$.  Substituting~\eqref{eq:slope.bounded}
into~\eqref{eq:ord.divF.DeltaT} gives
\begin{equation} \label{eq:dustins.divF}
  v(\pi) \,m(\div(F),\Delta_T) = 
  \sum_{\substack{\Delta_S \succ \Delta_T\\\text{bounded}}} F(w_S)
  - \sum_{u\in\Delta_T} \alpha(u,\Delta_T)\, F(u), 
\end{equation}
where $w_S$ is the unique vertex of $\Delta_S$ not contained in $\Delta_T$.
In Cartwright's situation, we have $v(\pi)=1$ and  $\Delta_T$ is a a canonical simplex of $S(\Xcal)$. Then $m(\div(F),\Delta_T)$ agrees with~\cite[Definition~1.5]{cartwright:tropical_complexes}.

\smallskip
We can now state the {\bf slope formula} for the
skeleton $S(\sX,H)$.

\begin{thm}[Slope formula] \label{PL for pair} 
  Let $f\in K(X)^\times$ be a rational function such that
  $\supp(\div(f))\subset\supp(H)_\eta$ and let
  $F = -\log|f|\big|_{S(\sX,H)}$.  Then $F$ is continuous and 
  integral $\Gamma$-affine on each canonical polyhedron of $S(\sX,H)$, and we have 
  \[ \hat\div(F) = 0. \]
\end{thm}

\smallskip
The identity $\hat\div(F)=0$ is a kind of balancing condition for $F$ on
$S(\sX,H)$, which is strongly analogous to the balancing condition for tropical
varieties. We require one lemma before proving Theorem~\ref{PL for pair}.

\begin{lem} \label{lem:sum.of.cartier}
  We use the notation in~\parref{par:slopes2}.  
  Suppose that $\Delta_T$ has positive-dimensional finite part and length
  $v(\pi)$.  
  Then there
  exists $\lambda\in K^\times$ such that the Cartier divisor
  $(\lambda\inv f)|_\fU$ is an integer linear combination of 
  $\{C_u\mid u\in\fU_\eta\text{ a vertex of } S(\sX,H)\}$ and 
  $\{H_i|_\fU\mid i=1,\ldots,S\}$.  More precisely,
  \begin{equation} \label{eq:lambda.f}
    \div(\lambda\inv f)|_\fU = \sum_u n_uC_u + 
    \sum_{i=1}^S \ord(f,\sH_i)\, H_i|_\fU 
  \end{equation}
  where $n_u = \frac 1{v(\pi)}(\ord(f,V_u)-v(\lambda))\in\Z$ and $u$ ranges over all vertices of $S(\sX,H)$ contained in $\fU_\eta$.
\end{lem}

\pf Let $\fU_x\subset\fU$ be a building block as in~\parref{par:slopes2}.
By the proof of Proposition~\ref{piecewise linear} and Remark \ref{order and representation},   we see that 
there exists $\lambda\in K^\times$ such that
$\frac 1{v(\pi)}(\ord(f,V_u)-v(\lambda))\in\Z$ for all
vertices $u$ of $S(\sX,H)$ contained in $(\fU_x)_\eta$.
It follows that for any two such vertices $u,u'$ we have 
$\ord(f,V_u)-\ord(f,V_{u'})\in v(\pi)\Z$.  This last statement is
independent of the choice of building block containing $u,u'$, so since
any two vertices of $S(\sX,H)$ contained in $\fU_\eta$ are connected by finite
faces of positive dimension, we see that
$\ord(f,V_u)-\ord(f,V_{u'})\in v(\pi)\Z$ for any two vertices $u,u'$ of
$S(\sX,H)$ contained in $\fU_\eta$.  
Choosing $\lambda$ as above (with respect to any choice of building
block), we have $\ord(\lambda\inv f,V_u)\in v(\pi)\Z$ for all vertices $u$
of $S(\sX,H)$ contained in $\fU_\eta$.  Letting 
$n_u = \frac 1{v(\pi)}\ord(\lambda\inv f,V_u)$, we
have the equality~\eqref{eq:lambda.f}, as both sides have the same Weil
divisor (use Proposition \ref{vertical Cartier divisors}).
\qed 

\medskip 

\proof[of Theorem~\ref{PL for pair}]
We showed that $F$ is continuous and integral $\Gamma$-affine on
canonical polyhedra  in Proposition~\ref{piecewise linear}. 
We have to prove that $m(\hat\div(F),\Delta_T) = 0$ for all
$(d-1)$-dimensional canonical polyhedra $\Delta_T$ of $S(\sX,H)$.
First suppose that $\Delta_T$ has positive-dimensional finite part and
length $v(\pi)$.  We use the notation in~\parref{par:slopes2}.  
Choose $\lambda\in K^\times$ as in Lemma~\ref{lem:sum.of.cartier}.  
Multiplying $f$ by a non-zero scalar does not change
$\hat\div(F)$, so we may replace $f$ with $\lambda\inv f$ to assume that we
have an equality of Cartier divisors
\[ \div(f)|_\fU = \sum_u n_u\, C_u + \sum_{i=1}^S \ord(f,\sH_i)\,H_i|_{\fU} \]
where $n_u = \frac 1{v(\pi)}\ord(f,V_u)\in\Z$ and $u$ ranges over all vertices of $S(\sX,H)$ contained in $\fU_\eta$.
As $\div(f)$ is a principal Cartier divisor, we have
\begin{equation} \label{eq:f.dot.T} 
  0 = \div(f)|_\fU\cdot\bar T = \sum_u n_u\, C_u\cdot\bar T
  + \sum_{i=1}^S \ord(f,\sH_i)\,H_i\cdot\bar T. 
\end{equation}
For $u$ a vertex of $S(\sX,H)$ contained in $\fU_\eta$, by definition we have 
$F(u) = \ord(f,V_u) = v(\pi)\, n_u$.
Substituting into~\eqref{eq:f.dot.T}, we have
\begin{equation} \label{eq:bounded.plus.unbounded.parts}
  0 = \frac 1{v(\pi)} \sum_u F(u)\, C_u\cdot\bar T 
  + \sum_{i=1}^S \ord(f,\sH_i)\,H_i\cdot\bar T. 
\end{equation}
If $u\notin\Delta_T$ then $C_u\cdot\bar T$ is equal to the number of
canonical polyhedra $\Delta_S$ containing $\Delta_T$ and $u$.  If
$u\in\Delta_T$ then $C_u\cdot\bar T = -\alpha(u,\Delta_T)$, so 
\begin{equation} \label{eq:bounded.part}
\begin{split} \frac 1{v(\pi)} \sum_u F(u)\,C_u\cdot\bar T
 &= \frac 1{v(\pi)}\bigg(\sum_{\substack{\Delta_S \succ \Delta_T\\\text{bounded}}} F(w_S)
- \sum_{u\in\Delta_T} \alpha(u,\Delta_T)\, F(u)\bigg) \\
&= \sum_{\substack{\Delta_S \succ \Delta_T\\\text{bounded}}}\frac 1{v(\pi)} \bigg(
F(w_S) - \frac 1{\deg_b(\Delta_T)} \sum_{u\in\Delta_T} \alpha(u,\Delta_T)\,F(u)\bigg)\\
&= \sum_{\substack{\Delta_S\succ \Delta_T\\\text{bounded}}} \slope(F;\Delta_T,\Delta_S),
\end{split}
\end{equation}
where the first sum runs over all canonical polyhedra $\Delta_S$ extending
$\Delta_T$ in a bounded direction, and $w_S$ is the vertex of $\Delta_S$ not
contained in $\Delta_T$.  

Recall that a ray $r$ contained in the recession cone $\rho(\Delta_R)$ of a canonical polyhedron $\Delta_R$ of $S(\sX,H)$ 
corresponds to a horizontal component $\sH_r$ containing the stratum $R$ in the special 
fibre. Moreover,  the multiplicity
$\ord(f,\sH_r)$ is equal to the derivative of $F$ along the primitive
vector on the ray of $\rho(\Delta_T)$ corresponding to $H_r$:
see~\parref{tropical divisor at infinity} and~\parref{par:slopes2}.
If $T$ is not contained in $\sH_i$ for some $i \in \{1, \dots, R+S\}$, then either 
$T \cap \sH_i = \emptyset$ or there is a $d$-dimensional canonical polyhedron $\Delta_S$ 
which extends $\Delta_T$ in the unbounded direction $r_i$ for a ray $r_i$ with $\sH_{r_i}=\sH_i$. 
The first case is equivalent to $H_i \cdot \bar T = 0$. In the second case, $H_i\cdot\bar T$ is equal to the number of $d$-dimensional canonical polyhedra $\Delta_S$
extending $\Delta_T$ in an unbounded direction corresponding to $\sH_i$.   
If $T$ is contained in $\sH_i$, then there is a ray $r_i \subset \rho(\Delta_T)$ with $\sH_{r_i}=\sH_i$ and $H_i\cdot\bar T = -\alpha(r_i,\Delta_T)$.  It follows that
\begin{equation} \label{eq:unbounded.part}
\begin{split}
  \sum_{i=1}^S \ord(f,\sH_i)\,H_i\cdot\bar T
  &= \sum_{\substack{\Delta_S \succ \Delta_T\\\text{unbounded}}} d_s F
  - \sum_{r\subset \rho(\Delta_T)} \alpha(r,\Delta_T)\, d_r F \\
  &= \sum_{\substack{\Delta_S \succ \Delta_T\\\text{unbounded}}} \bigg(
  d_s F - \frac 1{\deg_u(\Delta_T)}\sum_{r\subset \rho(\Delta_T)}
  \alpha(r,\Delta_T)\, d_r F \bigg) \\
  &= \sum_{\substack{\Delta_S \succ \Delta_T\\\text{unbounded}}} 
  \slope(F;\Delta_T,\Delta_S),
\end{split}
\end{equation}
where the first sum runs over the canonical polyhedra $\Delta_S$ extending
$\Delta_T$ in an unbounded direction, the second runs over all rays of
$\rho(\Delta_T)$, and $s$ is the ray of $\rho(\Delta_S)$ not contained in
$\rho(\Delta_T)$.  Combining~\eqref{eq:bounded.plus.unbounded.parts}
with~\eqref{eq:bounded.part} and~\eqref{eq:unbounded.part},
we obtain
\begin{equation} \label{concluding argument} 
 0 = \sum_{\Delta_S \succ \Delta_T} \slope(F;\Delta_T,\Delta_S) =
m(\hat\div(F),\Delta_T). 
\end{equation}

Now suppose that $\Delta_T$ has zero-dimensional finite part $\{u\}$.  A
separate argument is needed as the vertical components of $\sX_s$ are not
necessarily Cartier on a formal open subscheme containing $\bar T$.  In this case 
$\bar T = V_u\cap(\sH_{i_2})_s\cap\cdots\cap(\sH_{i_{d}})_s$ for some
$i_2,\ldots,i_{d}\in\{1,\ldots,S\}$ and  for a vertex $u$ of $S(\sX,H)$. 
Replacing $f$ by $\lambda^{-1}f$ for a suitable non-zero $\lambda \in K$, we may assume
that $\ord(f,V_u)=0$. 
Let 
us consider the  Cartier divisor $E \coloneq \sum_{j=1}^{S}\ord(f,\sH_{j}) H_{j}$ on $\sX$. 
Then the Cartier divisor  $D \coloneq \Div(f)- E$ has support in the special fibre of $\sX$, and we have
$$\ord(D,V_u)= \ord(f,V_u)=F(u)=0.$$
We conclude that $D$ intersects $\bar T$ properly which means that the intersection of $\bar T$ with the support of $D$ is zero-dimensional.
Then the intersection product $D.\bar T$ is a well-defined cycle on $\bar T$ supported in the union of all zero-dimensional strata $S$ such 
that $\Delta_S$ extends $\Delta_T$ in a bounded direction. 

The multiplicity $m_S$ of $D.\bar T$ in $S$ may be computed on a building block $\fU$ with distinguished stratum $S$. The  finite part of 
$\Delta_S$ is an edge from $u$ to another vertex $w$ of length $v(\pi)$. 
By Proposition~\ref{prop:comps.are.cartier}, there is a unique effective Cartier divisor $C$ on $\fU$ with $\cyc(C)=v(\pi) \cdot (V_w \cap \fU_s).$
We have 
$$\cyc(D|_\fU)= \ord(f,V_w) \cdot (V_w \cap \fU_s)$$
and hence we deduce from
Proposition \ref{vertical Cartier divisors} that
$$D|_\fU = \frac{\ord(f,V_w)}{v(\pi)} \cdot C.$$
Using that $\fU$ is a strictly semistable formal scheme, the multiplicity $m_S$ of $D.\bar T$ in $S$ is 
$$m_S = \frac{\ord(f,V_w)}{v(\pi)} = \frac{1}{v(\pi)}(F(w)-F(u))=  \slope(F;\Delta_T,\Delta_S).$$
We conclude that 
\begin{equation} \label{intersection formula 1} 
 D \cdot \bar T = \sum_{\substack{\Delta_S\succ \Delta_T\\\text{bounded}}} \slope(F;\Delta_T,\Delta_S).
\end{equation}
On the other hand, we have 
$$0= \div(f) \cdot \bar T = D \cdot \bar T + E \cdot \bar T.$$
We substitute \eqref{intersection formula 1} for $D \cdot \bar T$. 
Moreover, $E \cdot \bar T$ can be calculated as in \eqref{eq:unbounded.part} as we have not used positive dimensionality 
of the finite part of $\Delta_T$ in the argument.  Now the claim follows as in \eqref{concluding argument}. \qed

\paragraph One could imagine defining $\hat\div(F)$ as a formal
sum which includes the $(d-1)$-dimensional canonical polyhedra $\Delta_S'$ in the boundary of
the compactified skeleton $\hat S(\sX,H)$, with the multiplicity
$m(\hat\div(F),\Delta_S')$ being determined by an outgoing slope of $F$.
In this case the correct statement would be $\hat\div(F) + \hat\tau( f) = 0$, where $\hat\tau(f)$ is 
the retraction of the principal Weil divisor $\cyc(f)$ to $\hat S(\sX,H)$ (see \ref{tropical divisor at infinity}).
We leave this reformulation to the interested reader. 

As a consequence of Theorem~\ref{PL for pair},
we derive the {\bf slope formula} for the skeleton $S(\sX)$.

\begin{thm} \label{PL for bounded} 
  Let $f\in K(X)^\times$ be a rational function such that
  $\supp(\div(f))\subset\supp(H)_\eta$ and let
  $F = -\log|f|\big|_{S(\sX)}$.  Then $F$ is continuous and 
  integral $\Gamma$-affine on each canonical polyhedron of $S(\sX)$, and we have 
  \[ \div(F) + \tau(f) = 0, \]
	where $\tau(f)$ is the retraction of the principal Weil divisor $\cyc(f)$ to $S(\Xcal)$ from \ref{tropical divisor} .
\end{thm}

\proof Let $\Delta_T$ be a bounded $(d-1)$-dimensional canonical polyhedron. 
It remains to prove that $m(\div(F),\Delta_T) + m(\tau(f),\Delta_T)= 0$.
We have
\[\begin{split}
  0 = m(\hat\div(F),\Delta_T) 
  &= \sum_{\substack{\Delta_S \succ \Delta_T\\\text{bounded}}} \slope(F;\Delta_T,\Delta_S)
  + \sum_{\substack{\Delta_S \succ \Delta_T\\\text{unbounded}}} \slope(F;\Delta_T,\Delta_S)\\
  &= m(\div(F),\Delta_T) 
  + \sum_{\substack{\Delta_S \succ \Delta_T\\\text{unbounded}}} \slope(F;\Delta_T,\Delta_S),
\end{split}\]
so we need to argue that the last sum is equal to 
$m(\tau(f),\Delta_T)$.  If $\Delta_S$ extends $\Delta_T$ in an unbounded
direction then since $\rho(\Delta_T)$ contains no rays, we have
$\slope(F;\Delta_T,\Delta_S) = d_sF = \ord(f,\sH_i)$, where $s$ is the ray of
$\rho(\Delta_S)$ and $\sH_i$ is the horizontal component containing $S$.
Using~\parref{tropical divisor} and \eqref{eq:unbounded.part}, we have
$$m(\tau(f),\Delta_T)= \sum_k \ord(f, \sH_k) \cdot \#(\bar T\cap(\sH_k)_s)=  \sum_k \ord(f, \sH_k) H_k \cdot \bar T= \sum_{\substack{\Delta_S \succ \Delta_T\\\text{unbounded}}} \slope(F;\Delta_T,\Delta_S),$$ 
which finishes the proof.\qed

\begin{eg}
  Suppose that $\dim(X) = 1$, i.e.\ that $X$ is a curve.
  In this case $\hat S(\sX,H) = S(\sX,H)\cup\supp(H)$, and
  identifying canonical polyhedra on the boundary of
  $\hat S(\sX,H)$ with points of $X(K)$ identifies
  $\hat\tau(f)$ with $\div(f)$.  
  If $\div(f) = \sum n_x \cdot x$ then
  $\tau(f) = \sum n_x \cdot \tau(x)$, which is the retraction 
  $\tau_*(\div(f))$ of the divisor
  $\div(f)$ to the skeleton $S(\sX,H)$.  Hence the identity
  $\div(F) + \tau(f) = 0$ of Theorem~\ref{PL for bounded} reads
  \[ \div(-\log|f|) = -\tau_*(\div(f)). \]

  As explained in~\parref{par:explain.slopes}, if 
  $\Delta_T = \{x\}\subset S(\sX,H)$ is a zero-dimensional canonical polyhedron
  then for every segment or ray $\Delta_S$ containing $x$, 
  the slope of $F$ at $\Delta_T$ along $\Delta_S$ is simply the 
  outgoing slope of $F$ in the direction of $\Delta_T$.  
  The statement $\hat\div(F)=0$ of Theorem~\ref{PL for pair} therefore
  says that the sum of the outgoing slopes of $F$ at $x$ is zero.  
  This, along with the other results we have proved, essentially recovers
  the slope formula for curves~\cite[Theorem~5.15]{bpr:analytic_curves}.
\end{eg}

\begin{rem}
  It is not hard to show that in dimension $1$, the divisor
  $\tau(f) = \tau_*\div(f)$ and the harmonicity condition 
  $\div(F) + \tau(f) = 0$ of Theorem~\ref{PL for pair} 
  uniquely determine the piecewise linear function $F$ on $S(\sX)$ up to additive translation (see \cite[\S 1.2.3]{thuillier:thesis} or  \cite[Proposition~3.2(A)]{baker_rumely:book}).
  In higher
  dimensions, this Neumann problem in polyhedral geometry need not have a
  unique solution:  there exist piecewise linear spaces $S$ with
  $\alpha$-numbers $\alpha(\scdot,\scdot)$ and
  non-constant piecewise linear functions 
  $F:S\to\R$ as in Definition~\ref{defn:div.PL.func} such that
  $\div(F) = 0$.  Such an example can be found in~\cite[Example~6.11]{cartwright:tropical_complexes}.  
	Dustin Cartwright has pointed out to the authors that his example arises  as a skeleton $S(\Xcal)$ 
	of a strictly semistable pair $(\Xcal,H)$, where $\Xcal$ is a toric scheme over a (discrete) valuation ring. 
	Moreover, his non-constant piecewise linear functions 
  $F$ with $\div(F) = 0$ is the restriction of $-\log|f|$ to $S(\Xcal)$ for a non-constant rational function $f$ on $\Xcal$. 
  
	Note that this cannot happen on $S(\Xcal,H)$ for any strictly semistable pair $(\Xcal,H)$. If $f$ is a non-zero rational function on $\Xcal$ and 
	if $F$ is the restriction of $-\log|f|$ to $S(\Xcal,H)$, then 
	$\hat\div(F)=0$ yields that $f$ is constant as 
	$f$ is determined up to multiples by  $\cyc(f)$  and hence also by $\hat\tau(f)$. 
\end{rem}

\section{A two-dimensional example} \label{A two-dimensional example}

In this section, we give an example to illustrate the slope
formula for skeletons. We choose the square of a Tate elliptic curve which
fits well for analytic purposes. As before, $K$ is an algebraically closed
field complete with respect to the valuation $v$ and with non-zero value group
$\Gamma \subset \R$. We assume that the residue field $\ktilde$ of $K$ has characteristic $\neq 2$.  

\paragraph \label{totally degenerate}
Recall that an abelian variety $A$ over $K$ is \emph{totally degenerate} if
$A^{\rm an}=\Tan /P$ where $T=\Spec(K[M])$ is a multiplicative torus and
$P$ is a lattice in $\Tan$. Here, a lattice means a discrete subgroup $P$ of $\Tan$ contained
in $T(K)$ such that $\trop:\Tan \rightarrow N_\R$ maps $P$ isomorphically onto a
complete lattice $\Lambda$ of $N_\R$, where $N$ is the dual of the character
lattice $M$ of $T$. Passing to the quotient, we get a continuous proper map
$\tropbar: A^{\rm an} \rightarrow N_\R / \Lambda$.

\begin{eg} \label{Tate curve}
The totally degenerate abelian varieties of dimension $1$ are called Tate elliptic curves. For every $q \in K^\times$ with $v(q) > 0$, we have a Tate elliptic curve $E$ given analytically by $T_1^{\rm an}/q^\Z$ with torus $T_1=\Spec(K[\Z])$ and lattice $q^\Z$. Algebraically, the elliptic curve $E$ is given by the generalized Weierstrass equation 
$$y^2 + xy = x^3 + a_4x + a_6,$$
where $a_4,a_6$ are given by the following convergent power series
$$a_4= - 5\sum_{n=1}^\infty n^3q^n/(1-q^n), \quad a_6=-\frac{1}{12}\sum_{n=1}^\infty(7n^5+5 n^3)q^n/(1-q^n).$$
The isomorphism $T_1^{\rm an}/q^\Z \rightarrow E^{\rm an}$ is given by 
\begin{equation} \label{Laurent expansions}
 \begin{split}
  x(\zeta)&=\sum_{n=-\infty}^\infty q^n\zeta/(1-q^n\zeta)^2 - 2 \sum_{n=1}^\infty nq^n/(1-q^n)\\
  y(\zeta)&= \sum_{n=-\infty}^\infty q^{2n}\zeta^2(1-q^n\zeta)^3 + \sum_{n=1}^\infty nq^n/(1-q^n),
 \end{split}
\end{equation}
where $\zeta$ is the torus coordinate on $T_1$. In the following, we will identify $E^{\rm an}$ with $T_1^{\rm an}/q^\Z$ using this isomorphism of analytic groups. 
This is due to Tate~\cite{tate:elliptic_functions} (see also~\cite[Theorem C14.1]{silverman:I}).
\end{eg}

\paragraph \label{Mumford model}
Let $A$ be a totally degenerate abelian variety over $K$ as in~\parref{totally degenerate}. We define a {\it polytope} in $N_\R / \Lambda$ as a subset $\Deltabar$ of $N_\R / \Lambda$ given as the bijective image of a polytope $\Delta$ in $N_\R$ with respect to the quotient homomorphism. Recall that a {\it polytopal decomposition} of a set $S$ in $N_\R$ is a polytopal complex with support $S$. A {\it polytopal decomposition} $\Ccalbar$ of $N_\R/ \Lambda$ is a finite collection of polytopes in $N_\R / \Lambda$ induced by an infinite $\Lambda$-periodic polytopal decomposition $\Ccal$ of $N_\R$. We will assume always that $\Ccalbar$ is integral $\Gamma$-affine which means that all polytopes of $\Ccal$ are integral $\Gamma$-affine polytopes in $N_\R$. 

For $\Delta \in \Ccal$, we have a polytopal subdomain $U_\Delta = \trop^{-1}(\Delta)$ of $\Tan$ which is the Berkovich spectrum of the strictly affinoid algebra 
$$K \langle U_\Delta \rangle = \bigg\{ \sum_{u \in M} \alpha_u \chi^u  \mid \lim_{|u|\to \infty } v(\alpha_u) + \langle u, \omega \rangle = \infty \; \forall \omega \in \Delta\bigg\},$$
where $\chi^u$ is the character of $T$ corresponding to $u \in M$, the coefficients $\alpha_u$ of the power series are in $K$ and where $|u|$ uses any norm on $M_\R$. The supremum norm is given here by
$$\bigg\|\sum_{u \in M}\alpha_u \chi^u \bigg\|_{\rm sup} =  \max_{\omega \in \Delta, \, u \in M} |\alpha_u|\, e^{-\langle u, \omega \rangle }.$$
Recall that $K \langle U_\Delta \rangle^\circ$ is the subalgebra of $K \langle U_\Delta \rangle$ given by the Laurent series of supremum norm $\leq 1$. 
The canonical formal $\kcirc$-model of $U_\Delta$ is $\fU_\Delta \coloneq \Spf(K
\langle U_\Delta \rangle^\circ)$.  This is an affine admissible formal
$K^\circ$-scheme, as mentioned in~\ref{par:notation}.

The admissible formal schemes $(\fU_\Delta)_{\Delta \in \Ccal}$ glue together to an admissible formal scheme $\fE$ which is a $\kcirc$-model of $\Tan$. Passing to the quotient $\fA\coloneq\fE / P$, we get a $\kcirc$-model of $A$. It is called the {\it Mumford model associated to the polytopal decomposition $\Ccalbar$}. 
For more details about this construction, we refer the reader to \cite[\S4, \S6]{gubler:tropical}.

\paragraph \label{strictly semistable Mumford model}
A simplex $\Delta$ in $N_\R$ is called {\it regular} if there is a basis
of $N$ and $a>0, a\in\Gamma$ such that $\Delta$ is a translate of
$\{\omega \in N_\R \mid \omega_i \geq 0, \,\omega_1 + \dots + \omega_n
\leq a\}$, where $\omega_1, \dots, \omega_n$ are the coordinates of
$\omega$.  
If $\Delta$ is regular then $\fU_\Delta\cong\fU_{\Delta(n+1,\pi)}$, where
$v(\pi) = a$.
Suppose that $\Ccalbar$ is induced by an infinite $\Lambda$-periodic polytopal decomposition $\Ccal$ of $N_\R$ as above. 
If $\Ccal$ consists of regular simplices, then it is clear from the definitions that the Mumford model $\fA$ associated to $\Ccalbar$ is a strictly semistable formal scheme. Then the skeleton $S(\fA)$ may be identified with
 $N_\R / \Lambda$, the triangulation $\Ccalbar$ induces the decomposition of $S(\fA)$ into canonical simplices and 
$\tropbar$ is the canonical retraction $\tau: A^{\rm an} \rightarrow S(\fA)$. This follows from \cite[Proposition 6.3]{gubler:tropical} and the definitions.

\paragraph \label{generic fibre of running example} 
In the rest of this section, we will focus on the following special case of the above
construction. We consider the Tate elliptic curve $E$ associated to $q \in
K^\times$ which means $E^{\rm an}=T_1^{\rm an}/q^\Z$ for the
one-dimensional torus $T_1=\Spec(K[\zeta^{\pm 1}])$.
 For simplicity, we assume $v(q)=1$. Then  $A\coloneq E^2$ is a totally
 degenerate abelian variety over $K$ given analytically by $\Tan/P$, where $T=T_1 \times T_1$ and $P=q^\Z \times q^\Z$. For the above lattices, we have $M=\Lambda=\Z \times \Z$. Referring to~\parref{Tate curve}, we denote by $\zeta_1$ (resp. $\zeta_2$) the pull-back of the torus coordinate $\zeta$ and by $x_1,y_1$ (resp. $x_2,y_2$) the pull-back of the algebraic coordinates $x,y$ of the generalized Weierstrass equation with respect to the projection to the the first (resp.\ second) factor of $E^2$. 

\paragraph \label{ss pair of running example}
We choose the regular triangulation $\Ccal$ of $\R^2$ obtained by $\Lambda$-translations from the natural unit square by dividing it in four squares and drawing the two diagonals in the original unit square.
See Figure~\ref{fig:square}.

\genericfig{square}{The triangulation $\sC$ of~\parref{ss pair of running example}.}

 Let $\fA$ be the Mumford model of $A$ associated to $\Ccal$.  We have seen in~\parref{Mumford model} that $\fA$ is a strictly semistable $\kcirc$-model of $A=E^2$. By~\parref{strictly semistable Mumford model},   the skeleton and its canonical simplices are induced by the decomposition $\Ccalbar$ of $\R^2/\Lambda$. The unit square serves as a fundamental lattice and we number its vertices by $P_1=(0,0), \dots ,P_9=(1,1)$, where we identify for example the vertices $P_1$, $P_3=(1,0)$, $P_7=(0,1)$ and $P_9$ accoding to $\Lambda$-translation. An edge in $\Ccalbar$ with the two vertices $P_i$ and $P_j$ is denoted by $e_{ij}$ and a face in $\Ccalbar$ with the three vertices $P_i,P_j,P_k$ is denoted by $\Delta_{ijk}$. By the stratum--face correspondence in Proposition \ref{stratum-face correspondence}, $S(\fA)$ has 4 two-dimensional strata, 12 one-dimensional strata and 8 zero-dimensional strata.

\begin{prop} \label{running example is algebraic}
The Mumford model $\fA$ obtained from the regular triangulation above is algebraic.
\end{prop}

\proof Let $L$ be an ample line bundle on $A=E^2$. 
By \cite[6.5]{gubler:tropical}, $L$ induces a positive definite symmetric bilinear form $b$ on $\Lambda$ and a cocycle $\lambda \mapsto z_\lambda$ of $H^1(\Lambda, C(\R^2))$ with 
$$z_\lambda(\omega)=z_\lambda(0)+b(\omega,\lambda)$$
for $\lambda \in \Gamma$ and $\omega \in \R^2$. It follows from \cite[Proposition 6.6]{gubler:tropical} that every $f \in C(\R^2)$ with the following two conditions (a) and (b) induces a line bundle $\fL$ on $\fA$ with generic fibre $L$:
\begin{itemize}
 \item[(a)] The restriction of $f$ to $\Delta \in \Ccal$ is integral $\Gamma$-affine.
 \item[(b)] $f(\omega + \lambda)=f(\omega) + z_\lambda(\omega)$ for all $\omega \in \R^2$ and all $\lambda \in \Lambda$. 
\end{itemize}
We choose a root $\sqrt{q} \in K^\times$ leading to a $2$-torsion point
$P\coloneq[\sqrt{q}]$ of $E$. Then $L_0\coloneq \sO(P \times E + E \times P)$ is an ample line bundle on $A=E^2$. 
It is shown in \S 3.3 of Christensen's thesis~\cite{christensen:thesis}
that  there exists a strictly convex, piecewise linear, continuous real function
$f_0 \in C(\R^2)$ satisfying (b)  and a non-zero $m \in \N$ such that 
$f\coloneq mf_0$ satisfies (a). Moreover, the maximal domains of linearity for Christensen's $f_0$ are the two-dimensional simplices of $\Ccal$. We conclude that $L\coloneq L_0^{\otimes m}$ 
is an ample line bundle with a formal $\kcirc$-model $\fL$ on the Mumford model $\fA$. By \cite[Corollary 6.7]{gubler:tropical}, the restriction of $\fL$ to the special fibre $\fA_s$ is an ample line bundle. It follows from Grothendieck's algebraization criterion (see \cite[Theorem 5.4.5]{egaIII_1} in the case of discrete valuations and the generalization to arbitrary real valuations by Ullrich in~\cite[Proposition~6.9]{ullrich:direct_img}
that $\fA$ is algebraic. 
\qed

\smallskip
Let $\sA$ be the algebraization of $\fA$.  This is a strictly semistable
algebraic $\kcirc$-model of $A$.

\paragraph \label{rational function} 
In our running example $A=E^2$, we
choose the rational function $f=x_1-x_2$ using the algebraic coordinates
from~\parref{generic fibre of running example}. Then the divisor
of $f$ on $A$ is equal to the sum of the diagonal and the anti-diagonal in
$A=E \times E$ minus $E \times 0 + 0 \times E$. If we consider the horizontal divisor $H'$ given as  the sum of the closures of the diagonal, the anti-diagonal, $E \times 0$ and $0 \times E$, then 
$(\sA,H')$ is not a strictly semistable pair on $\Acal$ as the zero
element 
lies in four components (at most two would be allowed). We blow up $\Acal$ in
the closure $\Ycal$ of the zero element of $A$.  (Note that this is a closed
subscheme of $\sA$ but not of $\sA_s$.) This leads to a strictly
semistable pair $(\Xcal,H)$, where the horizontal components 
$H_1,\dots , H_5$ are given as follows: The strict transform of the closure of
the diagonal 
(resp. anti-diagonal) with respect to the blow up is denoted by $H_1$ (resp. $H_2$). The
strict transform of the closure of $E \times 0$ (resp. $0 \times E$) is denoted
by $H_3$ (resp. $H_4$). The exceptional divisor of the blow up is denoted
by $H_5$.  To be clear, $\sX$ is not a $K^\circ$-model of $A$, but rather of
a blowup of $A$.

The skeleton $S(\Xcal,H)$ is obtained from $S(\Acal)$ in the following
way. First, we note that $\Xcal_s$ still has four vertical components 
$V_1, V_2, V_4, V_5$ lying over the irreducible components of $\Acal_s$
corresponding to $P_1,P_2,P_4,P_5$. Hence $\sX\to\sA$ induces an identification
$S(\sX)\isom S(\sA)$.  We have described $S(\Acal)$ as the quotient
of  $\R^2$ by the group action of $\Lambda$. Now we add to the plane $\R^2$ five
new independent directions $b_1, \dots ,b_5$ corresponding to the
horizontal components $H_1, \dots, H_5$. Then we expand the edges $e_{15}$
and $e_{59}$ (resp.\ $e_{35}$ and $e_{57}$) in $S(\Acal)$ to half-stripes
in the $b_1$-direction (resp. $b_2$-direction). They correspond to the
strata of $D$ in the intersection of two vertical components with either
$H_1$ or $H_2$.  Similarly, we expand the edges $e_{12}$ and $e_{23}$
(resp.\ $e_{14}$ and $e_{47}$) to half-stripes in the $b_3$-direction
(resp. in the $b_4$-direction). They correspond to the strata of $D$ in the
intersection of two vertical components with either $H_3$ or $H_4$. 

Over $P_5$, we fill in two quadrants between $b_1$ and
$b_2$ which both  have the same two edges given by the halflines starting in $P_1$ in the directions $b_1$ and $b_2$. This corresponds to the two strata points in the intersection of $H_1, H_2$ and
$V_5$.  Note that 
we use here that the residue characteristic is not $2$.
Over
$P_1$, we add the $5$ quadrants filling in between $(b_1,b_5)$, $(b_2,b_5)$, $(b_3,b_5)$, $(b_4,b_5)$ and $(b_1,b_2)$. The first four quadrants correspond to the single point in the intersection of $H_i$, $H_5$ and $V_1$ for $i \in \{1, \dots, 4\}$. We note that $H_1$, $H_2$ and $V_1$ intersect only in one point and this corresponds to the last quadrant. There are no other intersections over $0$ as the blow up separates $H_i$ and $H_j$ for $i \neq j$ in $\{1, \dots ,4\}$.  

\paragraph \label{PL in running example} 
Our goal is to illustrate the
slope formula (Theorem~\ref{PL for bounded}) for
$F\coloneq-\log|f|$ on the skeleton $S(\Acal) = S(\Xcal)$. 
By \ref{tropical divisor} and \ref{rational function}, the retraction $\tau(f)$ to $S(\Xcal)$ is given by
\begin{equation} \label{eq:delta.f.eg}
  \tau(f) = e_{15} + e_{59} + e_{35} + e_{57} - e_{12} - e_{23} - e_{14} - e_{47}.
\end{equation}
This is the only part where we use the strictly semistable pair
$(\Xcal,H)$. The remaining computations can be done solely on the Mumford
model $\Acal$ and on the skeleton $S(\Acal)=S(\Xcal)$. In particular, the
projection formula in Proposition \ref{properties of refined intersection} shows that we may compute the occurring intersection
numbers on the model $\Acal$. The vertices $P_1,P_2,P_4,P_5$ correspond to
the irreducible components of $\Acal_s$ which we denote by
$Y_1,Y_2,Y_4,Y_5$. Let $D_i\coloneq D_{P_i}$ be the Cartier divisor
associated to the vertex $P_i$ by
Proposition~\ref{prop:comps.are.cartier}. 

We illustrate the slope formula by showing that
$-m(\div(F),e_{15}) = m(\tau(f), e_{15}) = 1$.  Here 
the edge $e_{15}$ corresponds to the one-dimensional stratum $T_{15}$. 
By~\eqref{eq:dustins.divF} we have
\begin{equation} \label{runningPL}
  \begin{split}
  \frac 12 m(\div(F),e_{15}) &= 
  F(P_2) + F(P_4) - \alpha(P_1, e_{15})\,F(P_1) - \alpha(P_5,e_{15})\,F(P_5)\\
  &= F(P_2) + F(P_4) + (D_1\cdot\bar{T_{15}})\,F(P_1) 
  + (D_5\cdot\bar{T_{15}})\,F(P_5).
  \end{split}
\end{equation} 
We must show that the right side of this equation is $-\frac 12$.

\paragraph \label{slopes in running example} 
We compute now the quantities on the right side of
\eqref{runningPL}. Consider a point $\xi$ in the skeleton $S(T_1^\an)$ of
the torus  $T_1=\Spec(K[\zeta^{\pm 1}])$.  If
$|q|^{1/2} < |\zeta(\xi)| < 1$ then the unique summand in the Laurent
expansion~\eqref{Laurent expansions} of $x(\zeta(\xi))$ with maximal absolute
value is $\zeta(\xi)$.  In particular, $|x(\zeta(\xi))| = |\zeta(\xi)|$.
The ultrametric inequality and continuity of $F$
then imply that $F(P_1) = F(P_2) = F(P_4) = 0$ and
$F(P_5) = 1/2$.  Therefore we only need to prove that
$D_5\cdot\bar{T_{15}} = -1$.

\paragraph \label{intersection numbers in the running example} 
To compute
the intersection number $D_5\cdot \overline{T_{15}}$ from
\eqref{runningPL}, we are going to use Kolb's relations (see
Proposition~\ref{Kolb's relations}). The problem is that the canonical
simplices of $S(\Acal)$ are not determined by their vertices: 
for instance, $\Delta_{125}$ and $\Delta_{578}$ have the same vertices. To
deal with 
that, we pass to a covering $\hat{\varphi}:\fA' \rightarrow \fA$, where
$\fA'$ is the Mumford model of $A=E^2$ induced by the regular
triangulation $\frac{1}{2}\Ccal$ of $\R^2$ and where $\hat{\varphi}$ on
the
generic fibres is multiplication by $2$. The fundamental lattice is still
the unit square, but it is now divided up into $16$ squares. We number
the vertices by $Q_1=(0,0), \dots ,Q_{25}=(1,1)$
as in~\parref{ss pair of running example}. 
See Figure~\ref{fig:square2}. 
Similarly as in Proposition \ref{running example is algebraic}, we can show that $\fA'$ is 
algebraic. Indeed, the same proof works using the convex function $f_0 \circ [2]$. We denote the strictly semistable algebraic $\kcirc$-model by $\Acal'$ and let 
$\varphi:\Acal' \rightarrow \Acal$ be the underlying algebraic morphism. 

\genericfig{square2}{The triangulation $\frac 12\sC$ of~\parref{intersection numbers in the running example}.}

Let $D_j'$
(resp. $Y_j'$) be the Cartier divisor (resp. the irreducible component of $\Acal_s'$)
associated to $Q_j$. Since we may use $\pi=q^{1/4}$ for the strictly
semistable model $\Acal'$, the Weil divisor associated to $D_j'$ is
\begin{equation} \label{assoweil}
\cyc(D_j') = v(q^{1/4})Y_j'=\frac{1}{4}Y_j'.
\end{equation}
Here and in the following, we use the refined intersection theory with Cartier divisors which we recall in Appendix~\ref{Refined intersection theory with Cartier divisors}. It follows easily from the definition of the Cartier divisors $D_i$ and $D_j'$ that 
\begin{equation} \label{pullbackd1}
\varphi^*(D_5)=2(D_7'+D_9'+D_{17}'+D_{19}').
\end{equation}
Let $T_{17}'$ be the stratum of $\Acal'$ corresponding to the edge between
$Q_1$ and $Q_7$. Note that $\varphi$ induces a surjective morphism from the orbit $T_{17}'$ onto the orbit $T_{15}$ of degree $2$ (see \cite[Proposition 6.4]{gubler:tropical}). 
We deduce that 
\begin{equation} \label{coveringdegree}
\varphi_*(\overline{T_{17}'})=2\,\overline{T_{15}}.
\end{equation}
The projection formula in Proposition \ref{properties of refined intersection}, \eqref{pullbackd1} and \eqref{coveringdegree} show that
$$D_5 \cdot \overline{T_{15}}= (D_7'+D_9'+D_{17}'+D_{19}') \cdot \overline{T_{17}'} .$$
The combinatorial nature of the triangulation shows that
$D_9',D_{17}',D_{19}'$ do not intersect $\overline{T_{17}'}$;
since we deal with normal crossing divisors we have
$\overline{T_{17}'} = D_7' . Y_1'$, so
$$D_5 \cdot \overline{T_{15}}= D_7' \cdot \overline{T_{17}'} 
= D_7' \cdot D_7'  \cdot Y_1' = 4D_7'\cdot D_7'\cdot D_1',$$
where the last equality is from~\eqref{assoweil}.
Using Kolb's  relation (b) in Proposition~\ref{Kolb's relations}, we have
$$D_7' \cdot D_7' \cdot D_1'=  - D_7' \cdot (D_2'+D_{12}'+D_{17}') \cdot D_1' 
= -D_7' \cdot D_2' \cdot D_1',$$
where we have used again the combinatorial nature of the triangulation in
the last step.  By~\eqref{assoweil} we have 
$-4D_7' \cdot D_2' \cdot D_1' = - D_7'\cdot D_2'\cdot Y_1'$, so
\begin{equation} \label{final intersection number}
D_5 \cdot \overline{T_{15}} = -D_7' \cdot D_2' \cdot Y_1' = -1,
\end{equation}
where in the last step we have used that $D_2'+D_7'$ restricts to a normal crossing divisor on $Y_1'$. 

\smallskip
\noindent {\it Conclusion:} We have shown that
$-m(\div(F),e_{15}) = 1 = m(\tau(f), e_{15})$.
\smallskip

\begin{rem} \label{rem:in.terms.of.slopes}
  Equation~\eqref{final intersection number} says that
  $\alpha(P_5, e_{15}) = 1$; by Lemma~\ref{eq:deg.alphas} we have
  $\alpha(P_1, e_{15}) = 1$ as well.  Therefore the ``weighted midpoint'' 
  \[ m = \frac 12 (P_1 + P_5) \]
  of~\eqref{eq:m.T} makes sense, and is equal to the midpoint of
  $e_{15}$.  We have $F(P_1) = 0$ and $F(P_5) = \frac 12$, so
  $F(m) = \frac 14$.  Since $v(\pi) = \frac 12$ and
  $F(P_2) = F(P_4) = 0$, we have
  \[\begin{split}
    \slope(F; e_{15}, \Delta_{125}) 
    &= \frac 1{v(\pi)} \big(F(P_2) - F(m)\big) = -2\cdot\frac 14 = -\frac 12 \\
    \slope(F; e_{15}, \Delta_{145}) 
    &= \frac 1{v(\pi)} \big(F(P_4) - F(m)\big) = -2\cdot\frac 14 = -\frac 12. 
  \end{split}\]
  Thus we have
  \[ m(\div(F),e_{15}) = 
  \slope(F; e_{15}, \Delta_{125}) + \slope(F; e_{15}, \Delta_{145})
  = -\frac 12 - \frac 12 =- 1, \]
  as above.  Notice that the slopes are not integers in this case.
\end{rem}

\section{The Sturmfels--Tevelev formula} \label{Alterations}

The original Sturmfels--Tevelev multiplicity formula 
relates tropical multiplicities of maximal cones of tropicalizations of
closed subvarieties of tori under a torus homomorphism.  
It is proved in~\cite[Theorem~1.1]{sturmfels_tevelev:elimination} for
fields with a trivial valuation and
in~\cite[Corollary~8.4]{bpr:trop_curves} in general.
A ``skeletal'' variant was proved for a smooth curve embedded as a
closed subscheme of torus in~\cite[Corollary~6.9]{bpr:trop_curves}.
In the special case of a trivially valued field in characteristic $0$, a higher
dimensional variant was also proved by Cueto \cite[Theorem 2.5]{cueto:implicitization}.

In this section we will prove a generalization of the ``skeletal'' variant which
works in any dimension and also for varieties equipped with a map to a
torus which is generically finite onto its image, but not necessarily a
closed immersion.  As our formula is formally very similar to the ones
mentioned above, we also call it a Sturmfels--Tevelev multiplicity
formula. 

\paragraph \label{par:st.setup}
We fix a strictly semistable pair $(\sX,H)$. Let $T = \Spec(K[M])$ be an algebraic torus, let $N = \Hom(M,\Z)$ be the
group of one-parameter subgroups of $T$, and let 
$\trop:T^\an\to N_\R$ be the tropicalization map, as in~\parref{par:tropicalization}. 
Let $U = X\setminus\supp(H)$ and let 
$\phi:U\to T$ be a morphism.  
Let $U'\subset T$ be the schematic image of $\phi$.  

\begin{prop} \label{cor:map.to.torus}
 The map $\trop\circ\phi:U^\an\to N_\R$ factors through the retraction
  $\tau:U^\an\to S(\sX,H)$, and the restriction of $\trop\circ\phi$ to any
  canonical polyhedron of $S(\sX,H)$ is an integral $\Gamma$-affine map.
  Moreover, $\Trop(U') = \trop\circ \varphi (S(\sX,H))$.
\end{prop}

\pf Choosing a basis for $M\cong\Z^n$, we may
write $\phi$ as a tuple $(\phi_1,\ldots,\phi_n):U\to\bG_{m,K}^n$; we may regard
each $\phi_i$ as a non-zero rational function on $X$ such that
$\supp(\div(\phi_i))\subset\supp(H)_\eta$.  
The first assertions follow by applying 
Proposition~\ref{piecewise linear} to each $\phi_i$.  
The difficulty in the final assertion is that the map $\phi:U\to U'$ needs
not be surjective, but it follows from~\cite[Lemma~4.9]{gubler:forms} that $\Trop(U')=\trop(\varphi(\Uan))$.   
We conclude that $\Trop(U')=\trop \circ \varphi (\tau(\Uan))=\trop\circ \varphi (S(\sX,H))$.
\qed 

\paragraph \label{par:st.setup2}
Suppose now that $\phi:U\to U'$ is generically finite, so 
$d \coloneq \dim(X) = \dim(U')$.  We denote the degree of this map by $[U:U']$. 
As explained in~\parref{par:tropicalization}, 
$\Trop(U')$ is the support of an integral $\Gamma$-affine polyhedral complex $\Sigma_1$
of pure dimension $d$.  Recall that for every maximal (i.e.\ $d$-dimensional)
polyhedron $\Delta\in\Sigma_1$ we have defined a tropical multiplicity
$m_\Trop(\Delta)\in\N \setminus \{0\}$.

Let $\phi_\aff:S(\sX,H)\to\Trop(U')$ denote the restriction of
$\trop\circ\phi$ to $S(\sX,H)$. This is an integral $\Gamma$-affine
map on each canonical polyhedron of $S(\sX,H)$.
Consider a polyhedron $\Delta \in \Sigma_1$ of dimension $d=\dim(X)$.
Choose a $\Gamma$-rational point $\omega \in \relint(\Delta)$ not
contained in a polyhedron $\phi_\aff(\Delta_S)$ of dimension $<d$ for
any canonical polyhedron $\Delta_S$ of $S(\sX,H)$.
Such points are dense in $\relint(\Delta)$.  Clearly $\phi_\aff\inv(\omega)$ is
finite, with each point contained in (the relative interior of) a unique canonical polyhedron
$\Delta_S$ of dimension $d$.  Let $\Delta_S$ be such a canonical
polyhedron.  The image of $\Delta_S$ under $\phi_\aff$ is contained in the
affine span of $\Delta$, so we get a lattice index $[N_\Delta:N_{\Delta_S}]$ 
as in~\parref{par:polyhedra}.

\begin{thm}[Skeletal Sturmfels--Tevelev multiplicity formula]
 \label{Sturmfels--Tevelev multiplicity formula for alterations}
  Using the notations and hypotheses above, we have the identity
  \[ [U:U']\, m_\Trop(\Delta) = \sum_{\Delta_S}  [N_\Delta:N_{\Delta_S}], \]
  where the sum ranges over all canonical polyhedra $\Delta_S$ of the
  skeleton $S(\Xcal,H)$ with 
  $\relint(\Delta_S) \cap \phi_\aff^{-1}(\omega) \neq \emptyset$.
\end{thm}

\proof The proof is similar to the proof of Theorem 8.2 in
\cite{bpr:trop_curves}.  To simplify the notation we set
$\trop_\phi\coloneq\trop\circ\phi:U^\an\to N_\R$.
Consider the affinoid space
$U'_\omega \coloneq \trop^{-1}(\omega)\cap U'^\an$.  
By Proposition~\ref{cor:map.to.torus} the map $\trop_\phi$ factors through
the retraction to the skeleton, so
$X_\omega\coloneq \phi^{-1}(U'_\omega) = \trop_\phi\inv(\omega)$ is the finite
disjoint union of the analytic domains 
$X_{\omega'}= \tau\inv(\omega')$ with $\omega'$ ranging over the finite set 
$\phi_\aff\inv(\omega)$. 
By Corollary~\ref{torus corollary} each $X_{\omega'}$ is affinoid, hence $X_\omega$ is affinoid.

We claim that $\phi:X_\omega\to U'_\omega$ is finite.
By~\cite[Corollary 2.5.13]{berkovich:analytic_geometry}, it suffices to show that the boundary $\del(X_\omega/U'_\omega)$ is empty. 
We have
$\del(U^\an/U'^\an)=\emptyset$ because the analytification of any scheme
is boundaryless~\cite[Theorem 3.4.1]{berkovich:analytic_geometry} (or \emph{closed} in Berkovich's terminology); therefore $\del(X_\omega/U'_\omega)=\emptyset$ by
pullback~\cite[Proposition~3.1.3]{berkovich:analytic_geometry}.
This proves that $\phi:X_\omega\to U'_\omega$ is finite, so the induced
morphism of canonical models $\varphi_\omega: \fX_\omega \to \fU_\omega'$ is
finite by~\cite[Theorem 3.17 and Proposition~3.13]{bpr:trop_curves}.

The canonical model $\fX_\omega$ of $X_\omega$ is the disjoint union of
the canonical models $\fX_{\omega'}$ of $X_{\omega'}$ for 
$\omega'\in\phi_\aff\inv(\omega)$. 
By Corollary~\ref{torus corollary}, the special fibre of 
$\fX_{\omega'}$ is isomorphic to $\bG_{m,\td K}^d$.
By the projection formula~\cite[Proposition 4.5]{gubler:local_heights}
applied to the Cartier divisor $\Div(\nu)$ on $\fU'_\omega$ for some
non-zero $\nu$ in the maximal ideal $K^\cc$ of the valuation ring of $K$, we get
\begin{equation} \label{pre-multiplicity-formula}
\deg(\varphi_\omega) = \sum_{\omega' \to Y} [(\fX_{\omega'})_s: Y]
\end{equation}
for every irreducible component $Y$ of the special fibre of
$\fU'_\omega$, where $\omega'$ ranges over all elements in
$\phi_\aff\inv(\omega)$ with 
$\varphi((\fX_{\omega'})_s)=Y$. By~\cite[Lemma~8.3]{bpr:trop_curves}, 
since $X_\omega\to U'_\omega$ is finite we have
$\deg(\varphi_\omega)=[U:U']$.  

Let $\fU'^\omega$ be the polyhedral formal model 
of 
 $U'_\omega$ as  in~\cite[Definition 4.14]{bpr:trop_curves}. It is the closure of $U'_\omega$ in the canonical 
model of $\trop^{-1}(\omega)$. 
Its special fibre  is 
$\In_\omega(U')$. 
We have a canonical finite surjective morphism
$\fU'_\omega\to\fU'^\omega$ which is an isomorphism on generic fibres (see \cite[Corollary 3.16]{bpr:trop_curves}).
Since the special fibre of $\fU'_\omega$ is reduced, for each irreducible
component $Z$ of $\In_\omega(U')$ we have
\begin{equation} \label{eq:pmf2}
  \sum_{Y \to Z} [Y:Z]= m_Z(\In_\omega(U')), 
\end{equation}
where the sum runs over all irreducible components $Y$ of
$(\fU'_\omega)_s$ mapping onto $Z$, and $m_Z(\In_\omega(U'))$ is the
multiplicity of $Z$ in $\In_\omega(U')$.  
This follows again from the projection formula; see~\cite[3.34(2)]{bpr:trop_curves}. 

Composing $\phi_\omega:\fX_\omega\to\fU'_\omega$ with
$\fU'_\omega\to\fU'^\omega$ gives a finite surjective morphism
$\fX_\omega\to\fU'^\omega$.  As explained in~\parref{par:tropicalization},
(the reduced scheme underlying) an irreducible
component $Z$ of $\In_\omega(U')$ is isomorphic to the multiplicative
torus of rank $d$ over $\td{K}$. 
 As the same is true for $(\fX_{\omega'})_s$ for
$\omega'\in\phi_\aff\inv(\omega)$, one checks as in the proof
of~\cite[Corollary~8.4]{bpr:trop_curves} that 
\begin{equation} \label{eq:pmf3}
  [(\fX_{\omega'})_s:Z] = [N_\Delta:N_{\Delta_{S(\omega')}}]
\end{equation}
for every $\omega'$ with $(\fX_{\omega'})_s$ lying over $Z$, where
$\Delta_{S(\omega')}$ is the unique canonical polyhedron of 
$S(\Xcal,H)$ with $\omega' \in \relint(\Delta_{S(\omega')})$.

Recall that $m_\Trop(\Delta)$ is the 
number of irreducible components of $\In_\omega(U')$ counted with
multiplicities.  
Since $\deg(\phi_\omega) = [U:U']$ we have
\begin{equation} \label{mc}
[U:U']\,m_\Trop(\Delta) =  \sum_Z \deg(\varphi_\omega)\, m_Z(\In_\omega(U')),
\end{equation}
where $Z$ runs over all irreducible components of $\In_\omega(U')$. 
Combining this with~\eqref{pre-multiplicity-formula}, \eqref{eq:pmf2},
and~\eqref{eq:pmf3} leads to 
 $$[U:U']m_\Trop(\Delta) =  \sum_Z \sum_{Y \to Z} \sum_{\omega' \to Y} 
 [(\fX_{\omega'})_s: Y]\,[Y:Z]=\sum_{\omega'} [N_\Delta:N_{\Delta_{S(\omega')}}],  $$
where  $Y$ ranges over all irreducible components of $(\fU'_\omega)_s$ lying over $Z$ and where $\omega'$ ranges over all elements in
$\phi_\aff\inv(\omega)$ with 
$\varphi((\fX_{\omega'})_s)=Y$. 
Since $\omega' \mapsto \Delta_{S(\omega')}$ is a bijection from
${\varphi}_{\rm aff}^{-1}(\omega)$ onto the set of canonical polyhedra
$\Delta_S$ of $S(\Xcal,H)$ with 
$\relint(\Delta_S) \cap {\varphi}_{\rm aff}^{-1}(\omega) \neq \emptyset$, 
we get the claim.  \qed

\begin{rem}
  It follows from the considerations in~\parref{par:st.setup2} that when $\phi$ is generically
  finite onto its image then $S(\sX,H)$ necessarily has dimension 
  $d = \dim(X)$.  This non-trivial condition on the strictly semistable
  pair $(\sX,H)$ is not obvious from the definitions.  As
  $\Trop(U')$ has pure dimension $d$ and is connected in codimension one,
  one might wonder if there exist natural additional conditions on $\phi$
  which guarantee that $S(\sX,H)$ has the same properties.
\end{rem}

\section{Faithful tropicalization}

In this section we fix a strictly semistable pair $(\sX,H)$.
As always we use the associated notation~\parref{notn:ssp.notation}.
For $n\geq 0$ we let
$\trop:\bG_{m,K}^{n,\an} = \Spec(K[x_1^{\pm 1},\ldots,x_n^{\pm n}])^\an\to\R^n$
denote the tropicalization map as in~\parref{par:tropicalization}, defined by
\[ \trop(p) = \big(-\log|x_1(p)|,\ldots,-\log|x_n(p)|\big). \]

Our goal is to prove that there is a
rational map $\phi$ from $X = \sX_\eta$ to a torus $T\cong\bG_{m,K}^n$
which takes $S(\sX,H)$ isomorphically onto its image.  Note that 
a rational map $\phi$ is always defined on $S(\sX,H)$ since the points
of $S(\sX,H)$ are norms on the function field of $X$: see
Remark~\ref{rem:abhyankar}.  By ``isomorphically'' we mean that we want $\phi$ to be
injective on $S(\sX,H)$ and we want it to preserve the integral affine
structure.  Roughly, in this situation one ``sees'' the entire skeleton 
in the tropicalization of (the image of) $X$; this is an important
compatibility between the intrinsic and embedded polyhedral structures of
$X$.  (See however Remark~\ref{rem:bigger.tropicalization}.)

\smallskip
We start with the following basic property which is an immediate consequence of Proposition \ref{lem:piecewise.affine}. 

\begin{prop} \label{rational map on the skeleton}
Let $f = (f_1,\ldots,f_n):X\dashrightarrow\bG_{m,K}^n$ be a rational map. Then $S(\sX,H)$ can be covered by 
 finitely many integral $\Gamma$-affine polyhedra $\Delta$ such that
  $\trop\circ f|_\Delta$ is an  integral $\Gamma$-affine map.
\end{prop}

\begin{defn}
A rational map $f = (f_1,\ldots,f_n):X\dashrightarrow\bG_{m,K}^n$ is said
to be \emph{unimodular} on a canonical polyhedron $\Delta_S$ of $S(\sX,H)$ provided that 
$\Delta_S$ can be covered by finitely many integral $\Gamma$-affine polyhedra $\Delta$ such that
$\trop\circ f|_\Delta$ is a unimodular integral $\Gamma$-affine map on $\Delta$  (see \ref{par:polyhedra}). 
  We call $f$ {\it unimodular on $S(\sX,H)$} if $f$ is unimodular on any canonical polyhedron of $S(\sX,H)$. 
We say that $f$ is a \emph{faithful tropicalization} of $S(\sX,H)$ if $f$ is 
  unimodular and $\trop\circ f$ is injective on $S(\sX,H)$.
\end{defn}

\smallskip
The next lemma is essentially~\cite[Lemma~6.17]{bpr:trop_curves}.

\begin{lem} \label{lem:add.functions}
  Let $f_1,\ldots,f_n,g$ be non-zero rational functions on $X$, and suppose
  that the rational map $f=(f_1,\ldots,f_n):X\dashrightarrow\bG_{m,K}^n$ is
  unimodular on the canonical polyhedron $\Delta_S$ of $S(\sX,H)$.  Then 
  $(f_1,\ldots,f_n,g):X\dashrightarrow\bG_{m,K}^{n+1}$ is also
  unimodular on $\Delta_S$.
\end{lem}

\pf It follows from Proposition~\ref{lem:piecewise.affine} that the skeleton $S(\sX,H)$ 
has a covering by finitely many integral $\Gamma$-affine polyhedra $\Delta$ such that $\trop \circ h|_\Delta$
is an  integral $\Gamma$-affine map $\Delta \rightarrow \R^{n+1}$ for $h\coloneq (f_1,\ldots,f_n,g)$. 
Since $\trop \circ h|_\Delta$ factors through $\trop \circ f|_\Delta$, transitivity of lattice indices shows easily 
that $h$ is unimodular.  
\qed

\begin{prop} \label{prop:unimodularity}
  For every canonical polyhedron $\Delta_S\subset S(\sX,H)$, there
  exist non-zero rational functions $f_1,\ldots,f_n\in K(X)$ such that 
  $\trop\circ(f_1,\ldots,f_n)|_{\Delta_S}:\Delta_S\to\R^n$ is a unimodular integral $\Gamma$-affine map
  (and therefore   injective). In particular, $f=(f_1, \dots ,f_n)$ is unimodular on $\Delta_S$.
\end{prop}

\proof
Let $S$ be the corresponding vertical stratum of $D$. 
Every point of $S$
has a neighbourhood $\sU$ that admits an \'etale morphism
$\psi:\sU\to\sS = \Spec(K^\circ[x_0,\ldots,x_d]/\angles{x_0\cdots x_r-\pi})$
as in~\eqref{eq:local.sss.pair}.  We can shrink
$\sU$ so that $\fU=\hat\sU$ is a building block with distinguished stratum
$S$.  In particular, $S$ is defined by 
$\psi^*(x_0) = \cdots = \psi^*(x_{r+s}) = 0$.
The canonical polyhedron $\Delta_S$ of the skeleton $S(\sX,H)$
is contained in $\fU_\eta\subset\sU^\an$, and
\[ \Val_S(p) =
\big(-\log |\psi^*x_0(p)|,\ldots,-\log |\psi^*x_{r+s}(p)|\big) \]
maps $\Delta_S$ homeomorphically onto $\Delta(r,\pi)\times\R_+^s$ by~\parref{par:skel.building.block}.  In fact, the structure of
integral $\Gamma$-affine polyhedron on $\Delta_S$ is defined by the map $\Val_S$,
so $\Val_S|_{\Delta_S}$ is by definition a unimodular integral $\Gamma$-affine map.
Interpreting $\psi^*(x_0),\ldots,\psi^*(x_{r+s})$ as rational functions on
$X$ and $\Val_S$ as the composition of 
$(\psi^*(x_0),\ldots,\psi^*(x_{r+s})):X\dashrightarrow\bG_{m,K}^{r+s+1}$
with $\trop:\bG_{m,K}^{r+s+1,\an}\to\R^{r+s+1}$, we obtain the claim. \qed

\begin{thm} \label{thm:faithful.ss}
  Let $(\sX,H)$ be a strictly semistable pair.
  Then there exists a finite collection $f_1,\ldots,f_n$ of
  non-zero rational functions on $X$ such that the associated
  rational map $X\dashrightarrow\bG_{m,K}^n$ is a faithful tropicalization
  of $S(\sX,H)$.
\end{thm}

\pf By Proposition~\ref{prop:unimodularity} and
Lemma \ref{lem:add.functions}, 
we can find $f_1,\ldots,f_n\in K(X)^\times$ such that the rational map 
$f=(f_1,\ldots,f_n)$ is unimodular
on every canonical polyhedron $\Delta_S$: indeed, we may take
any collection $(f_1,\ldots,f_n)$ which includes all rational functions from
Proposition~\ref{prop:unimodularity} for each $\Delta_S$. 
By construction, $\trop \circ f$ is injective on every $\Delta_S$. 
It remains to enlarge the collection
$(f_1,\ldots,f_n)$ so that $\trop\circ f$ is injective on
$S(\sX,H)$. 

By Chow's lemma \cite[Theorem 5.6.1]{egaII}, there is a  birational surjective morphism 
$\varphi: \sX' \rightarrow \sX$ for a projective variety $\sX'$ over $\kcirc$. There are open 
dense subsets $\sU$ of $\sX$ and $\sU'$ of $\sX'$ such that $\varphi$ restricts to an 
isomorphism $\sU' \rightarrow \sU$. For simplicity, we use this to identify $\sU$ with $\sU'$ 
and hence we have an identification $K(X)=K(X')$ of the function fields of the
generic fibres $X,X'$ of $\sX$ and $\sX'$. 
Since any element of the skeleton is an Abhyankar point (see Remark \ref{rem:abhyankar}), we have 
$S(\sX,H) \subset \sU_\eta^\an = \sU_\eta'^\an$. 

We have seen in \ref{stratification} that a strictly semistable pair has a
canonical stratification $\str(\sX_s,H)$ of the special fibre $\sX_s$. The
preimage $\varphi^{-1}(S)$ of $S \in \str(\sX_s,H)$ is not necessarily
irreducible, but it contains only finitely many generic points.  Let $\Fcal$ be
the collection of all such generic points for all strata $S$. Note that $\Fcal$
is a finite subset of $\sX'$.  Since $\sX'$ is projective over the affine scheme
$\Spec(K^\circ)$, any two points of
$\sX'$ are contained in a common affine open subset.  Using that $\sX'$ is
quasicompact, we conclude that for every $\zeta' \in \Fcal$, there are finitely
many affine open subsets $\sU_{\zeta'j}'$ containing $\zeta'$ and covering
$\sX'$. On every such $\sU_{\zeta'j}'$, there exists a finite collection of regular
functions $f_{\zeta'jk} \in \Ocal(\sU_{\zeta'j}')$ whose reductions to the
special fibre have zero set $\bar{\zeta'} \cap (\sU_{\zeta'j}')_s$. This means
that for every $x' \in (X')^{\rm an}$ with reduction
$\red_{\sX'}(x') \in (\sU_{\zeta'j}')_s$, we have $|f_{\zeta'jk}(x')|=1$ for
some $k$ if $\red_{\sX'}(x') \not \in \bar{\zeta'}$ and $|f_{\zeta'jk}(x')|<1$
for all $k$ if $\red_{\sX'}(x') \in \bar{\zeta'}$. Note that we have
$|f_{\zeta'jk'}(x')| \leq 1$ for all $k'$.

These finitely many  functions $f_{\zeta'jk}$ may be viewed as rational functions on $X$ and we add them to the 
collection $(f_1, \dots ,f_n)$ considered at the beginning. We claim that the resulting tropicalization is faithful. 
As remarked above, it is enough to show that this tropicalization is injective. We consider points $x \neq y$ from $S(\sX,H)$ and 
we have to look for a function from our extended collection such that the absolute value of the function separates $x$ and $y$. 
By Corollary~\ref{cor:stratum-genericpt}, $x$ is contained in the relative
interior of a unique canonical polyhedron $\Delta_S$ and  the reduction
$\red_\sX(x)$ is the generic point $\zeta_S$ of the vertical stratum $S \in
\str(\sX_s,H)$. A similar statement holds for $y$ and we denote the
corresponding vertical stratum by $T$ and reduction by $\zeta_T$.

If $T \subset \bar{S}$, then $\Delta_S$ is a closed face of $\Delta_T$ by \ref{canonical polyhedra}. We conclude that $x,y \in \Delta_T$ and $x=y$ follows from injectivity of $\trop \circ f$ on $\Delta_T$. So we may assume that $T \not \subset \bar S$ and hence $\zeta_T \not \in \bar{S}$. Using that $x,y \in S(\sX,H) \subset \sU_\eta^\an = \sU_\eta'^\an$, we may view $x$ and $y$ as points of $(X')^\an$. Note that 
$$\varphi \circ \red_{\sX'}(x)=\red_{\sX}(x)=\zeta_S \in S.$$
 There is a generic point $\zeta'$ of $\varphi^{-1}(S)$ with $\red_{\sX'}(x) \in \bar{\zeta'}$. By construction, there is a $j$ such that $\red_{\sX'}(x) \in \sU_{\zeta'j}'$. From $\zeta_S=\varphi\circ \red_{\sX'}(x)$, we deduce that $\zeta_S = \varphi \circ \red_{\sX'}(\zeta')$. Similarly, we show
$$\varphi \circ \red_{\sX'}(y)=\red_{\sX}(y)=\zeta_T \not\in \bar S$$
and hence $\red_{\sX'}(y) \not \in \bar{\zeta'}$. We conclude from the above considerations that $|f_{\zeta'jk}(y)|=1$ for some $k$. Using $\red_{\sX'}(x) \in \bar{\zeta'}$, we have  $|f_{\zeta'jk}(x)|<1$. This proves the claim. \qed

\begin{rem} \label{rem:bigger.tropicalization}
  In the statement of Theorem~\ref{thm:faithful.ss} it is not assumed that
  the divisor of each $f_i$ has support contained in $\supp(H)$.  In
  particular, the tropicalization map will not generally factor through
  retraction to the skeleton, so its image in $\R^n$ may be much larger
  than the image of $S(\sX,H)$.
\end{rem}

\begin{eg}
  When $\dim(X) = 1$ the skeleton $S(\sX,H)$ is a metric graph.
  Theorem~\ref{thm:faithful.ss} says that there exists a rational map $\phi$
  from $X$ to a torus whose restriction to $S(\sX,H)$ is an isometry onto its
  image, where the metric on the image is defined by the lattice length.
  Increasing the dimension of the torus, we may even assume that $\phi$ is a
  closed embedding on an open subscheme.  Therefore this extends the faithful
  tropicalization result of~\cite[Theorem~6.22]{bpr:trop_curves}, as well as
  considerably simplifying its proof.  Indeed, as $\sX$ is automatically
  projective when $\dim(X)=1$, the Chow's lemma argument used in
  Theorem~\ref{thm:faithful.ss} is not needed, so that the proof given in the
  present paper is shorter and
  more conceptual than the one in~\cite{bpr:trop_curves}.
\end{eg}

\section{Sections of tropicalizations}

\paragraph\label{par:setting-trop}\label{residue-norm}
 Let $\varphi: U \hookrightarrow T$ be a
closed immersion of a very affine variety $U$ into an algebraic $K$-torus
$T=\Spec(K[M])$.  As $U$ is a variety, it is an integral scheme.  We set
$\trop_\phi = \trop\circ\phi^\an:U^\an\to N_\R$, where $N = \Hom(M,\Z)$, so
$\Trop(U) = \trop_\phi(U^\an)$ is the corresponding tropicalization as defined
in section~\parref{par:tropicalization}.  
Fixing a basis of the character group $M$ of $T$, we identify $T$ with
the torus $\Spec(K[x_1^{\pm 1},\ldots,x_n^{\pm 1}])$ and $M_\R,N_\R$ with
$\R^n$.  Let $\fa$ be the ideal defining $U$ as a closed subscheme of $T$.  Fix
$\omega = (\omega_1, \ldots, \omega_n) \in \Trop(U)$, and put
$r_i = \exp(- \omega_i) \in \rdop$.  Then $U_\omega$ is the spectrum of the
affinoid algebra
\[ \sA_\omega = K\angles{r^{-1} x, r x^{-1}} / \mathfrak{a}
K\angles{ r^{-1} x, r x^{-1} }, \]
which carries the residue norm
\begin{equation}\label{eq:resnorm}
  \|f\|_\res = \inf_{\pi_\omega(g) = f} \|g\|_r = \min_{\pi_\omega(g) = f} \|g\|_r,
\end{equation}
where $\pi_\omega:K\angles{r^{-1} x, r x^{-1}}\to\sA_\omega$ is the quotient map
and $\| \scdot \|_r$ denotes the spectral norm on
$K\angles{ r^{-1} x, r x^{-1} }$, i.e.
\[\bigg\| \sum_{I \in \mathbf{Z}^n} a_I x^I\bigg\|_r = \max_I \big\{
|a_I|\,r^I\big\}. \]
When $\omega \in N_\Gamma$, then $\sA_\omega$ is strictly affinoid, so
the infimum in~\eqref{eq:resnorm} is a minimum by~\cite[Corollary~5.2.7/8]{bgr:nonarch}. In general, we can use 
the base change trick from \cite[Proof of Proposition 2.1.3]{berkovich:analytic_geometry} to reduce to the strictly 
affinoid case.

\paragraph\label{par:Shilov-setting}
Now let 
 $\omega\in\Trop(U)\cap N_\Gamma$. Then 
the analytic domain $U_\omega = \trop_\phi\inv(\omega) = \sM(\sA_\omega)$ is
strictly affinoid. 
Let $\fU^\omega$ be the \emph{ polyhedral formal model} of $U_\omega$  defined by 
$$\fU^\omega = \Spf\left(K\angles{r^{-1} x, r x^{-1}}^\circ / \left(\mathfrak{a} K\angles{r^{-1} x, r x^{-1}} \cap K\angles{r^{-1} x, r x^{-1}}^\circ \right)    \right)
.$$
It is an admissible formal scheme over $\kcirc$ with 
special fibre equal to the initial degeneration $\inn_\omega(U)$ (see \cite[Proposition 4.17]{bpr:trop_curves}). 
Let
$\fU_\omega^\can = \Spf(\sA_\omega^\circ)$ be the \emph{ canonical model} of the
affinoid space $U_\omega$, as defined in~\parref{par:notation}; its special
fibre is the \emph{ canonical reduction} of $U_\omega$.
As noted in the remark
after~\cite[Proposition 4.17]{bpr:trop_curves}, we have a canonical finite
morphism $\fU_\omega^\can\to\fU^\omega$, which is the identity on generic
fibres. 

By~\cite[Proposition~2.4.4(iii)]{berkovich:analytic_geometry} and its proof, the
supremum seminorm $\|\scdot\|_{\sup}$ on $\sA_\omega$ is equal to the maximum of
the norms contained in the Shilov boundary of $U_\omega$, which correspond to
the irreducible components of the canonical reduction
$\Spec(\sA_\omega^\circ\tensor_{K^\circ}\td K)$ of $U_\omega$.  The canonical
reduction is always a reduced scheme.  Therefore, the following
are equivalent:
\begin{enumerate}
\item $U_\omega$ has a unique Shilov boundary point;
\item the canonical reduction of $U_\omega$ is an integral scheme; and
\item the supremum seminorm $\|\scdot\|_{\sup}$ on $\sA_\omega$ is multiplicative.
\end{enumerate}
In particular, when the above conditions hold then the supremum seminorm is the
unique Shilov boundary point of $U_\omega$.
We will leverage this fact to construct a section of $\trop_\phi$ on the locus
$\{\omega\in\Trop(U)~:~m_\Trop(\omega) = 1\}$.

\begin{lem}\label{lem:unique_shilov}
  Let $\phi:U\inject T$ be a closed immersion of a very affine variety $U$ into
  an algebraic $K$-torus $T = \Spec(K[M])$.
  For $\omega\in\Trop(U)\cap N_\Gamma$, if $m_\Trop(\omega) = 1$ then 
  $U_\omega = \trop_\phi\inv(\omega)$ contains a unique Shilov boundary point.
\end{lem}

\proof  As noted above, we have a canonical finite
morphism $\fU_\omega^\can\to\fU^\omega$, which is the identity on generic
fibres. By~\cite[Corollary~3.16]{bpr:trop_curves}, the canonical morphism induces a surjective finite map between the equidimensional special fibres
$(\fU_\omega^\can)_s = \Spec(\sA_\omega^\circ\tensor_{K^\circ}\td K)$ and $\inn_\omega(U)$.
By hypothesis, $\inn_\omega(U)$ is irreducible and generically reduced. Therefore, by the projection
formula~\cite[3.34(2)]{bpr:trop_curves}, there is a unique irreducible component
of $(\fU_\omega^\can)_s$, and it maps birationally onto the reduction of
$\inn_\omega(U)$.  We conclude that $(\fU_\omega^\can)_s$ is an integral scheme and the equivalence of (1) and (2) above shows the claim.  \qed

\begin{rem}
  In Lemma~\ref{lem:unique_shilov} it is necessary to assume that $U$ is
  irreducible, or at least equidimensional.  As a counterexample, let $U$ be the
  closed subscheme of $T = \Spec(K[x_1^{\pm1}, x_2^{\pm 1}])$ given by the ideal
  $\fa = \angles{x_1-1,x_2-1}\cap\angles{x_1-1-\varpi}$ for $\varpi\in K^\times$
  with $|\varpi| < 1$.  Then $U$ is the disjoint union of the line
  $\{x_1 = 1+\varpi\}$ with the point $(1, 1)$.  The initial degeneration at
  $\omega=0$ is defined by the ideal
  $\inn_w(\fa) = \angles{(x_1-1)^2,(x_1-1)(x_2-1))}$ over $\td K$.  This is a
  generically reduced line with an associated point at $(1,1)$.  It has
  tropical multiplicity $1$, but the canonical reduction is the
  disjoint union of a point and a line, so that $U_\omega$ has \emph{two} Shilov
  boundary points.
\end{rem}

\begin{prop}\label{prop:resnorm} For $\omega \in \Trop(U) \cap N_\Gamma$, the following are equivalent:
\begin{itemize}
 \item[(a)] $\|\scdot\|_\res = \|\scdot\|_{\sup}$;
 \item[(b)] $\inn_\omega(U)$ is a
  reduced scheme.
\end{itemize}
If we assume additionally that $m_\Trop(\omega) = 1$, i.e.\ if $\inn_\omega(U)$ is irreducible and
  generically reduced, then (a) and (b) are also equivalent to:
\begin{itemize}
\item[(c)] the unique Shilov boundary point of $U_\omega$ is  equal to
  $\|\scdot\|_\res$;
\item[(d)] $\|\scdot\|_\res$ is multiplicative;
\item[(e)] $\inn_\omega(U)$ is an integral scheme;
\end{itemize}
 \end{prop}

\proof The polyhedral formal model $\fU^\omega$ is $\Spf(B_\omega)$ for
$B_\omega =
\pi_\omega(K\angles{r^{-1} x, r x^{-1}}^\circ)\subset\sA_\omega^\circ$.
For $f\in\sA_\omega$ we have $f\in B_\omega$ if and only if there exists
$g\in K\angles{r^{-1} x, r x^{-1}}^\circ$ such that $\pi_\omega(g) = f$, which
is true if and only if $\|f\|_{\res}\leq 1$.  Since $\|f\|_{\sup}$ and
$\|f\|_\res$ are in $|K^\times|$ in any case, we have
$\|\scdot\|_{\res} = \|\scdot\|_{\sup}$ if and only if
$\sA_\omega^\circ = B_\omega$.  Now the equivalence of (a) and (b) follows from~\cite[Proposition~3.18]{bpr:trop_curves}.

If  $m_\Trop(\omega) = 1$, then Lemma \ref{lem:unique_shilov} shows that $U_\omega$ has a unique Shilov boundary point. 
It follows from \ref{par:Shilov-setting} that  $\|\scdot\|_{\sup}$ is multiplicative and equal to the unique Shilov boundary point. 
This proves the equivalence of (a), (c) and (d). Since $\inn_\omega(U)$ is irreducible, the equivalence of (b) and (e) is obvious. 
\qed

\smallskip

The next result was proved for compact subsets of curves
in~\cite[Theorem~6.24]{bpr:trop_curves}. 
We prove it here in a very general situation.

\begin{thm}\label{thm:section}
  Let $\varphi: U \hookrightarrow T$ be a closed immersion of a very affine
  (integral) variety $U$ into an algebraic $K$-torus $T=\Spec(K[M])$, and let
  $\trop_\phi = \trop\circ\phi^\an:U^\an\to N_\R$.  Let $Z \subset \Trop(U)$ be a
  subset such that $m_\Trop(\omega) = 1$ for all $\omega \in Z$. Then for every
  $\omega \in Z$, the affinoid space $U_\omega=\trop_\varphi^{-1}(\omega)$ has a
  unique Shilov boundary point $s(\omega)$, and $\omega\mapsto s(\omega)$
  defines a continuous partial section $s:Z\to U^\an$ of the tropicalization map
  $\trop_\varphi : \Uan \rightarrow \Trop(U)$ on the subset $Z$.  Moreover, if
  $Z$ is contained in the closure of its interior in $\Trop(U)$, then $s$ is the
  unique continous section of $\trop_\phi$ defined on $Z$.
 
\end{thm}

\proof
First we prove that $\trop_\phi\inv(\omega)$ has a unique Shilov boundary point
for all $\omega\in Z$.  When the valuation map 
$v: K^\times \rightarrow \mathbf{R}$ is surjective this follows from
Lemma~\ref{lem:unique_shilov}.  
In the general case let $L$ be a non-archimedean extension field of $K$
such that the valuation map $L^\times \rightarrow \mathbf{R}$ is
surjective. Let $U_L$ denote the base change of $U$ to $L$, 
and let $p: U_L \rightarrow U$ be the projection. 
We have $\Trop(U) = \Trop(U_L)$, and the tropical multiplicities in the two
tropicalizations coincide essentially by definition.  See~\cite[Proposition~3.7,
Definition~13.4]{gubler:guide}. 
By the above, the affinoid space
$U_\omega\hat\tensor_K L = \trop_\phi\inv(\omega)\hat\tensor_K L =
\trop_{\phi\hat\tensor L}\inv(\omega) = (U_L)_\omega$
has a unique Shilov boundary point.
It follows directly from the definition of the Shilov boundary that the image of the Shilov boundary of 
$U_\omega\hat\tensor_K L$ with respect to $p^\an$ contains the Shilov boundary of $U_\omega$ as the former 
has more analytic functions than the latter. 
Hence $p^{\an}$  maps the unique Shilov
boundary point $s_L(\omega)$ of $U_\omega\hat\tensor_K L$ to the unique Shilov
boundary point $s(\omega)$ of $U_\omega$.  Clearly
$\omega\mapsto s(\omega)$ is a section of $\trop_\phi$.  Note that we have in
fact shown that the section $s$ respects base extension, in that 
$s = p^{\an} \circ s_L$, where $s_L:Z\to U_L^\an$ is the partial section defined
relative to $L$.  In particular, if $s_L$ is continuous, then $s$ is.

Next we prove continuity and uniqueness when $U = T$ (and $\phi$ is the identity
map).  In this case $s(\omega)$ is the
Gauss norm $\|\scdot\|_r:K[x_1^{\pm 1},\ldots,x_n^{\pm 1}]\to\R_+$
in the notation of~\parref{residue-norm}, where $\omega = (\omega_1,\ldots,\omega_n)$ and 
$r = (\exp(-\omega_1),\ldots,\exp(-\omega_n))$.  It is clear that $s$ is
continuous and is defined on all of $\R^n$; its image is by definition 
the \emph{skeleton} $S(T) = s(\R^n)$ of the torus $T$.
We now turn to uniqueness.
Let $\omega\in \R^n$, and suppose that there exists a continuous section $s':Z\to T^\an$
defined on an open neighbourhood $Z$ of $\omega$ such that $s(\omega)\neq
s'(\omega)$.   Let $u = s(\omega)$.  By hypothesis $u'\coloneq s'(\omega)\neq u$, so
there exists a (non-zero) Laurent polynomial $h = \sum_{I\in\Z^n} a_I x^I\in
K[x_1^{\pm 1},\ldots,x_n^{\pm1}]$ such that $|h(u')| < \|h\|_r$.  Since
$\trop(u') = \omega$ we have $|x_i(u')| = r_i$ for all $i$.  If there were a unique exponent
$I$ such that $|a_I|r^I = \|h\|_r$, then by the ultrametric inequality as
applied to the seminorm corresponding to $u'$, we would have $|h(u')| =
\|h\|_r$; therefore there are at least two exponents $I$ such that $|a_I|r^I$ is
maximal.  In other words, the initial degeneration of $h$ at $\omega$ is not a
monomial, so
$\omega\in\Trop(h)\coloneq\trop(V(h))$, where $V(h)$ is the zero set of $h$.
The maps $w\mapsto|h(w)|$ and $t\mapsto\|h\|_t$ are continuous, so there exists
a small open neighbourhood $W$ of $u'$ in $T^\an$ such that $|h(w)| < \|h\|_t$
for all $w\in W$, where $t = (|x_1(w)|,\ldots,|x_n(w)|)$.  By the above
argument, then, we have $\trop(W)\subset\Trop(h)$.  But
$s'^{-1}(W)\subset\trop(W)$ is an open neighbourhood of $\omega$ and $\Trop(h)$
has codimension one in $\R^n$, a contradiction.
Since the locus where two maps to a Hausdorff space coincide is closed, this
implies that $s$ is the unique continuous section defined on any subset $Z$
which is contained in the closure of its interior.

We treat the general situation by reducing to the case of a torus settled above.
Let $d = \dim(U)$; we may assume $d < n$.
The tropicalization $\Trop(U)$ is the support of an integral $\Gamma$-affine polyhedral complex $\Sigma_1$ of pure dimension $d$. It contains finitely many polyhedral faces $\Delta_i$ of maximal dimension $d$. 
Let $L_i$ be the ($d$-dimensional) linear span of $\Delta_i -v$ for any $v \in \Delta_i$. Arguing by induction, one can show that for all  $n>d$ there exists a $(d \times n)$-matrix with entries in $\mathbf{Z}$ such that  the corresponding linear map $f: \R^n \rightarrow \R^d$ is injective on each $L_i$, and hence on each $\Delta_i$.  Let $\alpha: \G_m^n \rightarrow \G_m^d$ be the homomorphism of tori such that the associated cocharacter map is $f$. Consider the morphism $\psi = \alpha \circ \varphi: U \rightarrow \G_m^d$. The diagram
\begin{equation} \label{eq:reduce.to.torus}
\xymatrix{ 
  {U^{\an}} \ar[r]^(.5){\psi} \ar[d]_{\trop_\varphi} & \G_m^{d, \an} \ar[d]^{\trop} \\
\Trop(U) \ar[r]^(.6){f} & \R^d}
\end{equation}
is commutative, where we also write $\psi$ for $\psi^\an$.
By construction, the map $f$ is finite-to-one on the subset $\Trop(U)$ of $\R^n$. 
Let $S(\G_m^d)$ be the skeleton of the torus $\G_m^{d}$, defined above.   Fix
$\omega' \in \R^d$ with coordinates in $\Gamma$, and write $\{\omega_1,
\ldots,\omega_\ell\} = f^{ -1}(\omega')$.  The affinoid domain
$U'_{\omega'} \coloneq \trop\inv(\omega')\subset\G_m^{d,\an}$ 
has a unique Shilov boundary point, namely the
unique point in the skeleton $S(\G_m^{d})$ mapping to $\omega'$. Now we use a
similar argument as in the proof of the Sturmfels-Tevelev formula
\ref{Sturmfels--Tevelev multiplicity formula for alterations}.
By the commutativity of~\eqref{eq:reduce.to.torus},
$\psi^{-1}( U'_{\omega'})$ is the disjoint union of the
finitely many affinoid subdomains $U_{\omega_i} = \trop_\varphi^{-1}(\omega_i)$ for
$i=1,\ldots,\ell$.  The analytification of $\psi$ is
boundaryless~\cite[Theorem~3.4.1]{berkovich:analytic_geometry}, hence by
pullback 
~\cite[Proposition~3.1.3]{berkovich:analytic_geometry} we find that
$\del(U_{\omega_i}/U'_{\omega'})$ is empty. By
\cite[Corollary~2.5.13]{berkovich:analytic_geometry} this implies that
$U_{\omega_i}\to U'_{\omega'}$ is finite. Therefore the associated map on
reductions $\td U_{\omega_i}\to\td U'_{\omega'}$ is 
finite~\cite[Theorem 6.3.4/2]{bgr:nonarch}, where
$U_{\omega_i} = \sM(\sA_{\omega_i})$ (resp.\ $U'_{\omega'} = \sM(\sA'_{\omega'})$)
and $\td U_{\omega_i} = \Spec(\td \sA_{\omega_i})$ 
(resp.\ $\td U'_{\omega'} = \Spec(\td \sA'_{\omega'})$).
Let $\omega = \omega_i$ for some $i$, and suppose that
$m_\Trop(\omega)=1$.  Then both reductions $\td U_{\omega}, \td U'_{\omega'}$
are irreducible,  $d$-dimensional
schemes over the residue field, so the generic point is mapped to the generic
point. This implies that for such $\omega$, the image of $s(\omega)$ under
$\psi$ lies in the skeleton $S(\G_m^{d})$ and, conversely, if a point in
$U_\omega = \trop_\varphi^{-1}(\omega)$ is mapped to the skeleton $S(\G_m^{d})$ under
$\psi$, then it is equal to $s(\omega)$. In other words,
$\{s(\omega)\} = \trop_\phi\inv(\omega) \cap \psi\inv(S(\G_m^d))$.

Now we prove that $s$ is continuous.  For this we may assume that
$v:K^\times\to\R$ is surjective, as remarked above.
It suffices to show that $s(Z)$ is closed in 
$\trop_\varphi^{-1} (Z)\subset U^{\an}$ (endowed with its relative topology),
since $\trop_\varphi: \trop_\varphi^{-1}(Z) \rightarrow Z$ is a proper map to a
metric space, and a proper map to a metric space is
closed~\cite{palais:proper_maps}.  In fact, since the image of a continuous section of a
continuous map between Hausdorff spaces is necessarily closed, showing
that $s$ is continuous is equivalent to proving $s(Z)$ is closed; in particular,
$S(\G_m^d)$ is closed in $\G_m^{d,\an}$, being the image of the continuous
section of $\trop:\G_m^{d,\an}\to\R^d$.
When the valuation is surjective we have shown that 
$\{s(\omega)\} = \trop_\phi\inv(\omega) \cap \psi\inv(S(\G_m^d))$ for all
$\omega\in Z$, so $s(Z) = \trop_\phi\inv(Z) \cap \psi\inv(S(\G_m^d))$ is
indeed closed in $\trop_\phi\inv(Z)$.

Finally, we prove that $s$ is unique when $Z$ is contained in the closure of its
interior in $\Trop(U)$, no longer under the assumption that the
valuation is surjective. 
Let $s':Z\to U^\an$ be another continuous partial section of $\trop_\phi$.
Let $Z'\subset Z$ be an open subset of $\Trop(U)$
contained in the relative interior of a $d$-dimensional integral $\Gamma$-affine
polyhedron in $\Trop(U)$.  Then $f(Z')$ is open in $\R^d$ and
$f:\Trop(U)\to\R^d$ is injective on $Z'$, so it has an
inverse $g:f(Z')\to Z'$.  Let $\omega\in Z'$ have $\Gamma$-rational coordinates.
Define
$\sigma,\sigma': f(Z')\to\G_m^{d,\an}$ by $\sigma = \psi\circ s\circ g$ and
$\sigma' = \psi\circ s'\circ g$.  These are both partial sections of
$\trop:\G_m^{d,\an}\to\R^d$ defined on $f(Z')$, so by the torus case, they are
equal.  Hence $\psi(s'(\omega)) = \psi(s(\omega))\in S(\G_m^d)$, so 
$s'(\omega)\in\trop_\phi\inv(\omega)\cap\psi\inv(S(\G_m^d)) =
\{s(\omega)\}$. Therefore $s'(\omega)=s(\omega)$, so since such $\omega$ are
dense in $Z'$, we conclude $s = s'$ on $Z'$.  Because $Z$ is contained in
the closure of its interior in $\Trop(U)$, the union of all such
$Z'$ is dense in $Z$, so since $s=s'$ on each $Z'$, we have $s=s'$ on $Z$.\qed

\begin{rem}
  Suppose that $\Trop(U)$ has multiplicity one everywhere.
  With the notation in the proof of Theorem~\ref{thm:section},
  we claim that $s(\Trop(U)) = \psi\inv(S(\G_m^d))$.
  For $\omega\in\Trop(U)\cap N_\Gamma$ we showed that
  $\{s(\omega)\} = \trop_\phi\inv(\omega)\cap\psi\inv(S(\G_m^d))$, which implies
  the claim if $\Gamma=\R$.  One easily reduces to this
  case by extending scalars to a non-archimedean extension field $L$ of $K$
  whose valuation map $L^\times\to\R$ is surjective, and using the fact that 
  $p\inv(S(\G_{m,K}^d)) = S(\G_{m,L}^d)$, where 
  $p:\G_{m,L}^{d,\an}\to\G_{m,K}^{d,\an}$ is the structural morphism.
  Therefore $s(\Trop(U))$ is a \emph{$c$-skeleton} in the sense
  of~\cite{ducros:image_reciproque, ducros:squelettes_modeles}.
  See Theorem~5.1 of~\cite{ducros:squelettes_modeles}.
\end{rem}

To conclude this section, we show that the image of the section of tropicalization
is contained in the skeleton in the case of a strictly semistable pair,
and that the section preserves integral affine structures in a suitable
sense.

\begin{prop}\label{prop:section-skeleton}
  Let $(\sX,H)$ be a strictly semistable pair of dimension $d$, let 
  $U = X\setminus\supp(H)_\eta$, and let $\varphi:U\to T\cong\G_m^n$ be 
  a closed immersion into an algebraic $K$-torus. Let $Z \subset \Trop(U)$
  be a subset such that $m_\Trop(\omega) = 1$ for all $\omega \in Z$. Then
  the image of the section $s: Z \rightarrow U^{\an}$ defined in 
  Theorem~\ref{thm:section} is contained in the skeleton
  $S(\sX,H)$. 
  
Moreover if $\Delta$ is an integral $\Gamma$-affine polyhedral face in $\Trop(U)$ of  dimension $d$ which is contained in $Z$, then 
$\Delta$ is covered by finitely many integral $\Gamma$-affine polyhedra $\Delta_i$ such that $s$ induces a unimodular integral $\Gamma$-affine 
map  $\Delta_i \rightarrow \Delta_{S_i}$ for a canonical polyhedron $\Delta_{S_i}$ of $S(\sX,H)$.
\end{prop}

\pf Define a partial ordering $\leq$ on $U^\an$ by declaring that
$x\leq y$ if $|f(x)|\leq|f(y)|$ for all $f\in K[M]$, where $M$ is the
character lattice of $T$.  This is indeed a partial ordering because
$U^\an\subset T^\an$ and $T^\an$ can be identified with a space of
seminorms on $K[M]$.
Let $z\in\Trop(U)$ be a point of tropical multiplicity one, let
$x = s(z)\in U^\an$, and let $y = \tau(x)\in S(\sX,H)$, where
$\tau:U^\an\to S(\sX,H)$ is the retraction map.  We want to show
that $x = y$.  Since $\trop_\phi$ factors through $\tau$ by Proposition~\ref{piecewise linear}, we have
$\trop_\phi\inv(z) = \tau\inv(\trop_\phi\inv(z)\cap S(\sX,H))$, so
$y\in\trop_\phi\inv(z)$ as well.  Since $x$ is by definition the Shilov
boundary point of $\trop_\phi\inv(z)$, we have $y\leq x$.

By the $\epsilon$-approximation argument used in the construction of the
skeleton of a strictly semistable pair in  
\S\ref{The skeleton of a strictly semistable pair}, 
there exists an affinoid neighbourhood $X'\subset U^\an$ of $y$ of the form
$X' = \tau\inv(X'\cap S(\sX,H))$ which is the generic fibre of a strictly
semistable formal scheme $\fX'$, such that the classical skeleton $S(\fX')$
coincides with $S(\sX,H)\cap X'$.  The retraction map $\tau:X'\to S(\fX')$
as defined by Berkovich also coincides with ours.  
Then~\cite[Theorem~5.2(ii)]{berkovich:locallycontractible1} gives
$x\leq y$, so $x=y$.

The unimodularity statement follows immediately from Proposition~\ref{piecewise linear} and
Theorem~\ref{Sturmfels--Tevelev multiplicity formula for alterations} .\qed

\appendix

\section{Refined intersection theory with Cartier divisors} \label{Refined intersection theory with Cartier divisors}

Let $K$ be an algebraically closed field  endowed with a non-trivial non-archimedean complete absolute value $|\scdot|$, corresponding valuation $v\coloneq -\log|\scdot|$, valuation ring $\kcirc$, residue field $\ktilde$ and value group $\Gamma \coloneq  v(K^\times)$. 
We will first recall the construction from \cite{gubler:local_heights} of the Weil divisor associated to a Cartier divisor on an admissible formal scheme over $\kcirc$. This will be useful in the 
paper for several local considerations. Then we will study the refined intersection product  of a Cartier divisor with a cycle on a proper flat variety $\sX$ over $\kcirc$. Note that we cannot use the algebraic intersection theory as in \cite{fulton:itheory} since the valuation ring $\kcirc$ is not noetherian. Instead we pass to the formal completion $\hat{\sX}$ of $\sX$ along the special fibre to get an admissible formal scheme. Then the refined intersection product is an easy consequence 
of the above construction of the associated Weil divisor. The reference for the refined intersection product is \cite[\S5]{gubler:lchs}.

We start with a Cartier divisor $D$ on a quasicompact admissible formal scheme $\fX$ over $\kcirc$. Our first goal is the construction of the Weil divisor $\cyc(D)$ on $\fX$ associated to $D$. 

\paragraph \label{cycles on analytic space}
Let $X$ be the generic fibre of $\fX$. We define cycles on $X$ as formal $\Z$-linear combinations of irreducible Zariski-closed  subsets of $X$. By definition, a Zariski-closed subset is 
the image of a closed immersion of analytic spaces over $K$. In rigid geometry, Zariski-closed subsets are called closed analytic subsets. 
Usually in our paper,
the generic fibre is algebraic and we have the basic operations for cycles as proper push-forward, flat pull-back and proper intersection with Cartier divisors from 
the first two chapters of Fulton's book \cite{fulton:itheory}. We note that this generalizes to quasicompact analytic spaces as they are covered by strictly affinoid subdomains $\sM(\Acal)$. Since $\Acal$ is a noetherian $K$-algebra, we may use the algebraic intersection theory on $\Spec(\Acal)$ and glue to get the corresponding operations on $X$ (see \cite[\S 2]{gubler:local_heights}). 
We should mention that $\Acal$ is not of finite type over $K$ as required in Fulton's book, but 
this assumption is not really necessary to develop the basic properties mentioned above (see~\cite{thorup:rational_equivalence}). Another issue in the analytic setting is the definition of irreducible components of $X$ which was handled in a paper by Conrad~\cite{conrad:irredcomps}.

\paragraph \label{cycle group} 
A {\it horizontal prime cycle} $\fY$ on 
$\fX$ is the closure of an irreducible   Zariski-closed analytic subset $Y$ of $X$
in $\fX$ as in \cite[Proposition 3.3]{gubler:local_heights}. 
That is, $\fY$ is the formal closed subscheme of $\fX$ defined by the ideal sheaf of regular formal functions on $\fX$ which vanish on $Y$.
A {\it horizontal cycle} on $\fX$ is a formal $\Z$-linear combination of horizontal prime cycles. 
A {\it vertical prime cycle} is an irreducible closed subset of  $\fX_s$. A {\it vertical cycle} is a formal $\Gamma$-linear combination of vertical prime cycles. 

A {\it cycle} $\fZ$ on $\fX$ is a formal sum of a horizontal cycle $\fY$ and a vertical cycle $V$. 
The prime cycles with non-zero coefficients are called the {\it prime components} of $\fZ$. 
A cycle is called {\it effective} if the multiplicities in its prime components are positive.

We say that a cycle $\fZ$ on $\fX$ is of  codimension $p$ if any horizontal prime component of  $\fZ$ is the closure of an irreducible closed analytic subset of $X$ of  codimension $p$ and if any vertical prime component of $\fZ$ has  codimension $p-1$ in $\fX_s$. 
A cycle on $\fX$ of codimension $1$ is called a  {\it Weil divisor}.

\begin{eg} \label{simplest example}
To understand why we need $\Gamma$-coefficients for vertical cycles, we look at the simplest example  $\fX=\Spec(\kcirc)$ and $D=\Div(f)$ for a non-zero $f \in K$. Then the valuation $v$ gives the multiplicity in the special fibre $\fX_s=\Spec(\ktilde)$ and we set $\cyc(D)\coloneq v(f)\fX_s$. 
\end{eg}
\medskip 

In general, the construction of the Weil divisor $\cyc(D)$ is based on the following local definition:

\paragraph \label{affinoid setting} 
Let $\fX=\Spf(\Acal^\circ)$ for a strictly affinoid algebra $\Acal$ and $D=\div(a/b)$ for $a,b \in\Acal^\circ$ which are not a zero-divisors. 
For an irreducible component $Y$ of $\fX_s$, there is a unique $\xi_Y$ in the generic fibre $X=\fX_\eta$ which reduces to the generic point of $Y$. Note that $X$ is the Berkovich spectrum of $\Acal$, therefore we get existence and uniqueness of $\xi_Y$ from \cite[Proposition 2.4.4]{berkovich:analytic_geometry}. Then we define the multiplicity  of $D$ in $Y$ as $\ord(D,Y) \coloneq \log|b(\xi_Y)|-\log|a(\xi_Y)|$.

\paragraph \label{associated Weil divisor} 
We deal now with the general case. We assume first that the generic fibre $X$ of $\fX$ is irreducible and reduced. This assumption is satisfied in all our applications.
To define the associated Weil divisor $$\cyc(D)=\sum_\fY \ord(D,\fY)\fY$$ on $\fX$, we have to define the multiplicity $\ord(D,\fY)$ of $D$ in a prime cycle $\fY$ of $\fX$ of codimension $1$. If $\fY$ is horizontal, then $\fY$ is the closure of an irreducible Zariski closed subset $Y$ of the generic fibre $X$ of codimension $1$. The restriction $D_\eta$ of $D$ to  $X$ is a Cartier divisor  and hence we may use \ref{cycles on analytic space} to define $\ord(D,\fY)\coloneq \ord(D_\eta,Y)$. 

If $\fY$ is vertical, then $\fY$ is equal to an irreducible component $V$ of $\fX_s$. We choose a formal affine open subset $\fU=\Spf(A)$ such that $\fU_s$ is a non-empty subset of $V$ and such that $D$ is given on $\fU$ by $a/b$ for $a,b \in A$ which are not a zero-divisors. Then $\Acal\coloneq  A \otimes_{\kcirc} K$ is a strictly affinoid algebra and we set $\fX'\coloneq \Spf(\Acal^\circ)$. We have a canonical morphism $\fX' \rightarrow \Spf(A)$ of admissible formal schemes with the same generic fibre which induces
a finite surjective morphism on special fibres (see 4.13 of \cite{gubler:guide} for the argument). Let $D'$ be the Cartier divisor on $\fX'$ given by the pull-back of $D$. Then $D'$ is given by $a/b$ on $\fX'$ as well. 
By using~\parref{affinoid setting}, we define the multiplicity of $D$ in $V$ by
$$\ord(D,V)\coloneq  \sum_{Y} [\ktilde(Y):\ktilde(V)] \ord(D',Y),$$
where $Y$ ranges over all irreducible components of $\fX_s'$ (see \cite[\S 3]{gubler:local_heights} for more details). 

If the generic fibre $X$ is not irreducible or not reduced, then we use the cycle $\cyc(X)$ associated to $X$ from \cite[2.7]{gubler:local_heights} and proceed by linearity in the prime components of $\cyc(X)$ to define $\cyc(D)$. 

\begin{rem} \label{meromorphic functions}
If $f$ is a meromorphic function on $\fX$ which is invertible as a meromorphic function, then $f$ defines a Cartier divisor $\div(f)$ on $\fX$. This notion has 
to be distinguished from the associated Weil divisor which we denote by $\cyc(f)$. We will use the same notation in the algebraic setting below.
\end{rem}

\begin{prop} \label{vertical Cartier divisors} 
  Let $D$ be a vertical Cartier divisor on the admissible formal scheme $\fX$
  over $\kcirc$; that is, $D$ is a Cartier divisor whose restriction to the generic fibre $X
  =\fX_\eta$ is trivial. We assume that the special fibre $\fX_s$ is
  reduced. Then the following properties hold:
\begin{itemize}
\item[(a)] The union of the prime components of $\cyc(D)$ is equal to $\supp(D)$.
\item[(b)] The Cartier divisor $D$ is effective if and only if $\cyc(D)$ is effective.
\item[(c)] We have $D = 0$ if and only if $\cyc(D) = 0$
\end{itemize}
In particular, the map $D \mapsto \cyc(D)$ is an injective homomorphism from the group of 
vertical Cartier divisors on $\fX$ to the group of (vertical) Weil divisors on $\fX$.
\end{prop}

\proof 
Recall that the support of the vertical Cartier divisor $D$ is the union of points of 
$\fX$ for which the restriction of $D$ to some neighbourhood is non-trivial. This gives 
a closed subset $\supp(D)$ of $\fX$. 
It follows from the definition of $\cyc(D)$ in \ref{associated Weil divisor} that every prime 
component of $\cyc(D)$ is contained in $\supp(D)$. 
Moreover, if $D$ is an effective Cartier divisor, then $\cyc(D)$ is an effective Weil divisor. 

First we prove (b). This claim is local and so we may assume that $\fX=\Spf(A)$ for a $\kcirc$-admissible algebra $A$ and that $D$ is given by $f=a/b$ for $a,b \in A$ which are not a zero-divisors. 
Now we use that $\fX_s$ is reduced. By a result of Bosch and L\"utkebohmert (see \cite[Proposition 1.11]{gubler:local_heights}), this implies that $A =\Acal^\circ$ for a strictly affinoid algebra $\Acal$ over $K$. 
To prove (b), it remains to show that $f \in \Acal^\circ$ if $\cyc(f)$ is an effective Weil divisor. Since $D$ is a vertical Cartier divisor, we know that $f$ is an invertible element of $\Acal$. As we assume now that $\cyc(f)$ is effective, we deduce from \ref{affinoid setting} that $|f(\xi_Y)| \leq 1$ for every irreducible 
component $Y$ of $\fX_s$. Since the supremum norm of $f$ on $X$ is equal to $\max_Y |f(\xi_Y)|$ \cite[Proposition 2.4.4]{berkovich:analytic_geometry}, we conclude that $f \in \Acal^\circ$, proving (b). 

Next we prove (a). Let $\supp(\cyc(D))$ be the union of the prime components of
$\cyc(D)$. We have seen at the beginning of the proof that $\supp(\cyc(D))
\subset \supp(D)$; we have to show equality. By passing to the open subset $\fX \setminus \supp(\cyc(D))$, we may assume that $\cyc(D)=0$. 
Then (b) implies that $D$ and $-D$ are both effective Cartier divisors, which means that $D$ is trivial. This proves (a).

Finally, (c) is an easy consequence of~(b).
\qed

\medskip

Now we switch to the algebraic setting. Our goal is to define a refined intersection theory of  Cartier divisors with cycles on a proper flat variety $\sX$ over $\kcirc$ with generic fibre $X=\Xcal_\eta$. Let $Z_k(\Xcal,\Gamma)=Z^p(\sX,\Gamma)$ be the group of cycles on $\sX$ of topological dimension $k$ (resp. of codimension $p$ with $p=\dim(\sX)-k$), where again the horizontal cycles have $\Z$-coefficients and the vertical cycles have coefficients in $\Gamma$.

\paragraph \label{proper intersection}  
We consider a Cartier divisor $D$ on $\sX$ and $\sZ \in Z_k(\sX,\Gamma)$ which
{\it intersect properly}.  This means that no prime component of $\sZ$ is
contained in the support $\supp(D)$ of the Cartier divisor $D$. In this
situation, the intersection product $D.Z$ is well-defined in
$Z_{k-1}(\sX,\Gamma)$ by the construction in~\parref{associated Weil
  divisor}. Indeed, proceeding by linearity in the prime components of $Z$, we
may assume that $\sZ$ is a prime cycle; then:

If $\sZ$ is horizontal, then $\sZ$ is the closure of a closed subvariety $Z$ of $X$.  By properness of the intersection, $D|_Z$ is a well-defined Cartier divisor on $Z$ and we define the 
horizontal part of $D.\sZ$ to be the cycle on $\sX$ induced from $\cyc(D|_Z)$ by passing to the closures of the prime components. 
Let $\fZ\coloneq \hat{\sZ}$ be the formal completion along the special fibre. Then $D$ induces 	
 a well-defined Cartier divisor $\hat{D}|_\fZ$ on the admissible formal scheme $\fZ$ over $\kcirc$. Since $\fZ$ and $\sZ$ have the same special fibre, it makes sense to define the vertical part of $D.\sZ$ as the vertical part of the Weil divisor associated to $\hat{D}|_\fZ$. We end up with a cycle $D.\sZ$ on $\sX$ with  support in $\supp(\sZ) \cap \supp(D)$. It follows easily from the construction that $\cyc(\hat{D}|_\fZ)$ is the formal completion of $D.\sZ$ defined componentwise. 

If $\sZ$ is vertical, then it is a closed subvariety of $\Xcal_s$ and we define $D.\sZ$ as the Weil divisor associated to $D|_\sZ$. This is a vertical cycle on $\sX$ with support in $\supp(\sZ) \cap \supp(D)$.

\begin{rem} \label{normal setting} 
In the applications, $\Xcal$ is often normal. For example, every $\kcirc$-toric variety is normal (see \cite[Proposition 6.11]{gubler:guide}) and hence it follows from \cite[Proposition 18.12.15]{egaIV_4} that every strictly semistable variety over $\kcirc$ is normal. If $\Xcal$ is normal, then the multiplicity $\ord(D,V)$ in the irreducible component $Y$ of $\Xcal_s$ has also an algebraic  description: It follows from results of Knaf that the local ring $\Ocal_{\Xcal,\zeta_V}$ in the generic point $\zeta_V$ is a valuation ring for a unique real valued valuation $w_V$ extending $v$. Then $\ord(D,V)= w_V(a)$ for any local equation $a$ of $D$ in $\zeta_V$. For details, we refer to \cite{gubler_soto:normal_toric_schemes}, Proposition 2.11. 
\end{rem}

\paragraph \label{rational equivalence}  
To get a refined intersection theory, we have to consider rational equivalence on a closed subset $S$ of $\Xcal$. Let  $R(S,\Gamma)$ be the subgroup of $Z(\Xcal,\Gamma)$ generated by all $\cyc(f|_\sY)$ and $\gamma\, \cyc(g|_V)$, where $\sY$ (resp.\ $V$) ranges over all horizontal (resp.\ vertical) closed subvarieties of $S$, where $f$ (resp.\ $g$) are non-zero rational functions of $\sY$ (resp.\ $V$) and where $\gamma$ ranges over the value group $\Gamma$. The {\it local Chow group of $\Xcal$ with support in $S$} is defined by
$$CH_S^*(\Xcal,\Gamma)\coloneq  Z^*(\Xcal,\Gamma)/R^*(S,\Gamma)$$
and it will be graded by codimension.

\begin{defn} \label{refined intersection product}
Let $D$ be a Cartier divisor on the proper flat scheme $\Xcal$ over $\kcirc$ and let $\sZ \in Z^p(\Xcal,\Gamma)$. For a closed subset $S$ of $\Xcal$ containing the support of $\sZ$, we define 
the {\it refined intersection product} $$D.\sZ \in CH_{\supp(D) \cap S}^{p+1}(\Xcal,\Gamma)$$ as follows: By linearity, we may assume that $\sZ$ is a prime cycle. If $D$ intersects $\sZ$ properly, then $D.\sZ$ is even well-defined as a cycle of codimension $1$ in $\sZ$ by~\parref{proper intersection}. If $\sZ$ is contained in $\supp(D)$, then we choose a linearly equivalent Cartier divisor $D'$ which intersects $\sZ$ properly and we define $D.\sZ$ as the class of $D'.\sZ$ in $CH_{\supp(D)\cap S}^{p+1}(\Xcal,\Gamma)$. 
\end{defn}
    
In the following result, we use proper push-forward and flat pull-back of cycles on flat varieties over $\kcirc$. The definitions are the same as in Fulton's book \cite{fulton:itheory}. 

\begin{prop} \label{properties of refined intersection}
The construction in Definition~\ref{refined intersection product} leads to a well-defined refined intersection product 
$$CH_S^{p}(\Xcal,\Gamma) \rightarrow CH_{\supp(D) \cap S}^{p+1}(\Xcal,\Gamma)$$
which maps the class of a cycle $\sZ$ with support in $S$ to $D.\sZ$. It has the following properties:
\begin{itemize}
 \item[(a)] The refined intersection product is bilinear using the union of supports.
 \item[(b)] If $\varphi:\Xcal' \rightarrow \Xcal$ is a morphism of flat proper varieties over $\kcirc$ and if $S'$ is a closed subset of 
$\Xcal'$ with $\varphi(S') \subset S$, then the projection formula 
$$\varphi_*(\varphi^*D.\alpha') = D. \varphi_*(\alpha') \in CH_{\supp(D) \cap S}^{p+1}(\Xcal,\Gamma)$$ 
holds for every $\alpha' \in CH_{S'}^{p}(\Xcal',\Gamma)$.
 \item[(c)] For Cartier divisors $D,E$ and $\alpha \in CH_S^{p}(\Xcal,\Gamma)$ on $\Xcal$, we have the commutativity law
$$D.E.\alpha= E.D. \alpha \in CH_{\supp(D) \cap \supp(E) \cap S}^{p+2}(\Xcal,\Gamma).$$
 \item[(d)] If $\varphi:\Xcal' \rightarrow \Xcal$ is a flat morphism of flat proper varieties over $\kcirc$ and if $\alpha \in CH_S^{p}(\Xcal,\Gamma)$, then 
$$\varphi^*(D.\alpha)=\varphi^*D.\varphi^*\alpha \in CH_{\varphi^{-1}(\supp(D) \cap S)}(\Xcal',\Gamma).$$
\end{itemize}
\end{prop}

\proof By passing to the formal completion of $\Xcal$ along the special fibre, this follows easily from \cite[Proposition 5.9]{gubler:lchs}. \qed

\begin{rem} \label{pseudo divisors}
The pull-back of the Cartier divisor $D$ with respect to the morphism $\varphi:\Xcal' \rightarrow \Xcal$ is only well-defined as a Cartier divisor if $\varphi(\Xcal')$
 is not contained in $\supp(D)$. However, the pull-back is well-defined as a pseudo-divisor in the sense of \cite[\S 2]{fulton:itheory} and the refined intersection product makes sense in this more general setting (see \cite[\S 5]{gubler:lchs} for details). We will not use it in our paper.
\end{rem}

\paragraph \label{local degree} 
We have a degree map on $0$-dimensional cycles of $Z(\Xcal,\Gamma)$. It is compatible with vertical rational equivalence, i.e.\ it induces a homomorphism 
$$\deg:CH_{\Xcal_s}^{d+1}(\Xcal,\Gamma) \rightarrow \Gamma$$
where $d\coloneq  \dim(X)=\dim(\Xcal)-1$. 
Let $D_0, \dots, D_k$ be Cartier divisors on $\Xcal$ and let $\Zcal \in Z_{k+1}(\Xcal, \Gamma)$ with horizontal part $Z$. We assume that
$\supp(D_0|_X) \cap \dots \cap \supp(D_k|_X)$ does not intersect the support of $Z$. 
Then we get a well-defined {\it intersection number}
$$D_0 \cdots D_k \cdot \Zcal \coloneq  \deg(D_0 \dots D_k.\Zcal) \in \Gamma.$$
Often we will consider the special case of Cartier divisors  $D_0, \dots, D_d$ on $\Xcal$ with 
$$\supp(D_0|_X) \cap \dots \cap \supp(D_d|_X) = \emptyset.$$ 
Using the cycle $\cyc(\Xcal)$ induced by the cycle  associated to $X$ by passing to the closure of the components, we get a well-defined {\it intersection number}
$$D_0 \cdots D_d \coloneq D_0 \cdots D_d \cdot \cyc(\Xcal)\coloneq  \deg(D_0 \dots D_d \cdot \cyc(\Xcal)) \in \Gamma.$$

\begin{rem} \label{refined intersection and formal model}
The refined intersection theory considered above works also for an admissible formal scheme $\fX$ over $\kcirc$ if the generic fibre is the analytification of a proper algebraic variety $X$. 
The same arguments apply (see \cite[\S 5]{gubler:lchs} for details). Thus it is
not necessary to check  if a given formal $\kcirc$-model of $X$ is algebraic. This will be used in the following example.
\end{rem}

\begin{eg} \label{product and strictly semistable}
Let $C$ be a smooth projective curve over $K$ with strictly semistable $\kcirc$-model $\Ccal$. 
For simplicity, we assume that every $1$-dimensional canonical simplex has the same length $v(\pi)$ (see \ref{canonical polyhedra}). 
Then $\Ycal\coloneq \Ccal \times \Ccal$ is a strictly polystable $\kcirc$-model of $X\coloneq C \times C$ such that 
the canonical polyhedra of the skeleton $S(\Ycal)$ are squares with edges of uniform length $v(\pi)$. We choose a diagonal in every square to get a triangulation of $S(\Ycal)$. The preimages of the triangles with respect to the retraction $\Xan \rightarrow S(\Ycal)$ form a formal analytic atlas of $\Xan$ inducing a strictly semistable formal scheme $\fX$ with skeleton $S(\fX)=S(\Ycal)$ as a set, but with canonical simplices given by the chosen triangulation (see \cite[Proposition 5.5, Remark 5.6, Remark 5.19]{gubler:canonical_measures} for this construction). 
The triangulation yields that $\fX$ lies over the formal completion of $\Ycal$. 
For $i=1,2$,  the projection $p_i:X=C \times C \rightarrow C$ to the $i$-th factor extends to a morphism $p_i:\fX \rightarrow \fC$ for the formal completion $\fC$ of $\Ccal$. Let $(p_i)_{\rm aff}:S(\fX) \rightarrow S(\Ccal)$ be the composition of $p_i:S(\fX) \rightarrow C^{\rm an}$ with the retraction $C^{\rm an} \rightarrow S(\Ccal)$. This is integral $\Gamma$-affine on every canonical simplex of $S(\fX)$ (see Proposition \ref{Functoriality}).  
By the stratum--face correspondence  in Proposition \ref{stratum-face correspondence}, we have a bijective correspondence between vertices $u$ of the canonical triangulation of $S(\fX)$ and irreducible components $Y_u$ of $\fX_s$. Using Proposition \ref{prop:comps.are.cartier}, there is a unique effective  vertical Cartier divisor $D_u$ with $\cyc(D_u)=v(\pi)\cdot V_u$. 
\end{eg}

Kolb has shown in \cite{kolb:calculus} that the intersection numbers of these vertical Cartier divisors $D_u$ can be computed by the following relations:

\begin{prop} \label{Kolb's relations}
Let $\fX$ be as in Example~\ref{product and strictly semistable} and assume that
every $1$-dimensional canonical simplex in $S(\Ccal)$ is determined
by its two vertices, i.e.\ that $S(\Ccal)$ is a graph with no multiple edges. For vertices $a,b$ of $S(\fX)$, we have the following two relations:
\begin{itemize} 
 \item[(a)] $D_a \cdot  D_b \cdot \sum_{c} D_c  = 0$, where $c$ runs over all vertices of $S(\fX)$.
 \item[(b)] If $(p_1)_{\rm aff}(a) \neq (p_1)_{\rm aff}(b)$ in $S(\Ccal)$, then 
$$D_a \cdot D_b \cdot \sum_{(p_1)_{\rm aff}(c) = (p_1)_{\rm aff}(b)} D_c = 0,$$
where $c$ runs over all vertices of $S(\fX)$ with $(p_1)_{\rm aff}(c) = (p_1)_{\rm aff}(b)$.
\end{itemize}
 \end{prop}

\proof As Kolb's Propositions 4.10 and 4.11 in  \cite{kolb:calculus} are formulated algebraically and over a discrete valuation ring, we reproduce the argument for convenience. The first relation is obvious from the 
fact that $\sum_c D_c = \Div(\pi)$. By construction, $a_1\coloneq (p_1)_{\rm aff}(a)$ and $b_1\coloneq (p_1)_{\rm aff}(b)$ are two distinguished vertices of $S(\Ccal)$. It is clear that $p_1(Y_a)$ (resp.\ $p_1(Y_b)$) is the irreducible component of $\Ccal_s=\fC_s$ corresponding to $a_1$ (resp.\ $b_1$). We may assume that $p_1(Y_a) \cap p_1(Y_b) \neq \emptyset$, as otherwise the intersection number in (b) is obviously zero. Hence we have an edge between $a_1$ and $b_1$. By assumption, this edge is completely determined by its vertices $a_1$ and $b_1$ which means geometrically that  $p_1(Y_a) \cap p_1(Y_b)=\{S\}$ for a single point $S \in \Ccal_s(\ktilde)$. Let $D_{b_1}$ be the unique effective vertical Cartier divisor on $\fC$ with $\cyc(D_{b_1})=v(\pi) \cdot Y_{b_1}$, where $Y_{b_1}$ is the irreducible component of $\Ccal_s$ corresponding to the vertex $b_1$ of $S(\Ccal)$ (see Proposition \ref{prop:comps.are.cartier}). In a neighbourhood of $S$, the Cartier divisor $D_{b_1}$ is given by a rational function $f$. We conclude that the Cartier divisor $D_{b_1}-\Div(f)$ has support in the complement of this neighbourhood of $S$ showing
$$0=D_a.D_b.p_1^*\Div(f )=D_a.D_b. p_1^*D_{b_1} \in CH_{\fX_s}^3(\fX,\Gamma).$$
It follows easily from the constructions that 
$$p_1^*(D_{b_1})= \sum_{(p_1)_{\rm aff}(c) = (p_1)_{\rm aff}(b)} D_c.$$
We conclude that
$$0=D_a.D_b. \sum_{(p_1)_{\rm aff}(c) = (p_1)_{\rm aff}(b)} D_c  \in CH_{\fX_s}^3(\fX,\Gamma),$$
proving the claim. \qed

\bibliographystyle{grw}
\bibliography{grw}

\providecommand{\noopsort}[1]{}\def\cprime{$'$}
\providecommand{\bysame}{\leavevmode\hbox to3em{\hrulefill}\thinspace}
\providecommand{\MR}{\relax\ifhmode\unskip\space\fi MR }
\providecommand{\MRhref}[2]{%
  \href{http://www.ams.org/mathscinet-getitem?mr=#1}{#2}
}
\providecommand{\href}[2]{#2}
\begin{thebibliography}{EGAIV${}_4$}

\bibitem[Ber90]{berkovich:analytic_geometry}
V.~G. Berkovich, \emph{Spectral theory and analytic geometry over
  non-{A}rchimedean fields}, Mathematical Surveys and Monographs, vol.~33,
  American Mathematical Society, Providence, RI, 1990.

\bibitem[Ber93]{berkovich:etalecohomology}
\bysame, \emph{{\'E}tale cohomology for non-{Archimedean} analytic spaces},
  Inst. Hautes {\'E}tudes Sci. Publ. Math. \textbf{78} (1993), 5--161.

\bibitem[Ber99]{berkovich:locallycontractible1}
\bysame, \emph{Smooth $p$-adic analytic spaces are locally contractible},
  Invent. Math. \textbf{137} (1999), no.~1, 1--84.

\bibitem[Ber04]{berkovich:locallycontractible2}
\bysame, \emph{Smooth $p$-adic analytic spaces are locally contractible.
  {II.}}, Geometric aspects of Dwork theory, Walter de Gruyter and Co. KG,
  Berlin, 2004, pp.~293--370.

\bibitem[BGR84]{bgr:nonarch}
S.~Bosch, U.~G{\"u}ntzer, and R.~Remmert, \emph{Non-{A}rchimedean analysis},
  Grundlehren der Mathematischen Wissenschaften [Fundamental Principles of
  Mathematical Sciences], vol. 261, Springer-Verlag, Berlin, 1984.

\bibitem[Bos77]{bosch:rigid_raume}
S.~Bosch, \emph{Zur {K}ohomologietheorie rigid analytischer {R}\"aume},
  Manuscripta Math. \textbf{20} (1977), no.~1, 1--27.

\bibitem[BPR11]{bpr:trop_curves}
M.~Baker, S.~Payne, and J.~Rabinoff, \emph{Nonarchimedean geometry,
  tropicalization, and metrics on curves}, 2011, preprint available at
  \url{http://arxiv.org/abs/1104.0320v2}; to appear in \textit{Algebraic
  Geometry} (Foundation Compositio Mathematica).

\bibitem[BPR13]{bpr:analytic_curves}
\bysame, \emph{On the structure of nonarchimedean analytic curves}, Tropical
  and {N}on-{A}rchimedean {G}eometry, Contemp. Math., vol. 605, Amer. Math.
  Soc., Providence, RI, 2013, pp.~93--121.

\bibitem[BR10]{baker_rumely:book}
M.~Baker and R.~Rumely, \emph{Potential theory and dynamics on the {B}erkovich
  projective line}, Mathematical Surveys and Monographs, vol. 159, American
  Mathematical Society, Providence, RI, 2010.

\bibitem[Car13]{cartwright:tropical_complexes}
D.~Cartwright, \emph{Tropical complexes}, 2013, Preprint available at
  \url{http://arxiv.org/abs/1308.3813}.

\bibitem[CD12]{chambert_loir_ducros:forms_currents}
A.~Chambert{--}Loir and A.~Ducros, \emph{Forms diff\'erentielles r\'eeles et
  courants sur les espaces de berkovich}, 2012, Preprint available at
  \url{http://arxiv.org/abs/1204.6277}.

\bibitem[Chr13]{christensen:thesis}
C.~Christensen, \emph{Erste {C}hernform und {C}hambert--{L}oir {M}a\ss e auf
  dem {Q}uadrat einer {T}ate-{K}urve.}, Ph.D. thesis, Universit\"at T\"ubingen,
  2013, available at \url{http://tobias-lib.uni-tuebingen.de}.

\bibitem[CHW14]{CHW}
M.~A. Cueto, M.~H{\"a}bich, and A.~Werner, \emph{Faithful tropicalization of
  the {G}rassmannian of planes}, Math. Ann. \textbf{360} (2014), no.~1-2,
  391--437.

\bibitem[Con99]{conrad:irredcomps}
B.~Conrad, \emph{Irreducible components of rigid spaces}, Ann. Inst. Fourier
  (Grenoble) \textbf{49} (1999), no.~2, 473--541.

\bibitem[Cue12]{cueto:implicitization}
M.~A. Cueto, \emph{Implicitization of surfaces via geometric tropicalization},
  2012, Preprint available at \url{http://arxiv.org/abs/1105.0509}.

\bibitem[dJ96]{dejong:alterations}
A.~J. de~Jong, \emph{Smoothness, semi-stability and alterations}, Inst. Hautes
  \'Etudes Sci. Publ. Math. (1996), no.~83, 51--93.

\bibitem[DP16]{dp}
J.~Draisma and E.~Postinghel, \emph{Faithful tropicalisation and torus
  actions}, Manuscripta Math. \textbf{149} (2016), no.~3-4, 315--338.
  \MR{3458171}

\bibitem[Duc03]{ducros:image_reciproque}
A.~Ducros, \emph{Image r\'eciproque du squelette par un morphisme entre espaces
  de {B}erkovich de m\^eme dimension}, Bull. Soc. Math. France \textbf{131}
  (2003), no.~4, 483--506.

\bibitem[Duc12]{ducros:squelettes_modeles}
\bysame, \emph{Espaces de {B}erkovich, polytopes, squelettes et th\'eorie des
  mod\`eles}, Confluentes Math. \textbf{4} (2012), no.~4, 1250007, 57.

\bibitem[EGAII]{egaII}
A.~{\noopsort{EGAII}}Grothendieck, \emph{\'{E}l\'ements de g\'eom\'etrie
  alg\'ebrique. {II}. \'{E}tude globale \'el\'ementaire de quelques classes de
  morphismes}, Inst. Hautes \'Etudes Sci. Publ. Math. (1961), no.~8, 222.

\bibitem[EGAIII${}_1$]{egaIII_1}
A.~{\noopsort{EGAIII_1}}Grothendieck, \emph{\'{E}l\'ements de g\'eom\'etrie
  alg\'ebrique. {III}. \'{E}tude cohomologique des faisceaux coh\'erents. {I}},
  Inst. Hautes \'Etudes Sci. Publ. Math. (1961), no.~11, 167.

\bibitem[EGAIV${}_4$]{egaIV_4}
A.~{\noopsort{EGAIV_4}}Grothendieck, \emph{\'{E}l\'ements de g\'eom\'etrie
  alg\'ebrique. {IV}. \'{E}tude locale des sch\'emas et des morphismes de
  sch\'emas {IV}}, Inst. Hautes \'Etudes Sci. Publ. Math. (1967), no.~32, 361.

\bibitem[Ful98]{fulton:itheory}
W.~Fulton, \emph{Intersection theory}, second ed., Springer-Verlag, Berlin,
  1998.

\bibitem[GS15]{gubler_soto:normal_toric_schemes}
W.~Gubler and A.~Soto, \emph{Classification of normal toric varieties over a
  valuation ring of rank one}, Documenta Math. (2015), no.~20, 171--198.

\bibitem[Gub98]{gubler:local_heights}
W.~Gubler, \emph{Local heights of subvarieties over non-{A}rchimedean fields},
  J. Reine Angew. Math. \textbf{498} (1998), 61--113.

\bibitem[Gub03]{gubler:lchs}
\bysame, \emph{Local and canonical heights of subvarieties}, Ann. Sc. Norm.
  Super. Pisa Cl. Sci. (5) \textbf{2} (2003), no.~4, 711--760.

\bibitem[Gub07]{gubler:tropical}
\bysame, \emph{Tropical varieties for non-{A}rchimedean analytic spaces},
  Invent. Math. \textbf{169} (2007), no.~2, 321--376.

\bibitem[Gub10]{gubler:canonical_measures}
\bysame, \emph{Non-{A}rchimedean canonical measures on abelian varieties},
  Compos. Math. \textbf{146} (2010), no.~3, 683--730.

\bibitem[Gub13a]{gubler:forms}
\bysame, \emph{Forms and currents on the analytification of an algebraic
  variety (after {C}hambert--{L}oir and {D}ucros)}, 2013, Preprint available at
  \url{http://arxiv.org/abs/1303.7364}.

\bibitem[Gub13b]{gubler:guide}
\bysame, \emph{A guide to tropicalizations}, Algebraic and combinatorial
  aspects of tropical geometry, Contemp. Math., vol. 589, Amer. Math. Soc.,
  Providence, RI, 2013, pp.~125--189.

\bibitem[Kol16]{kolb:calculus}
J.~Kolb, \emph{A simplicial calculus for local intersection numbers at
  non-archimedian places on products of semi-stable curves}, Abhandlungen aus
  dem Mathematischen Seminar der Universit{\"a}t Hamburg (2016), 1--36.

\bibitem[Mik05]{Mikhalkin:enumerative}
G.~Mikhalkin, \emph{{Enumerative tropical algebraic geometry in $\mathbb
  R^2$}}, J. Amer. Math. Soc. \textbf{18} (2005), no.~2, 313--377.

\bibitem[MS15]{maclagan_sturmfels:book}
D.~Maclagan and B.~Sturmfels, \emph{Introduction to tropical geometry},
  Graduate Studies in Mathematics, vol. 161, American Mathematical Society,
  Providence, RI, 2015.

\bibitem[Pal70]{palais:proper_maps}
R.~S. Palais, \emph{When proper maps are closed}, Proc. Amer. Math. Soc.
  \textbf{24} (1970), 835--836.

\bibitem[Pay09]{payne:analytification}
S.~Payne, \emph{Analytification is the limit of all tropicalizations}, Math.
  Res. Lett. \textbf{16} (2009), no.~3, 543--556.

\bibitem[Sil09]{silverman:I}
J.~H. Silverman, \emph{The arithmetic of elliptic curves}, second ed., Graduate
  Texts in Mathematics, vol. 106, Springer, Dordrecht, 2009.

\bibitem[ST08]{sturmfels_tevelev:elimination}
B.~Sturmfels and J.~Tevelev, \emph{Elimination theory for tropical varieties},
  Math. Res. Lett. \textbf{15} (2008), no.~3, 543--562.

\bibitem[Tat95]{tate:elliptic_functions}
J.~T. Tate, \emph{A review of non-{A}rchimedean elliptic functions}, Elliptic
  curves, modular forms, \& {F}ermat's last theorem ({H}ong {K}ong, 1993), Ser.
  Number Theory, I, Int. Press, Cambridge, MA, 1995, pp.~162--184.

\bibitem[Tho90]{thorup:rational_equivalence}
A.~Thorup, \emph{Rational equivalence theory on arbitrary {N}oetherian
  schemes}, Enumerative geometry ({S}itges, 1987), Lecture Notes in Math., vol.
  1436, Springer, Berlin, 1990, pp.~256--297.

\bibitem[Thu05]{thuillier:thesis}
A.~Thuillier, \emph{Th{\'e}orie du potentiel sur les courbes en
  g{\'e}om{\'e}trie analytique non archim{\'e}dienne. {A}pplications {\`a} la
  th{\'e}orie d'{A}rakelov}, Ph.D. thesis, University of Rennes, 2005, Preprint
  available at
  \url{http://tel.ccsd.cnrs.fr/documents/archives0/00/01/09/90/index.html}.

\bibitem[Tyo12]{tyomkin}
I.~Tyomkin, \emph{Tropical geometry and correspondence theorems via toric
  stacks}, Math. Ann. \textbf{353} (2012), no.~3, 945--995.

\bibitem[Ull95]{ullrich:direct_img}
P.~Ullrich, \emph{The direct image theorem in formal and rigid geometry}, Math.
  Ann. \textbf{301} (1995), no.~1, 69--104.

\end{thebibliography}

\bigskip
\hrule
\bigskip

\end{document}